\documentclass[11pt]{amsart}
\usepackage{amsmath,amssymb,vmargin,enumerate}

\title{Entropy boundary layers}
\author{Franck SUEUR}
\date{\today}

\newcommand{\eps}{\varepsilon}
\newcommand{\X}{\textbf{X}}

\renewcommand{\L}{\textbf{L}}
\renewcommand{\S}{\textbf{S}}

\providecommand{\Phi}{\textbf{\Phi}}

\newcommand{\un}{\underline}
\newcommand{\til}{\tilde}
\newcommand{\rond}{\mathring}
\newcommand{\M}{\mathcal{M}}
\renewcommand{\H}{\mathcal{H}}

\renewcommand{\P}{\mathcal{P}}

\providecommand{\R}{\mathbb{R}}
\providecommand{\N}{\mathbb{N}}
\providecommand{\Nrond}{\mathcal{N}}

\newtheorem{theo}{Theorem}[section]
\newtheorem{coro}[theo]{Corollary}
\newtheorem{lem}[theo]{Lemma}

\newtheorem{prop}[theo]{Proposition}
\newtheorem{rem}[theo]{Remark}

\newtheorem{step}{Step}
\newtheorem{step2}{Step}

\begin{document}

\begin{abstract}

We consider the Euler system of compressible and entropic gaz dynamics in a bounded open domain of $\R^d$ with 
 wall boundary condition. 
We prove the existence and the stability of families of solutions which correspond to a ground state 
plus a large entropy boundary layer.
The ground state is a solution of the Euler system which satisfies some explicit additional conditions on the boundary.
These conditions are used in a reduction of the system.
We construct BKW expansions at all order. 
The profile problems are linear thanks to a transparency property.
We prove the stability of these expansions by proving $\eps$-conormal estimates
 for a characteristic boundary value problem.
 
\end{abstract}

\maketitle

\section{Introduction}

  We consider the  \textbf{Euler system of compressible and entropic gaz dynamics}. 
We respectively denote by  $\mathtt{s}$, $\mathtt{v} $  and  $\mathtt{p}$ the  entropy, the speed and the pressure. 
 The volumic density, noted $\mathtt{\rho}$ , is a function of  $\mathtt{p}$ and $\mathtt{s}$.  
 We assume that  $\mathtt{\rho}(\mathtt{p},\mathtt{s}) > 0$ for all $\mathtt{p}$ and  $\mathtt{s}$. 
We also introduce the function $\alpha (\mathtt{p},\mathtt{s}) := 
(\mathtt{\rho} (\mathtt{p} ,\mathtt{s}))^{-1} \,  \mathtt{\rho}'_\mathtt{p} (\mathtt{p},\mathtt{s}) $.
  We assume that $\alpha (\mathtt{p},\mathtt{s}) > 0$ for all $\mathtt{p}$ and  $\mathtt{s}$.
  We denote by $t$ the time variable and by $x=(x_1 , ...,x_d )$ the space variable. 
  We denote by $\nabla$ the gradient with respect to $x$.
   The Euler system is:
 \begin{eqnarray*}
  \X_{\mathtt{v}}  \,  \mathtt{v} + \mathtt{\rho}^{-1} \,  \nabla  \mathtt{p} =0 ,
\quad \X_{\mathtt{v}} \,  \mathtt{p} +  \alpha^{-1} \, \mathrm{div} \  \mathtt{v}  =0 \quad \X_{\mathtt{v}} \,  \mathtt{s}=0 ,
\end{eqnarray*}
  where  $\X_{\mathtt{v}} \, $ is the particle derivative $\X_{\mathtt{v}} \,  := \partial_t + \mathtt{v} .\nabla $.
The previous system  is a nonconservative form of a system of conservation laws. 
Moreover, this system is  hyperbolic symmetrizable. It has three characteristic fields.
One of them is linearly degenerate (cf. \textit{section $\ref{setting}$}). 
 We consider these equations in a bounded open domain  $\Omega \subset \R^d$ lying on one side
of its $C^\infty$  boundary $\Gamma$. 
More precisely, since we will need an equation of
the boundary $\Gamma $, we fix once for all a function
$\varphi \in C^{\infty}(\R^d, \R)$ and we assume that
$\Omega = \{ \varphi > 0 \}$, $ \Gamma = \{ \varphi = 0 \}$ 
and  $| \nabla \varphi (x) | = 1 $ in an open neighborhood $\mathcal{V}_\Gamma $ 
of $\Gamma $\footnote{ Hence for $x \in \Omega \cap \mathcal{V}_\Gamma$: 
$\varphi(x) = dist(x,\Gamma )$.}. 
We consider the natural boundary condition $ \mathtt{v}. \mathfrak{n} =0$, where $\mathfrak{n}$ is
the unit outward normal to $\Gamma$. 
For $T>0$, the \textbf{boundary value problem} reads:
\begin{eqnarray}
\label{aller}
\left\{
\begin{array}{cc}
\left.
\begin{array}{l}
  \X_{\mathtt{v}} \,  \mathtt{v} + \mathtt{\rho}^{-1} \,  \nabla  \mathtt{p} =0 
\\  \X_{\mathtt{v}} \,  \mathtt{p} +  \alpha^{-1} \, \mathrm{div} \ \mathtt{v}  =0
\\ \X_{\mathtt{v}} \,  \mathtt{s}=0
\end{array}
\right\}
\quad &\mathrm{where} \, (t,x) \in (0,T) \times \Omega
\\  \mathtt{v}. \mathfrak{n} = 0 \quad &\mathrm{where} \, (t,x) \in (0,T) \times \Gamma
\end{array}
\right.
\end{eqnarray}
The boundary $\Gamma$ is  characteristic for the linearly degenerate field. 
We introduce the tangential velocity $\mathtt{v}_\mathfrak{t} $ and 
the  normal velocity $ \mathtt{v}_\mathfrak{n}$. 
Thus, we have $ \mathtt{v} = \mathtt{v}_\mathfrak{t} + \mathtt{v}_\mathfrak{n} $. 
The choice of the set of thermodynamic variables $\mathtt{v}$, $\mathtt{p}$ and  $\mathtt{s}$ is particularly well adapted 
to boundary problem (cf. \cite{moi3}).
The existence of \textbf{local regular solutions} of $(\ref{aller})$ is given in \cite{cho} and \cite{Gu90}.
If $\mathcal{O}$ is an open subset of $\R^d$, we denote by $H^\infty (\mathcal{O})$ the set of $u \in L^2 (\mathcal{O})$ 
such that all the derivatives of $u$ are in $L^2 (\mathcal{O})$. 
From now on, we will assume that a real $T_0 > 0$ and 
a solution $u^0 :=(\mathtt{v}^0 , \mathtt{p}^0 ,\mathtt{s}^0 ) \in H^\infty ((0,T_0) \times \Omega)$ of $(\ref{aller})$ are given.

An interesting question about the Euler system is the study of the convergence of a more acute model: 
the Navier-Stokes system, which includes a viscosity term,  when the amplitude of the viscosity goes to $0$.  
The difficulty is linked to the existence of a \textbf{boundary layer} i.e. of a rapid variation 
of the solutions of the viscous model
 near the boundary. 
There is a sensitivity on the type of boundary conditions imposed for the viscous model. 
The most delicate is the homogeneous Dirichlet condition. 
In this case, there are in general some large characteristic boundary layers
 of large amplitude\footnote{For other boundary conditions the amplitude of the boundary layers 
 can be weaker or the boundary layers can be noncharacteristic.
Then the problem is simpler (see \cite{fg}, \cite{cmr}, \cite{moi3}, \cite{Temam}).}.
 In one space  dimension, a simple case is the isentropic one since there is no boundary layer 
 and the solutions of the Navier-Stokes
system are regular perturbations of the solutions of the Euler system. 
For the entropic Navier-Stokes  equations, 
an answer was given by F.Rousset in
\cite{Rousset} using boundary layers analysis. 
In several space dimensions, the analysis is quite more complicated even for the incompressible Navier-Stokes equations.
There is a huge literature about the tangential velocity boundary layers which appear (see, for example, 
the papers of Z.Xin and T.Yanagisawa \cite{Xin}, M.Sammartino and R.E.Caflisch
\cite{Samm1}, \cite{Samm2}, E.Grenier \cite{Gr1}, \cite{Gr2},...).
An attempt of analysis involves Prandtl equations
(see the surveys of W.E \cite{E} and E.Grenier \cite{Grenier}). 

Here we do not consider the  Navier-Stokes equations but only the Euler system.
The idea to investigate first the stability of boundary layers type solutions for the Euler equations is a classical approach in
fluid mechanics (see the books of P.G.Dranzin and W.H.Reid \cite{Drazin}, of C.Marchioro and M.Pulvirenti \cite{Marchioro},
S.Schochet \cite{Sch}. 
This idea was followed  more recently by E.Grenier for the study of velocity boundary layers  \cite{Gr1}, \cite{Gr2}.
Such a strategy is also possible for entropy boundary layers since they are characteristic, like the velocity ones.  
Thus our goal in this paper is to study  \textbf{entropy boundary layers} for the Euler system.
As far as we know there was no mathematical study of entropy boundary layer in \textbf{several space dimensions}. 

\section{Overview of the results}

We introduce the space $\Nrond (T) :=  H^\infty  ( (0,T) \times  \Omega , \mathcal{S}(\R_+)$ 
where we  denote by $\mathcal{S} (\R_+ )$ the Schwartz space of $C^\infty$ rapidly decreasing functions. 
 Thus a function $U (t,x,X) \in  \Nrond (T) $ is $C^\infty$ rapidly decreasing with respect to $X$.
 Let us begin to look naively for solutions $ u^\eps :=( \mathtt{v}^0 , \mathtt{p}^0 ,\mathtt{s}^\eps )_{0 < \eps \leqslant 1}$
 of $(\ref{aller})$ with $\mathtt{s}^\eps$ of the form 
  \begin{eqnarray}
 \label{elec} 
 \mathtt{s}^\eps (t,x) :=  \mathtt{s}^0  (t,x) + \tilde{\mathtt{S}}^0 (t,x,\frac{\varphi(x)}{\eps})
 \end{eqnarray}
 where the function $\tilde{\mathtt{S}}^0 $ is in $\Nrond (T_0)$ and
the function $u^0 := (\mathtt{v}^0 , \mathtt{p}^0 ,\mathtt{s}^0 ) $ is the ground state  given above.
Replacing in the third equation of the system $(\ref{aller})$ leads to the following linear transport equation
 that  necessarily the function $\tilde{\mathtt{S}}^0 $ verifies:
 \begin{eqnarray}
 \label{cle} 
  (\partial_t  + \mathtt{v}^0 .\nabla
 +  \frac{\mathtt{v}^0_{\mathfrak{n}}}{\varphi(x)}  .X \partial_X ) \tilde{\mathtt{S}}^0 =
0 \quad \mathrm{where} \, (t,x,X) \in (0,T_0) \times \Omega  \times \R_+ .
 \end{eqnarray}
 Thanks to the boundary condition (see $(\ref{aller})$), the function $\mathtt{v}^0_{\mathfrak{n}} \, /\varphi(x)$ is 
 $C^\infty$\footnote{Notice that for this transport equation, the boundary $\Gamma \times \{0\}$ of the domain $\Omega \times \R_+ $ is totally
 characteristic. As a consequence,  none boundary condition needs to be prescribed for $\tilde{\mathtt{S}}^0 $.}.
In general these necessary conditions are not sufficient to insure 
that the functions $( \mathtt{v}^0 , \mathtt{p}^0 ,\mathtt{s}^\eps )_{0 <\eps \leqslant 1}$ are solutions of $(\ref{aller})$.
Indeed, because the functions $\mathtt{\rho}$ and $\alpha$ depend on $\mathtt{p}$ and $\mathtt{s}$, the two first equations of
$(\ref{aller})$ are not satisfied.
However if in addition we assume that the ground state $(\mathtt{v}^0 , \mathtt{p}^0 ,\mathtt{s}^0 ) $ satisfies 
\begin{eqnarray}
\label{pro}
\X_{\mathtt{v}^0} \,  \,  \mathtt{v}^0 =0 , \quad \mathrm{div}_x \ \mathtt{v}^0 =0 , \quad \nabla_{t,x} \, \mathtt{p}^0 =0 \quad \mathrm{where} \,
(t,x) \in (0,T_0) \times \Omega,
\end{eqnarray}
then it is easy to check that the functions $(\mathtt{v}^0 , \mathtt{p}^0 ,\mathtt{s}^\eps )_{0 <\eps \leqslant 1}$ 
are solutions of $(\ref{aller})$.
 \textit{In section $\ref{ground}$}, we will prove the existence of ground states 
 $(\mathtt{v}^0 , \mathtt{p}^0 ,\mathtt{s}^0 ) \in H^\infty$ solutions of
 $(\ref{aller})$ and verifying the conditions $(\ref{pro})$.

 Our goal is to relax the conditions $(\ref{pro})$ into   \textbf{conditions localized on the boundary $\Gamma$}.
  In fact the conditions $(\ref{pro})$ were first introduced in a paper of C.Cheverry, O.Gu\`es and G.M\'etivier  \cite{cgm2}
   where the existence and the stability of large amplitude high frequency
  entropy waves are shown.  
  Here we look for entropy  boundary layers which are local singularities. 
  It seems rather natural that local conditions are sufficient to deal with  boundary layers.
  We reach our goal and claim the following theorem.
  
  \begin{theo}
  \label{int}
  
  Assume that a solution  $u^0 :=(\mathtt{v}^0 , \mathtt{p}^0 ,\mathtt{s}^0 ) \in H^\infty ((0,T_0) \times \Omega)$ of
  $(\ref{aller})$ verifying  
  \begin{eqnarray}
  \label{cond}
\X_{\mathtt{v}^0} \,  \mathtt{v}^0 =0 \ \mathrm{and} \ \X_{\mathtt{v}^0} \,  \mathtt{p}^0 =0
 \quad  \mathrm{for} \ (t,x) \in (0,T_0) \times \Gamma
  \end{eqnarray}
   and a solution $ \tilde{\mathtt{S}}^0 \in \Nrond (T_0)$ of  $(\ref{cle})$ are given.
  For $0 < \eps \leqslant 1$, we denote by  $ u^\eps $ the function 
  $ u^\eps :=( \mathtt{v}^0 , \mathtt{p}^0 ,\mathtt{s}^\eps )$ where
  $\mathtt{s}^\eps$ is given by $(\ref{elec})$.
  Then there exists a family of solutions  $(\tilde{u}^\eps )_{0 < \eps \leqslant 1}$ in $H^\infty ((0,T_0) \times \Omega) $
  of $(\ref{aller})$ such that  $\tilde{u}^\eps - u^\eps$ tends to $0$ in $H^1 ((0,T_0) \times \Omega)$ when $\eps \rightarrow 0^+ $.
  
   \end{theo}
  
  In fact Theorem $\ref{int}$ is a corollary of more acute results involving \textbf{WKB} (Wentzel-Kramers-Brillouin) expansions.
 We will prove the \textbf{existence and stability} of families of solutions 
  $( \mathtt{v}^\eps , \mathtt{p}^\eps ,\mathtt{s}^\eps)_{0 < \eps \leqslant 1}$
of the Euler system with a large amplitude
entropy boundary layer i.e.  of the form 
\begin{eqnarray}
\label{pre}  \quad
   \left\{
   \begin{array}{l}
  \mathtt{v}_\mathfrak{t}^\eps  (t,x) := \mathtt{v}_\mathfrak{t}^0  (t,x) + \eps
(\mathtt{V}_\mathfrak{t}  (t,x)  
+ \tilde{\mathtt{V}}_\mathfrak{t} (t,x,\frac{\varphi(x)}{\eps})) + O( \eps^2 ) ,
\vspace{0.2cm} \\  \mathtt{v}_\mathfrak{n}^\eps (t,x) := \mathtt{v}_\mathfrak{n}^0  (t,x)+ \eps \mathtt{V}_\mathfrak{n} (t,x)  
+ O( \eps^2 ) ,
\vspace{0.2cm} \\  \mathtt{p}^\eps  (t,x):= \mathtt{p}^0 (t,x) + \eps \mathtt{P} (t,x)  + O( \eps^2 ) ,
\vspace{0.2cm} \\  \mathtt{s}^\eps (t,x) := \mathtt{s}^0 (t,x)  +  \tilde{\mathtt{S}}^0 (t,x,\frac{\varphi(x)}{\eps})  + O(\eps ) ,
\end{array}
  \right.
\end{eqnarray}
where $u^0 :=(\mathtt{v}^0 , \mathtt{p}^0 ,\mathtt{s}^0 )$ is a solution of $(\ref{aller})$  verifying the conditions $(\ref{cond})$.
 This analysis is inspired by \cite{cgm2}, where the propagation of large amplitude high frequency entropy waves is shown 
 for ground states $u^0 $ solutions of $(\ref{aller})$  
 verifying the conditions $(\ref{pro})$
 (in the first equation of section $2.3$ of \cite{cgm2}, read $ \nabla_{t,x} \mathtt{p} =0$ instead of $\nabla_{x} \mathtt{p}=0$).
 In our analysis the condition on the particle derivative is localized on the boundary. 
 Considering directly the Euler system, we use implicitly the  general structure conditions of \cite{cgm2}.
 In particular, the choice of the set of thermodynamic variables $\mathtt{v}$, $\mathtt{p}$ and  $\mathtt{s}$ is a key point. 

 An important feature of the expansion $(\ref{pre})$ is that there are some boundary layers not only on the entropy 
 but also on the
other components, ponderated by some $\eps$.
More accurately, boundary layer appears on tangential velocity with an amplitude $\eps$ whereas boundary layer
appears on the normal velocity and the pressure with an amplitude $\eps^2$. 
 The conditions $(\ref{cond})$ on the ground state play a main role in the fact that the large boundary layer keeps polarized on
 the entropy. 
 \textit{In section $\ref{redbull}$}, we reduce the system, thanks to  a change of unknown singular with respect to
 $\eps$. 
 This \textbf{reduction} is inspired by \cite{cgm2}. 
 The conditions $(\ref{cond})$ on the ground state will be used at this step.
 Because the localized conditions $(\ref{cond})$ are weaker than the conditions $(\ref{pro})$ of \cite{cgm2}, our reduction is much
 more delicate than in \cite{cgm2}. We now describe
in more detail the rest of the contents of the paper. 
 \\
 \\ \textit{In section $\ref{secform}$}, we look for \textbf{formal WKB solutions} of the problem  $(\ref{aller})$. 
This means that we construct WKB expansions of infinite order.
Let us precise this.
We introduce the profile space  
\begin{eqnarray}
\label{defp}  \P (T):= \{ \mathcal{U} \in L^2 ((0,T) \times \Omega \times \R_+  )  \ \mathrm{such} \ \mathrm{that} \  \mathrm{there} \ 
\mathrm{exist} \ \un{ \mathcal{U}} \in H^\infty ((0,T) \times \Omega ) \ \mathrm{and}
\\  \nonumber \ \tilde{ \mathcal{U}} 
\in  \Nrond (T) \   \mathrm{such} \ \mathrm{that} \ \mathrm{for} \ \mathrm{all} \ (t,x,X)
 \in (0,T) \times \Omega \times \R_+ ,
  \  \mathcal{U}  (t,x, X)= \un{ \mathcal{U}} (t,x) 
+ \tilde{ \mathcal{U}} (t,x,  X)   \}. 
\end{eqnarray}
  The function  $\un{ \mathcal{U}}$ is the regular part and  $\tilde{ \mathcal{U}}$ is
   the characteristic boundary layer term. 
   We will  split $\mathcal{U} \in \P(T)$ into $\mathcal{U} =(\mathtt{V},\mathtt{P} ,\mathtt{S})$ and $\mathtt{V}$ 
 into $\mathtt{V} := \mathtt{V}_\mathfrak{t} + \mathtt{V}_\mathfrak{n}$ where $\mathtt{V}_\mathfrak{n} :=  (\mathtt{V}  .
 \mathfrak{n})  \mathfrak{n}$.
 The function  $\mathtt{V}$ (respectively  $\mathtt{P}$) takes its values in $\R^{d}$ (resp. $\R$).
 The function  $\mathtt{V}_\mathfrak{t} $ (respectively  $\mathtt{V}_\mathfrak{n}$) takes its values in $\R^{d-1}$ (resp. $\R$).
By abuse of notations, we will say that $\mathtt{V}$, $\mathtt{V}_\mathfrak{t} $,  $\mathtt{V}_\mathfrak{n} $ 
and $\mathtt{P}$ are in $\P(T)$ even if they do not take values in $\R^{d+2}$. 

   We look for formal solutions  $(u^\eps )_\eps $ of $(\ref{aller})$ of the form
\begin{eqnarray}
\label{sep}
 u^\eps (t,x) = \sum_{n \geqslant 0} \eps^n \, \mathcal{U}^n (t,x, \frac{\varphi(x) }{\eps})
\end{eqnarray}
 where each  $\mathcal{U}^n$ belongs to  $\P (T) $ and $\un{\mathcal{U}}^0 = u^0 $. 
 Let us explain what is meant by \textit{formal} \textit{solutions}.
 Plugging the expansion $(\ref{sep})$ into the system, using Taylor expansions and ordering the terms in powers of $\eps$, 
  we get a formal expansion in power series of $\eps$:
\begin{eqnarray*}
\sum_{n \geqslant -1} \eps^n \,  \Phi^n (t,x, \frac{\varphi(x) }{\eps})
\end{eqnarray*}
 where the  $(\Phi^n )_{n \geqslant -1}$ are in  $\P (T) $.
 We say that $(u^\eps )_\eps $ is a formal solution when all the resulting $\Phi^n $ are identically zero.
 Theorem $\ref{try}$ will sum up the main results of  section $\ref{secform}$.
 It  states that the system has formal solutions of the form
 \begin{eqnarray}
 \label{pre2}  \quad
   \left\{
   \begin{array}{l}
  \mathtt{v}_\mathfrak{t}^\eps (t,x) = \mathtt{v}_\mathfrak{t}^0 (t,x) 
  + \displaystyle{\sum_{j \geqslant 1}} \eps^j \,  \mathtt{V}_\mathfrak{t}^{j-1} (t,x, \frac{\varphi(x) }{\eps}) ,
\vspace{0.2cm}  \\ \mathtt{v}_\mathfrak{n}^\eps (t,x) =  \mathtt{v}_\mathfrak{n}^0 (t,x) 
  + \eps \,  \mathtt{V}_\mathfrak{n}^0 (t,x) 
  + \displaystyle{\sum_{j \geqslant 2}} \eps^j \,  \mathtt{V}_\mathfrak{n}^{j-1} (t,x, \frac{\varphi(x)}{\eps}) ,
\vspace{0.2cm}  \\ \mathtt{p}^\eps (t,x) =  \mathtt{p}^0 (t,x) + \eps \, \mathtt{P}^0 (t,x) 
  +  \displaystyle{\sum_{j \geqslant 2}} \eps^j \,  \mathtt{P}^{j-1} (t,x, \frac{\varphi(x)}{\eps}) ,
  \vspace{0.2cm}  \\ \mathtt{s}^\eps (t,x) = \displaystyle{\sum_{j \geqslant 0}} \eps^j \,  \mathtt{S}^{j} (t,x,
  \frac{\varphi(x)}{\eps}) ,
    \end{array}
\right.
    \end{eqnarray}
  and that we can prescribe arbitrary initial values to the $(\mathtt{S}^j |_{t=0})_{j \in \N } $  
  and to the $(\mathtt{V}_\mathfrak{t}^j |_{t=0})_{j \in \N } $. 
   The profiles $\mathtt{V}_\mathfrak{t}^0$, $\mathtt{V}_\mathfrak{n}^0$ and $\mathtt{P}^0$ involved in $(\ref{pre2})$ are the
  profiles  $\mathtt{V}_\mathfrak{t}$, $\mathtt{V}_\mathfrak{n}$ and $\mathtt{P}$ involved in $(\ref{pre})$. 
  In $(\ref{pre2})$, we use the index $0$ to match with notations of section  $\ref{secform}$.
  In $(\ref{pre})$, we did not write the index in order to avoid heavy notations. 
   
  The two first equations of $(\ref{aller})$ involve the entropy $\mathtt{s}$, through the functions $\rho$ et $\alpha$.
  A priori the large entropy boundary layer could contaminate the 
  velocity\footnote{which can create a destabilizing effect  (see \cite{Gr1}, \cite{Gr2})}.
  We prove that since we consider some ground states  $u^0 := (\mathtt{v}^0 , \mathtt{p}^0 ,\mathtt{s}^0 )$ solutions of 
 $(\ref{aller})$  verifying the conditions $(\ref{cond})$, there is no contamination at
 order $\eps^0$: when looking at the expansion $(\ref{pre2})$, we notice that there is no large pressure boundary layer and
 no large  velocity boundary layer.
 To use  the conditions $(\ref{cond})$, we reduce the system (cf. section $\ref{redbull}$). 
 This step is difficult.
 Let us briefly mention here some key points of our strategy.
 For sake of clarity, we  begin with the first equation of $(\ref{aller})$ only.
 We look for solutions $(u^\eps )_\eps $ 
   where $\mathtt{v}^\eps$ and $\mathtt{p}^\eps$ are of the form
     $\mathtt{v}^\eps:= \mathtt{v}^0 + \eps \mathtt{V}^\eps $, 
     $\mathtt{p}^\eps := \mathtt{p}^0 + \eps \mathtt{P}^\eps $.
 We split the particle derivative $\X_{\mathtt{v}^\eps} \, \mathtt{v}^\eps$ into
 \begin{eqnarray*}
   \X_{\mathtt{v}^\eps} \, \mathtt{v}^\eps =  \X_{\mathtt{v}^0} \, \mathtt{v}^0 
   + \eps  \X_{\mathtt{v}^\eps} \, \mathtt{V}^\eps
   +  \eps \mathtt{V}^\eps . \nabla  \mathtt{v}^0 .
   \end{eqnarray*}
%%%
  Thanks to the conditions $(\ref{cond})$, there exists a function $\phi \in C^\infty ((0,T_0 ) \times \Omega )$ 
  such that for all  $(t,x) \in (0,T_0 ) \times \Omega $, 
   $(\X_{\mathtt{v}^0} \, \mathtt{v}^0 )(t,x)= \varphi(x)\, \phi(t,x) $.
   Remark that if $\tilde{ \mathcal{U}} \in   \Nrond (T_0 ) $, then 
 \begin{eqnarray*}
  \X_{\mathtt{v}^0} \, \mathtt{v}^0 \, . \tilde{\mathcal{U}} (t,x, \frac{\varphi(x)}{\eps}) =
  \eps (\phi(t,x) X \tilde{ \mathcal{U}} (t,x,X))|_{X=\frac{\varphi(x)}{\eps}} 
  \quad \mathrm{and} \quad X \tilde{ \mathcal{U}} \in  \Nrond (T_0 ).
  \end{eqnarray*}
%%%
To use this remark, we develop the term $\rho(\mathtt{p}^\eps , \mathtt{s}^\eps )$. 
Thus, we write  $\mathtt{s}^\eps$ under the form  
 \begin{eqnarray}
 \label{SSS}
 \mathtt{s}^\eps (t,x) = \mathtt{s}^0 (t,x) + \tilde{\mathtt{S}}^0  (t,x ,\frac{\varphi(x)}{\eps})  
   + \eps \mathtt{S}^{\eps,\flat} (t,x),
 \end{eqnarray}
with $\tilde{\mathtt{S}}^0  \in  \Nrond (T_0 ) $.
We get for the term $\rho( \mathtt{p}^\eps , \mathtt{s}^\eps )$ the following expansion:
 \begin{eqnarray*}
\rho( \mathtt{p}^\eps , \mathtt{s}^\eps ) &=&  \rho( \mathtt{p}^0 , \mathtt{s}^0 )
+  (\rho( \mathtt{p}^0 ,  \mathtt{s}^0 + \tilde{\mathtt{S}}^0 ) -  \rho( \mathtt{p}^0 , \mathtt{s}^0 ))
+  \eps D \rho( \mathtt{p}^0 , \mathtt{s}^0 ). ( \mathtt{P}^\eps , \mathtt{S}^{\eps,\flat} )
+ O(\eps^2 ) .
\end{eqnarray*}
%%%
Because  $u^0$ satisfies the equation  $(\ref{aller})$, we get the equation:
 \begin{eqnarray*}
 \rho( \mathtt{p}^\eps , \mathtt{s}^\eps ) .\X_{\mathtt{v}^\eps } \,  \mathtt{V}^\eps + \nabla   \mathtt{P}^\eps   
   + \mathfrak{P}_\mathtt{v} = O(\eps ) \quad
   \mathrm{with}  \quad
  \mathfrak{P}_\mathtt{v} := \mathfrak{P}_\mathtt{v}^a + \mathfrak{P}^b_\mathtt{v} + \mathfrak{P}^c_\mathtt{v}
 \\ \nonumber  
  \mathrm{where}  \quad 
  \left\{
\begin{array}{l}
 \mathfrak{P}_\mathtt{v}^a :=
  \rho( \mathtt{p}^0 , \mathtt{s}^\eps ). \mathtt{V}^\eps . \nabla \mathtt{v}^0 , 
\vspace{0.2cm}\\  \mathfrak{P}^b_\mathtt{v} := 
\eps^{-1} \, (\rho( \mathtt{p}^0 ,  \mathtt{s}^0 + \tilde{\mathtt{S}}^0 ) -  \rho( \mathtt{p}^0 , \mathtt{s}^0
)).\X_{\mathtt{v}^0} \, \mathtt{v}^0 , 
\vspace{0.2cm}\\  \mathfrak{P}^c_\mathtt{v} :=
  D \rho( \mathtt{p}^0 , \mathtt{s}^0 ). ( \mathtt{P}^\eps , \mathtt{S}^{\eps,\flat} ) .\X_{\mathtt{v}^0} \,
\mathtt{v}^0 .
\end{array}
\right.
 \end{eqnarray*}
With the previous remark, we get
 \begin{eqnarray*}
 \mathfrak{P}^b_\mathtt{V} =  \{ \phi X (\rho( \mathtt{p}^0 ,  \mathtt{s}^0 + \tilde{\mathtt{S}}^0 ) 
 -  \rho( \mathtt{p}^0 , \mathtt{s}^0 )) \}|_{X=\frac{\varphi(x)}{\eps}} .
  \end{eqnarray*}
%%%
Proceed in the same way for the second equation of $(\ref{aller})$, we get the equation:
 \begin{eqnarray*}
  \alpha ( \mathtt{p}^\eps , \mathtt{s}^\eps ) .\X_{\mathtt{v}^\eps } \,  \mathtt{P}^\eps + \mathrm{div} \,  \mathtt{V}^\eps 
+ \mathfrak{P}_\mathtt{p}  = O(\eps ) \quad
  \mathrm{with}  \quad 
  \mathfrak{P}_\mathtt{p} := \mathfrak{P}^a_\mathtt{p} + \mathfrak{P}^b_\mathtt{p} + \mathfrak{P}^c_\mathtt{p} 
  \\ \nonumber  
  \mathrm{where}  \quad 
  \left\{
\begin{array}{l}
\mathfrak{P}^a_\mathtt{p} :=  \alpha ( \mathtt{p}^0 , \mathtt{s}^\eps ). \mathtt{V}^\eps . \nabla \mathtt{p}^0 ,
\vspace{0.2cm}\\ \mathfrak{P}^b_\mathtt{p} := \eps^{-1} \, ( \alpha  ( \mathtt{p}^0 ,  \mathtt{s}^0 + \tilde{\mathtt{S}}^0 ) 
- \alpha   ( \mathtt{p}^0 , \mathtt{s}^0
)).\X_{\mathtt{v}^0} \, \mathtt{p}^0
\vspace{0.2cm}\\  \qquad  =    \{ \phi X (\alpha  ( \mathtt{p}^0 ,  \mathtt{s}^0 + \tilde{\mathtt{S}}^0 ) 
 - \alpha  ( \mathtt{p}^0 , \mathtt{s}^0 )) \}|_{X=\frac{\varphi(x)}{\eps}} ,
\vspace{0.2cm}\\ \mathfrak{P}^c_\mathtt{p} := 
D \alpha (\mathtt{p}^0 , \mathtt{s}^0 ). (\mathtt{P}^\eps ,\mathtt{S}^{\eps,\flat}) .\X_{\mathtt{v}^0} \,
\mathtt{p}^0 .
\end{array}
\right.
 \end{eqnarray*}
In other words, the unknown $U^\eps := (\mathtt{V}^\eps ,\mathtt{P}^\eps , \mathtt{s}^\eps )$ verify the Euler system $(\ref{aller})$
except  the perturbation terms  $\mathfrak{P}_\mathtt{V}$, $  \mathfrak{P}_\mathtt{P}$, and 
 $O(\eps )$.
We see that the terms $\mathfrak{P}_\mathtt{v}^a$, $\mathfrak{P}_\mathtt{v}^b$, $\mathfrak{P}_\mathtt{p}^a$ and 
$\mathfrak{P}_\mathtt{p}^b$ do not have any singular factor with respect to $\eps$.
This is a consequence of  $(\ref{cond})$.
These terms are expressed in function of the unknown $U^\eps$, of the ground state $u^0$  and of the boundary layer $\tilde{\mathtt{S}}^0$.
The terms  $\mathfrak{P}^c_\mathtt{v}$ and $\mathfrak{P}^c_\mathtt{p}$, them, involve $\mathtt{S}^{\eps,\flat}$.
If we try to eliminate $\mathtt{S}^{\eps,\flat}$ via $(\ref{SSS})$, we involve  the unknown $U^\eps$ in a singular way  
via the term $\eps^{-1} \, \mathtt{s}^\eps$. 
We overcome this difficulty  in Lemma  $\ref{psg}$  using that the terms  $\X_{\mathtt{v}^0}
\, \mathtt{p}^0$ and  $\X_{\mathtt{v}^0} \, \mathtt{v}^0$ are respectively in factor of  $\mathfrak{P}^c_\mathtt{v}$ and
$\mathfrak{P}^c_\mathtt{p}$.
 Moreover the terms  $\mathfrak{P}_\mathtt{V}$ and $  \mathfrak{P}_\mathtt{P}$ are affine with respect to
  $ (\mathtt{V}^\eps ,\mathtt{P}^\eps )$. 
  
   The profile equations are linear, thanks to some original \textbf{transparency} properties of the Euler system. 
On one hand, the entropy boundary layer profile $\tilde{\mathtt{S}}^0$ verifies a transport equation
 which is linear with respect to the entropy (cf. equation $(\ref{clem})$). 
On the other hand, the amplitude of the boundary layer on the tangential velocity is weak (of order $\eps$) and the boundary is
characteristic for a linearly degenerate field. 
Thanks to this, the tangential velocity  boundary layer  profile $\tilde{\mathtt{V}}_\mathfrak{t}$ satisfies a linear equation,
 without Burgers-like nonlinearity
(cf. equation $(\ref{poli})$).
This is a transparency phenomenon analogous to the one observed in \cite{moi3}.
A interesting point is that such transparency phenomena does not occur for large amplitude high frequency entropy waves (see
Theorem $3.9$ of  \cite{cgm2}).
\\
\\  \textit{In section $\ref{stab}$}  we are interested in the \textbf{existence} (cf. Theorem $\ref{di}$) and
 the \textbf{propagation} (cf. Theorem $\ref{propa}$) of \textbf{exact solutions}
 of $(\ref{aller})$ asymptotic to approximate solutions obtained by truncating formal solutions 
 constructed in section $\ref{secform}$.
 Theorem $\ref{int}$ given in the introduction is a consequence of Theorem $\ref{try}$ and Theorem  $\ref{di}$. 
  After a reduction (cf. Prop. $\ref{Reduction2}$, subsection $\ref{Reduction}$), 
 we will face a singular perturbation problem because of boundary layers which corresponds to variations
in $\frac{\varphi(x)}{\eps}$. 
 More precisely we deal with
  a family of quasi-linear symmetric hyperbolic boundary value problem.
 As for the originating Euler problem, the boundary is conservative and characteristic of constant multiplicity.
 To tackle this characteristic problem  we get inspired by the paper \cite{Gu90} of O.Gu\`es which uses 
 the notion of conormal regularity and  the spaces 
\begin{eqnarray*}
E^m (T) := \{   u  \in L^2  (T)  /  \quad 
\sum_{0 \leqslant 2k+|\alpha| \leqslant m}  || \partial_\mathfrak{n}^k  \, Z^\alpha \, u ||_{L^2 ( T )}  < \infty  \} ,
\end{eqnarray*}
 with $\alpha := (\alpha_0 ,..., \alpha_d ) \in \N^{d+1}$, $|\alpha|:= \alpha_0 +...+ \alpha_d $, $Z^\alpha := Z_0^{\alpha_0} ...
Z_d^{\alpha_d} $ where $( Z_0,...,Z_d) $ generates the  algebra of $C^\infty$ tangent vector to $\Gamma$.
To simplify, we  denote $L^2 (T) :=  L^2 ((0,T) \times \Omega ) $.
 For these spaces  one normal derivative corresponds to two conormal derivatives.
 We adapt the method of \cite{Gu90} by substituting the derivative $\eps \partial_\mathfrak{n} $ 
 to the derivative $\partial_\mathfrak{n} $ 
in order to obtain uniform estimates and will use the following subsets of $L^2 (T)^{]0,1]}$:
\begin{eqnarray*}
\mathbf{E}^m (T) := \{   (u^\eps)_\eps   \in L^2 (T ) /  \quad 
 \mathrm{sup}_{0 < \eps \leqslant 1}  
 \sum_{0 \leqslant 2k+|\alpha| \leqslant m}  ||  (\eps \partial_\mathfrak{n} )^k \, Z^\alpha \, u^\eps ||_{L^2 ( T )}  < \infty  \} .
\end{eqnarray*}
This idea to  use some derivatives with $\eps$ in factor for some  singular perturbation problems is natural and was also used in
the papers of  \cite{Gu93},  \cite{Gu92} with the $\eps$-stratified notion, \cite{cgm2} with
 the $\boldsymbol{\eps}$\textbf{-conormal} notion. 
Here, this idea is applied to (characteristic) boundary value problem and anisotropic Sobolev spaces.

  At first look, this system we obtained is singular with respect to $\eps$ but a trick allows to overcome this false singularity
(see subsection $\ref{Reduction}$).
 We will use a family of iterative schemes. 
Thus we will supply in subsection $\ref{Linesti}$ linear estimates which are the core the proof. 
We will successively perform $L^2$ estimates, conormal estimates and normal estimates.
A main difficulty lies in the way to deal with commutators (cf. Proposition $\ref{lapin}$).
This strategy yields exact solutions till $T_0$.
The proof of Theorem  $\ref{di}$ needs carefulness about the existence of compatible initial data.      
Subsection  $\ref{bobo}$ is devoted to this question.
  
 It is possible to obtain $L^\infty$ estimates, in spite of the fact that $d \geqslant 1$. 
We refer to  papers \cite{m2}, \cite{sonique} of G.M\'etivier, paper \cite{RR} of J.Rauch and M.Reed and paper \cite{Gu90}.
This idea is still relevant when adapting $\eps$-conormal regularity to characteristic boundary value problem.
Therefore we can weaken the regularity of the solution and prove a propagation result for some solutions admitting only
one normal derivative in $L^2$. 
We introduce the following subsets of $L^2 (T)^{]0,1]}$:
\begin{eqnarray*}
\mathbf{A}^m (T) := \{   (u^\eps)_\eps   /  \quad 
 \mathrm{sup}_{0 < \eps \leqslant 1}  
 \sum_{0 \leqslant |\alpha| \leqslant m}  ||  Z^\alpha  \, u^\eps ||_{L^2 ( T)} 
 +  \sum_{0 \leqslant  |\alpha| \leqslant m-2}  ||  \eps \partial_\mathfrak{n}  \, Z^\alpha  \, u^\eps ||_{L^2 ( T)} < \infty  \} .
\end{eqnarray*}
We will also use some norms built on $L^\infty$. 
Because the boundary is characteristic, we will need not only the Lipschitz norms but higher order $L^\infty$ control, 
as O. Gu\`es in \cite{Gu90} and G. M\'etivier in \cite{sonique}.
We will denote by $L^\infty (T) $ the space     $ L^\infty (T)  :=  L^\infty ((0,T) \times \Omega ) $.
 We will introduce the norms
\begin{eqnarray*}
 ||u||^*_{\eps,T}  &:=& \sum_{0 \leqslant |\alpha| \leqslant 2} ||Z^\alpha \, u||_{L^\infty (T)}
 +  \sum_{0 \leqslant |\alpha| \leqslant 1} ||Z^\alpha  \, \eps \partial_\mathfrak{n}  \, u||_{L^\infty (T)} ,
\end{eqnarray*}
and the following subsets of $L^\infty (T)^{]0,1]}$:
\begin{eqnarray*}
\mathbf{\Lambda}^m (T) := \{  (u^\eps)_\eps    /  \quad 
 \mathrm{sup}_{0 < \eps \leqslant 1} \, ||u^\eps ||^*_{\eps,T} < \infty  \} . 
\end{eqnarray*}
Theorem $\ref{propa2}$ states a propagation result in the spaces $\mathbf{A}^m (T) \cap \mathbf{\Lambda}^m (T) $.  
  
 One quality of our method is that we need approximate solutions with only a few profiles.
The minimum number of profiles required is linked to the lost of a factor $\eps^{\frac{1}{2}}$ in a Sobolev embedding Lemma (Lemma
$\ref{sobo}$). 
Let us explain one motivation to minimize the number of profiles needed.
 In this paper, we consider a ground state $u^0$ in $H^\infty ((0,T_0) \times \Omega )$  and formal solutions with $H^\infty$
regularity. It could also be possible to extend to ground states of high but finite regularity. 

 \section{Forthcoming}
 \label{coming}

We plan to show in a further work that for more general ground states $u^0$ which are solutions of $(\ref{aller})$,
which verify the condition $\X_{\mathtt{v}^0} \,  \mathtt{v}^0 =0$ for $(t,x) \in (0,T_0) \times \Gamma$
 but which do not verify the
condition $\X_{\mathtt{v}^0} \,  \mathtt{p}^0 =0$ for $(t,x) \in (0,T_0) \times \Gamma$, 
it is still possible to construct, in small time, some nontrivial formal solutions of $(\ref{aller})$ but
of the more general form
 \begin{eqnarray*}
    \left\{
   \begin{array}{l}
  \mathtt{v}_\mathfrak{t}^\eps (t,x) = 
   \displaystyle{\sum_{j \geqslant 0}} \eps^j \,  \mathtt{V}_\mathfrak{t}^{j} (t,x, \frac{\varphi(x) }{\eps}) ,
\vspace{0.2cm}  \\ \mathtt{v}_\mathfrak{n}^\eps (t,x) =  \mathtt{v}_\mathfrak{n}^0 (t,x) 
  + \displaystyle{\sum_{j \geqslant 1}} \eps^j \,  \mathtt{V}_\mathfrak{n}^{j} (t,x, \frac{\varphi(x)}{\eps}) ,
\vspace{0.2cm}  \\ \mathtt{p}^\eps (t,x) =  \mathtt{p}^0 (t,x) + \eps \mathtt{p}^1 (t,x) 
  +  \displaystyle{\sum_{j \geqslant 2}} \eps^j \,  \mathtt{P}^{j} (t,x, \frac{\varphi(x)}{\eps}) ,
  \vspace{0.2cm}  \\ \mathtt{s}^\eps (t,x) = \displaystyle{\sum_{j \geqslant 0}} \eps^j \,  \mathtt{S}^{j} (t,x,
  \frac{\varphi(x)}{\eps}) .
    \end{array}
\right.
    \end{eqnarray*}
Moreover we will show  that we can prescribe arbitrary initial values for the $(\mathtt{V}_\mathfrak{t}^{j} )_{j \in \N}$ 
and the $(\mathtt{S}^{j} )_{j \in \N}$.
  However these formal solutions can be unstable.

\section{Setting of the notations}
\label{setting}

To simplify and avoid heavy notations, we will consider from now on that the domain $\Omega$ 
is the half-space $ \Omega := \{ x \in \R^d / \quad x_d \geqslant 0 \}$.
 This assumption does not change the mathematical analysis of the problem. 
 We fix a notation. 
If $A$ is a $d+2$ by $d+2$ square matrix, we denote by $A^\star$ the $d+1$ by $d+1$ extracted square matrix which
contains the $d+1$ first rows of the $d+1$ first lines.
With $u=(\mathtt{v}, \mathtt{p},\mathtt{s})$, the Euler system is of the form
 \begin{eqnarray}
 \label{po}
 \X_\mathtt{v} \, u + \M (u,\partial_x ) u =0 
 \end{eqnarray}
 with
 \begin{eqnarray*}
\M (u,\xi) := 
 \begin{bmatrix}
 \M^\star (u,\xi) &0 
 \\ 0 & 0 
  \end{bmatrix}
  ,\quad 
  \M^\star (u,\xi) := 
 \begin{bmatrix}
  0 &  \mathtt{\rho}^{-1} \ \xi
\\  \mathtt{\alpha}^{-1} \  {}^{t} \xi & 0 
\end{bmatrix},
\end{eqnarray*}
We recall that we assume that  $\alpha (\mathtt{p},\mathtt{s}) > 0$ for all $\mathtt{p}$ and  $\mathtt{s}$. 
We  denote by
\begin{eqnarray*}
S(u) :=
\begin{bmatrix}
 S^\star (u) & 0
\\   0 &  1
  \end{bmatrix}
   ,\quad 
 S^\star (u) :=
\begin{bmatrix}
     \mathtt{\rho} I_d & 0 
\\    0  & \alpha 
 \end{bmatrix} .
\end{eqnarray*}
Multiply Eq. $(\ref{po})$ by  the  matrix $S(u)$
to get the equation 
\begin{eqnarray}
 \label{po2}
\mathfrak{L}(u,\partial) u =0 
 \end{eqnarray}
where  $\mathfrak{L}(u,\partial) := S(u)  \X_\mathtt{v}  + \L (\partial_x )$ and  for all $\xi \in \R^d$, 
\begin{eqnarray*}
\L(\xi) :=  S(u) \M (u,\xi) =
\begin{bmatrix}
 \L^{\star} (\xi) &0 
 \\ 0 & 0 
\end{bmatrix} 
 ,\quad 
 \L^{\star} (\xi) :=
\begin{bmatrix}
   0 &  \xi
\\    {}^{t} \xi  & 0 
  \end{bmatrix}.
\end{eqnarray*}
We also introduce the operator  $\mathfrak{L}^{\star} (u,\partial) := S^{\star}(u)  \X_\mathtt{v}  + \L^{\star} (\partial_x )$.
The matrix $S(u)$ is symmetric  definite nonnegative and depends in a $C^\infty $ way of  $u$.
The system $(\ref{po2})$ is therefore symmetric hyperbolic.
For this system, the boundary conditions $\mathtt{v}. \mathfrak{n} = 0$  are conservative.
We denote by $v:=  (\mathtt{v}, \mathtt{p})$ and by   $ w:=\mathtt{s}$.
We introduce for all $\xi \in \R^d - \{ 0 \}$ the subspace $ \mathbb{F} (u,\xi) := \ker \M(u,\xi)$ of $ \R^N $.
We denote by $(e_1,...,e_{d+2})$ the canonical basis of $\R^{d+2}$.
Notice that $ \mathbb{F} (u,\xi) = \xi^{\bot} \oplus \R e_{d+2} $ and that
$\mathbf{\lambda} (u,\xi):= \mathtt{v}. \xi$ is a linearly degenerate eigenvalue with constant multiplicity $d> 0 $ i.e.
\begin{eqnarray*}
\left.
\begin{array}{c}
r. \nabla_u \mathbf{\lambda} (u,\xi) = 0 ,\quad \mathrm{for} \  \mathrm{all} \ r \in \mathbb{F} (u,\xi) ,
  
\\  \dim \mathbb{F} (u,\xi) = d ,
\end{array}
\right\}
\quad \mathrm{for} \  \mathrm{all} \  (u,\xi) \in \R^{d+2} \times (\R^d - \{ 0 \}) .
\end{eqnarray*}
Notice that $\{ v=0 \} \subset \ker \M(u,\xi)$, for all  $ \xi \in \R^d - \{ 0 \}$.
We  denote by $P_0$ the  orthogonal projector on $\ker \L_d$ i.e.  
\begin{eqnarray*}
P_0 := 
\begin{bmatrix}
P_0^{\star} & 0
\\ 0 & 1 
\end{bmatrix},
\end{eqnarray*}
where $P_0^{\star}$ is a $(d+1) \times (d+1) $ matrix 
\begin{eqnarray*}
P_0^{\star} := 
\begin{bmatrix}
 I_{d-1} & 0 & 0
\\ 0 & 0 & 0
\\ 0 & 0 & 0
\end{bmatrix}.
\end{eqnarray*}
We will often split $\mathcal{U} \in \P(T)$ into $\mathcal{U} =( \mathcal{V},  \mathcal{W})$.
The function  $\mathcal{V}$ (respectively  $\mathcal{W}$) takes its values in $\R^{d+1}$ (resp. $\R$).
 We will sometimes split $\mathcal{V}$ into  $\mathcal{V}=(\mathtt{V},\mathtt{P})$ and $\mathtt{V}$ 
 into $\mathtt{V} := (\mathtt{V}_\mathfrak{t} , \mathtt{V}_d)$. 
 The function  $\mathtt{V}$ (respectively  $\mathtt{P}$) takes its values in $\R^{d}$ (resp. $\R$).
 The function  $\mathtt{V}_\mathfrak{t} $ (respectively  $\mathtt{V}_d$) takes its values in $\R^{d-1}$ (resp. $\R$).
By abuse of notations, we will say that $\mathcal{V}$, $\mathcal{W}$, $\mathtt{V}$, $\mathtt{V}_\mathfrak{t} $,  $\mathtt{V}_d$ 
and $\mathtt{P}$ are in $\P(T)$ even if they do not take values in $\R^{d+2}$. 

\section{Overdetermined ground states} 
\label{ground}

In this section, we  prove the existence of ground states $(\mathtt{v}^0 , \mathtt{p}^0 ,\mathtt{s}^0 ) \in H^\infty$ solutions of
 $(\ref{aller})$ and verifying the conditions $(\ref{pro})$.
 
\begin{theo}
\label{tt}
Given $h \in C^1 (\R^d ,\R^d )$ such that $h'(x)$ is nilpotent in a neighborhood of $0$ and such that $h_d (x) = 0$ when $x_d =0$,
there exists a local $C^1$ solution $\mathtt{v}$ of the initial boundary value problem:
\begin{eqnarray*}
\partial_t \mathtt{v} + (\mathtt{v} . \nabla_x ) \mathtt{v} = 0 , \quad \mathrm{div}_x \mathtt{v} =0 ,\quad \mathtt{v}_d |_{x_d =0 }
=0 ,\quad \mathtt{v} |_{t =0 } = h .
\end{eqnarray*}
\end{theo}
\begin{proof}
Local existence of $C^1$ solution of multidimensional Burgers equation can be achieved with characteristic method. 
The classical relation $\mathtt{v} (t,x+t h(x)) = h(x)$ holds in a neighborhood of $0$. 
Because $h_d (x) = 0$ when $x_d =0$, we get $\mathtt{v}_d  (t,x) = 0$ for $x_d =0$.
The divergence free relation is a consequence of the nilpotence of $h'(x)$ as proved in Theorem $2.6$ of \cite{cgm2}.
\end{proof}
Assume that $d=2$ for a moment and let us give some examples of initial velocity $h$
 which satisfy the assumptions of Theorem $\ref{tt}$.
If $h$ satisfies $h_1 (x)= x_2$ and $h_2 (x) =0$ in a neighborhood  of $0$, then $h$ is convenient. 
We now detail a more general process to get convenient $h$.
Let $F$ be a function in $ C^1 (\R_+ ,\R )$ such that $F(0)=0$ and $F'(x) > 0$ for all $x \in \R_+ $.
Let $a_{init} \in C^1 (\R_+ , \R)$. 
There exists a local solution $a \in C^1 (\R \times \R_+ ,\R )$ of the scalar initial boundary value problem  
\begin{eqnarray*}
\partial_1 a + \partial_2 a =0 , \quad a|_{x_2 =0} =0  , \quad a|_{x_1 =0} = a_{init} .
\end{eqnarray*}
Consider $h(x) := (a(x),F \circ a (x))$. 
Then $h$ satisfies the assumptions of Theorem $\ref{tt}$.

 Of course, because conditions $(\ref{pro})$ imply conditions  $(\ref{cond})$, what precedes proves the existence of $T_0 > 0$ and
of a ground state $u^0 := (v^0 , w^0 ) \in H^\infty ((0,T_0 ) \times \Omega )$ solution of  $(\ref{aller})$ and verifying the
conditions $(\ref{cond})$.
From now on, we will assume that such a ground state is given.

\section{Reduction of the system} 
\label{redbull}

We look for solutions of  $(\ref{aller})$  of the  form
\begin{eqnarray*}
u^\eps =(v^\eps := v^0 + \eps V^\eps , w^\eps := W^\eps ).
\end{eqnarray*}
We are going to realize a singular change of unknown, by looking for an equation for $U^\eps := (V^\eps , W^\eps )$. 

\begin{prop}
\label{Re1}

The function $u^\eps $ is solution of the equation $(\ref{po})$ if and only if the function $U^\eps$ verify the equation 
 \begin{eqnarray}
 \label{Re2}
\mathfrak{L}(u^\eps ,\partial)  U^\eps +  K^\eps  =0 ,
 \end{eqnarray}
 where
 \begin{eqnarray}
 \label{K1} 
  K^\eps  := \frac{1}{\eps} 
 \begin{bmatrix}   
 K_1 ( v^0 , \partial v^0  , u^\eps) 
 \\  0  
 \end{bmatrix} \quad  \mathrm{with} \  
 K_1 (v^0 , \partial v^0 , u^\eps) =  \mathfrak{L}^\star (u^\eps ,\partial ) v^0 .
  \end{eqnarray}

\end{prop}
\begin{proof}
The equation $(\ref{po})$ is equivalent to the two following equations:
\begin{eqnarray}
 \label{Re3}
 \X_{\mathtt{v}^\eps} \,   v^\eps + \M^\star ( u^\eps , \partial_x ) v^\eps =0 ,
\quad  \X_{\mathtt{v}^\eps} \,   w^\eps  =0 ,
\end{eqnarray} 
By dividing  the first equation by $\eps$, it's also  equivalent to 
\begin{eqnarray*}
 \X_{\mathtt{v}^\eps} \,   V^\eps + \M^\star ( u^\eps , \partial_x ) V^\eps 
 + \frac{1}{\eps} (  \X_{\mathtt{v}^\eps} \,   v^0 + \M^\star ( u^\eps , \partial_x ) v^0 ) = 0 ,
\quad  \X_{\mathtt{v}^\eps} \,   W^\eps =0 .
\end{eqnarray*}
 We get:
\begin{eqnarray}
\label{pet0}
 \X_{ \mathtt{v}^\eps} U^\eps +
  \mathcal{M} U^\eps + 
 \begin{bmatrix} 
 \frac{1}{\eps} (\X_{ \mathtt{v}^\eps} v^0 +   \mathcal{M}^\star ( u^\eps , \partial_x) v^0 )
 \\ 0 
 \end{bmatrix}
 =0  .
\end{eqnarray}
 We multiply $(\ref{pet0})$ to the left by the matrix  $S (u^\eps )$ to end the proof. 
\end{proof}

The underlying idea of the previous proposition is that because $v^0$ is given as a ground state, $U^\eps$ is the real unknown.
For each $\eps$, we obtain that $U^\eps$ satisfies a hyperbolic system. However in order to achieve an asymptotic analysis for
$\eps \rightarrow 0^+$, we can be \textit{a priori} worried about the singular factor $\eps^{-1}$ within $K^\eps$.
Let us recall the way \cite{cgm2} deals with this term, expliciting their calculus for this particular case of the Euler system.
First because $\X_{ \mathtt{v}^\eps} v^0 = \X_{\mathtt{v}^0} \,  v^0 + \eps \mathtt{V}^\eps . \nabla v^0 $, we get 
\begin{eqnarray*} 
K^\eps := \begin{bmatrix}
\frac{1}{\eps} S^\star (u^\eps ) \X_{\mathtt{v}^0} \,  v^0  + \frac{1}{\eps} \L^\star (\partial_x ) v^0 
+  S^\star (u^\eps )  \mathtt{V}^\eps . \nabla v^0 
 \\ 0 
 \end{bmatrix}.
 \end{eqnarray*}
 Remember that \cite{cgm2} use the condition  $(\ref{pro})$ on the ground state $u^0$
  which are stronger than the condition $(\ref{cond})$
  we use in the present paper.
 The condition  $(\ref{pro})$ reads  
$ \X_{\mathtt{v}^0} \,  v^0 = \L^\star (\partial_x ) v^0 = 0$ and 
  \begin{eqnarray*}
     K^\eps = \begin{bmatrix} S^\star (u^\eps )  \mathtt{V}^\eps . \nabla 
v^0 \\ 0 
 \end{bmatrix}
 \end{eqnarray*}
  does not contain any singular factor $\eps^{-1}$ anymore.
The following proposition will show that under condition $(\ref{cond})$ the term $K^\eps$ does not contain 
any singular factor $\eps^{-1}$ too.
In order to help the reader, we give some hints about our strategy.
We will take into account that we search $W^\eps$ of the form 
\begin{eqnarray*} 
W^{\eps} (t,x)  = \mathcal{W}^0 (t,x,\frac{x_d}{\eps}) + \eps W^{\eps,\flat} (t,x)
\end{eqnarray*}
 with $\mathcal{W}^0
\in \P (T )$ (cf. $(\ref{defp})$ for the definition) and $ \underline{\mathcal{W}}^0 =w^0$. 
Under the conditions $(\ref{pro})$, we have only expanded $v^\eps$.
Because we assumed that $\X_{ \mathtt{v}^0} v^0 =  0$, the term  $S^\star (u^\eps ) \X_{\mathtt{v}^0} \,  v^0 $ vanished.
Under the localized condition $(\ref{cond})$, the idea is tricker.
We will expand $W^\eps$ too, into $W^\eps = w^0 + \tilde{\mathcal{W}}^0 +\eps W^{\eps,\flat} $ within $S^\star (u^\eps )$.
Imagine in a first time that $\tilde{\mathcal{W}}^0$  is identically zero.
Then we obtain by a Taylor first order expansion that
\begin{eqnarray*} 
K^\eps =
 \begin{bmatrix} 
 \frac{1}{\eps} \mathfrak{L}(u^0,\partial)  v^0  + \mathrm{terms} \ \mathrm{of} \  \mathrm{order} \ \eps^0 
  \\ 0 
 \end{bmatrix} .
  \end{eqnarray*}
 Because the ground state $u^0$ satisfies $(\ref{aller})$, the first term in the right member is equal to zero and $K^\eps$ does
 not appear singular with respect to  $\eps$ anymore.
 Of course we want to deal with some nonvanishing function $\tilde{\mathcal{W}}^0$. 
A key difference with paper \cite{cgm2} is that here $\tilde{\mathcal{W}}^0$ denotes a boundary layer and we will see in the
 following proposition that its singular contribution within $K^\eps$ contains the trace of $\X_{ \mathtt{v}^0} v^0 $ 
 on the boundary $\Gamma := \{ x_d =0 \}$ as a factor.
Therefore the localized conditions $(\ref{cond})$ allow to conclude.
In prevision of the following sections, we will be careful with the way the term $K^\eps$ depends of $V^\eps$ and
$W^{\eps,\flat}$.

In order to avoid heavy notations, we will omit the two first arguments and write $ K^\eps_1 (u^\eps)$ 
instead of  $ K^\eps_1  (v^0 , \partial v^0 , u^\eps)$.
Furthermore, we want now to use the special form of the solutions $u^\eps$ we are looking for. 
Thus we will write $ K^\eps_1  (v^\eps , w^\eps)$.
We also introduce some notations.
A Taylor expansion proves that there exist two $C^\infty$ functions $v^{0,\flat}$ and $w^{0,\flat}$
 such that $ v^0 = \rond{v}^0 + x_d \, v^{0,\flat} $ 
and $ w^0 = \rond{w}^0 +  x_d \,  w^{0,\flat} $, where $\rond{u}$ is the trace of $u$ on $x_d =0$. 

\begin{prop} 
\label{sev}

Assume that $W^{\eps}$ is of the form 
\begin{eqnarray} 
\label{qz}
W^{\eps} (t,x)  = \mathcal{W}^0 (t,x,\frac{x_d}{\eps}) + \eps W^{\eps,\flat} (t,x) .
\end{eqnarray}
 where $ \mathcal{W}^0 \in \P(T) $
 with  $ \un{\mathcal{W}}^0 = w^0$.
Then there is a  $C^\infty $ matrix  $K_1^{\flat} $ such that  
\begin{eqnarray} 
 \label{singul}   K_{1} ( u^\eps ) = 
 \eps K_{1}^{\flat} (\eps, x_d ,\rond{v}^0 , v^{0,\flat} , V^\eps , \eps V^\eps , \rond{w}^0 , w^{0,\flat} ,  \tilde{\mathcal{W}}^0
 ,\frac{x_d}{\eps} \tilde{\mathcal{W}}^0 , W^{\eps,\flat} , \eps W^{\eps,\flat})
 \end{eqnarray}
 where $K_{1}^{\flat}$ is affine with respect to its fifth argument $V^\eps$ and affine with respect to its eleventh argument
 $W^{\eps,\flat}$ with $\X_{\mathtt{v}^0} \,  v^0$ as factor in the leading coefficient.  In  $(\ref{singul})$, to avoid heavy
 notations we denote
 $\tilde{\mathcal{W}}^0$ instead of $\tilde{\mathcal{W}}^0  (t,x,\frac{x_d}{\eps})$.

\end{prop}
\begin{proof} 
We begin giving a technical lemma. We will use it twice.

\begin{lem}
\label{tech}

There exist some $d+1$  square matrices $K^\flat_{1,1} (u_1 , u_2 )$  and some  $K^\flat_{1,2} (u_1
, u_2 ) \in \R^{d+1}$ both $C^\infty$ with respect to  to their arguments
$u_1 := (v_1 , w_1 )$,  $u_2 := (v_2 , w_2 ) \in \R^{d+1} \times \R   $ such that $K^\flat_{1,2}$ 
has  $\X_{\mathtt{v}_1} v^0$ as factor and that 
\begin{eqnarray} 
 \label{zou} \quad K_1 (v_1 + v_2 , w_1 + w_2 ) = K_1 (u_1) +  K_{1,1}^{\flat} (u_1 , u_2 ) .v_2
 + w_2 \,  K_{1,2}^{\flat} ( u_1 , u_2 ).
 \end{eqnarray}
\end{lem}

\begin{proof} 
We will proceed in two steps.  

 \begin{enumerate}[$1$.]

\item  A Taylor first order expansion yields the existence of matrices $S_{1}^{\star,\flat} (u_1 , u_2 )$, 
$S_{2}^{\star,\flat} (u_1 , u_2)$
 such that for all $(u_1 ,u_2)$, 
\begin{eqnarray}
\label{tuok} \quad  S^{\star} (v_1 +v_2 , w_1 + w_2) = S^{\star} (u_1 ) + v_2 . S_{1}^{\star,\flat} (u_1 , u_2) 
+  w_2 . S_{2}^{\star,\flat} (u_1 ,u_2)  .
 \end{eqnarray}
Here $S_{1}^{\star,\flat} (u_1 , u_2)$ is a $(d+1)^2 \times (d+1)$ matrix and 
$S_{2}^{\star,\flat} ( u_1 , u_2)$  is a $d+1 $ square matrix.
We write in a block form the matrix 
\begin{eqnarray*}
S_{1}^{\star,\flat} := 
\begin{bmatrix}
S_{1,1}^{\star,\flat}
\\ .
\\ .
\\ .
\\ S_{1,d+1}^{\star,\flat}
\end{bmatrix} ,
 \end{eqnarray*}
where the $(S_{1,j}^{\star,\flat} )_{j=1,...,d+1}$ are some $d+1$ square matrices. 
In Eq. $(\ref{tuok})$, 
the notation $ v. S_{1}^{\star, \flat}$ stands for the $d+1$ square matrix
 $v . S_{1}^{\star,\flat} := v_1 S_{1,1}^{\star,\flat}  + ... +  v_{d+1} S_{1,d+1}^{\star,\flat} $.
  
  \item We introduce, for all $(u_1 , u_2)$, the $d+1$ by $d+1$ matrix $K_{1,1}^{\flat} ( u_1 ,u_2) $ 
 such that for all $z=(\mathtt{z},z_{d+1}) \in \R^{d} \times \R$,
\begin{eqnarray*}
 K_{1,1}^{\flat} (u_1 , u_2).z &:=& z. S_{1}^{\star,\flat} ( u_1 , u_2) . \X_{ \mathtt{v}_1 + \mathtt{v}_2   }  \, v^0 
+ ( S^{\star} ( u_1 ) + w_2 . S_{2}^{\star,\flat} ). \mathtt{z}. \nabla v^0 ,
 \end{eqnarray*}
and the vector
\begin{eqnarray*}
K_{1,2}^{\flat} ( u_1 ,u_2) :=   S_{2}^{\star,\flat} (u_1 , u_2) . \X_{\mathtt{v}_1 } \,  v^0 \in \R^{d+1} .
\end{eqnarray*}
 As by definition,  for all $(v ,w ) \in \R^d \times \R$, 
 \begin{eqnarray*}
  K_1 (v ,w)  :=  S^\star (v ,w ) \X_{\mathtt{v}} \,  v^0 + \L^\star (\partial_x ) v^0 , 
  \end{eqnarray*}
   we get for all $(u_1 , u_2)$, 
 \begin{eqnarray*}
  K_1 (v_1 + v_2 , w_1 + w_2) - K_1 (v_1  , w_1 ) =  
  S^\star (v_1 + v_2 ,w_1 + w_2  ) \X_{\mathtt{v}_1 + \mathtt{v}_2}  \, v^0 
  \\ - S^\star (v_1 ,w_1) \X_{\mathtt{v}_1}  \, v^0 .
  \end{eqnarray*}
 Thanks to Eq. $(\ref{tuok})$ and because $\X_{\mathtt{v}_1 + \mathtt{v}_2   } v^0 = \X_{\mathtt{v}_1 } v^0 + \mathtt{v}_2 .\nabla
 v^0 $,  we get $(\ref{zou})$.
\end{enumerate}
\end{proof}

This technical lemma given, we will proceed in three steps. 

\begin{step}[First order expansion]

We are going to prove that there exists a $C^\infty $ function $ K_1^{\flat,1}$, affine
with respect to its third variable and affine
with respect to its sixth variable with $\X_{\mathtt{v}^0} \,  v^0$ as factor in the leading coefficient such that 
\begin{eqnarray} 
 \label{zoub}
  K_1 ( u^\eps ) =  K_1 (v^0  , \mathcal{W}^0 )
   + \eps  K_1^{\flat,1} (\eps , v^0 , V^\eps , \eps V^\eps , \mathcal{W}^0 ,  W^{\eps,\flat} ,\eps W^{\eps,\flat} ).
 \end{eqnarray}
\end{step}

We use a first time Lemma $\ref{tech}$ and apply $(\ref{zou})$ to 
\begin{eqnarray*}
(v_1 , w_1 ,v_2 , w_2)  =(v^0 ,\mathcal{W}^0 ,\eps V^\eps ,\eps W^{\eps,\flat}), 
 \end{eqnarray*}

we get $(\ref{zoub})$ with 
\begin{eqnarray*}
 K_1^{\flat,1} (\eps , v^0 , V^\eps , \eps V^\eps , \mathcal{W}^0 ,  W^{\eps,\flat}, \eps W^{\eps,\flat} ) &:=&
 V^\eps  K_{1,1}^{\flat}   (v^0  , \mathcal{W}^0 , \eps V^\eps , \eps W^{\eps,\flat})
\\  &&+ W^{\eps,\flat}  K_{1,2}^{\flat}   (v^0  ,  \mathcal{W}^0 ,\eps V^\eps ,\eps W^{\eps,\flat}). 
 \end{eqnarray*}
It is clear that the function $K_1^{\flat,1}$ is affine
with respect to its third variable and affine
with respect to its sixth variable with $\X_{\mathtt{v}^0} \,  v^0$ as factor in the leading coefficient.

\begin{step}[Do $x_d =0$]

  We are going to prove that there exists a $C^\infty $ function $K_1^{\flat,1}$, 
 affine with respect to its fifth argument and affine with respect to its eleventh argument with $\X_{\mathtt{v}^0} \, 
 v^0$ as factor in the leading coefficient, such that       
\begin{eqnarray} 
 \label{sing}   K_{1} ( u^\eps ) &= &  K_{1} (u^0 )
  + (  K_{1} (\rond{v}^0   ,\rond{w}^0 + \tilde{\mathcal{W}}^0) - K_{1} (\rond{v}^0  , \rond{w}^0 ))
 \nonumber \\  & +& \eps   K_{1}^{\flat} (\eps, x_d ,\rond{v}^0 , v^{0,\flat} , V^\eps , \eps V^\eps , \rond{w}^0 , 
 w^{0,\flat} ,  \tilde{\mathcal{W}}^0 , \frac{x_d}{\eps} \tilde{\mathcal{W}}^0 , W^{\eps,\flat} , \eps W^{\eps,\flat}) .
 \end{eqnarray}
\end{step} 

 Thus we have to do the analysis of  $K_1 (v^0 ,\mathcal{W}^0)$. 
We write naively 
\begin{eqnarray}
\label{N}
  K_1 (v^0  , \mathcal{W}^0 ) =  K_1 (u^0  ) + (  K_1 (v^0  , \mathcal{W}^0) -  K_1 (u^0 )) .
 \end{eqnarray}
and do an estimate for the error 
by substituting 
\begin{eqnarray*}
 K_1 (\rond{v}^0  , \rond{w}^0 + \tilde{\mathcal{W}}^0 ) -  K_1 (\rond{v}^0 ,\rond{w}^0 ) \quad
\mathrm{to} \  K_1 (v^0  , \mathcal{W}^0) -  K_1 ( u^0 )  .
 \end{eqnarray*}
\begin{lem}

There exists a  $C^\infty $ function $ K_1^{\flat,2}$ such that 
\begin{eqnarray}
\label{zoub2}
  K_1 (v^0  , \mathcal{W}^0) -  K_1 (u^0  ) &=& K_1 (\rond{v}^0 , \rond{w}^0 + \tilde{\mathcal{W}}^0 ) 
  - K_1 (\rond{v}^0  , \rond{w}^0 ) 
  \\ \nonumber && + \eps K_1^{\flat,2} (x_d ,\rond{v}^0 , v^{0,\flat} , \rond{w}^0  , w^{0,\flat}  ,\tilde{\mathcal{W}}^0 ,
  \frac{x_d}{\eps} \tilde{\mathcal{W}}^0 ).
 \end{eqnarray}
 \end{lem}
 \begin{proof}
 We will proceed in three steps.
   \begin{enumerate}[$1$.]
 
 \item We use use  Lemma $\ref{tech}$ again and apply respectively $(\ref{zou})$ to 
\begin{eqnarray*} 
(v_1 , w_1 , v_2 , w_2) &=& ( \rond{v}^0 ,  \rond{w}^0 + \tilde{\mathcal{W}}^0  ,  x_d  v^{0,\flat} ,   x_d  w^{0,\flat} ) 
\quad \mathrm{and} \  \mathrm{to} 
\\   (v_1 , w_1 , v_2 , w_2)  &=& ( \rond{v}^0 ,  \rond{w}^0  ,  x_d  v^{0,\flat} ,   x_d  w^{0,\flat} )
 \end{eqnarray*}
 we get 
\begin{eqnarray} 
\label{mlp} \qquad \left\{
\begin{array}{ccc}
  K_1 (v^0  , \mathcal{W}^0) &=&  K_1 (\rond{v}^0  , \rond{w}^0 + \tilde{\mathcal{W}}^0) 
  + x_d   K_{1}^{a} (x_d ,\rond{v}^0 , v^{0,\flat} ,  \rond{w}^0 ,  w^{0,\flat} ,  \tilde{\mathcal{W}}^0  ) ,
  \\   K_1 (v^0  , w^0 )  &=&  K_1 (\rond{v}^0  , \rond{w}^0 ) 
  + x_d   K_{1}^{b} (x_d ,\rond{v}^0 ,  v^{0,\flat} ,  \rond{w}^0 ,  w^{0,\flat}  ),
  \end{array}
  \right.
 \end{eqnarray}
with
\begin{eqnarray*} 
 K_{1}^{a} (x_d , \rond{v}^0 ,  v^{0,\flat} ,  \rond{w}^0 , w^{0,\flat}  ,  \tilde{\mathcal{W}}^0  )
    &:=&  v^{0,\flat}  K_{1,1}^{\flat}   ( \rond{v}^0  ,  \rond{w}^0  + \tilde{\mathcal{W}}^0  , x_d v^{0,\flat},  x_d  w^{0,\flat} ) 
 \\ &+&  w^{0,\flat}  K_{1,2}^{\flat}   (  \rond{v}^0 ,  \rond{w}^0 + \tilde{\mathcal{W}}^0 , x_d   v^{0,\flat} ,  x_d   w^{0,\flat} ),
\\  K_{1}^{b} (x_d , \rond{v}^0 , v^{0,\flat}  ,   \rond{w}^0  , w^{0,\flat} )
    &:=&   v^{0,\flat}  K_{1,1}^{\flat}  (\rond{v}^0 , \rond{w}^0 , x_d v^{0,\flat} , x_d  w^{0,\flat}) 
 \\ &+& w^{0,\flat}  K_{1,2}^{\flat}   ( \rond{v}^0  , \rond{w}^0 , x_d  v^{0,\flat},  x_d  w^{0,\flat} ).
 \end{eqnarray*}
 and we denote by $ K_1^{\flat,r,2} :=  K_{1}^{a} -  K_{1}^{b}$.
 
  \item By a first order Taylor expansion, we obtain that there exist some $(d+1) \times (d+1)$ matrices $G_{1,1}^{\flat}$
  and some $G_{1,2}^{\flat} \in \R^{d+1}$, $C^\infty$ with respect to their arguments
   \begin{eqnarray*} 
  u_1 :=  (v_1 , w_1 ) \ ,  u_2 := (v_2 , w_2 ) \in  \R^{d+1} \times \R \ ,  w_3 \in \R  
  \end{eqnarray*}

such that

\begin{eqnarray*}
K_{1,1}^{\flat} (v_1 ,w_1 + w_3 , u_2 ) - K_{1,1}^{\flat} (u_1 , u_2 ) = 
 w_3 . G_{1,1}^{\flat} (u_1 ,u_2 , w_3),
\\ K_{1,2}^{\flat} (v_1 ,w_1 + w_3 , u_2 ) - K_{1,2}^{\flat} (u_1 , u_2 ) = 
 w_3 .G_{1,2}^{\flat}  (u_1 ,u_2 , w_3).
\end{eqnarray*}
 
  \item Thanks the two previous points, we obtain $(\ref{zoub2})$ with
 \begin{eqnarray*}
 K_1^{\flat,2} (x_d ,\rond{v}^0 , v^{0,\flat} , \rond{w}^0  , w^{0,\flat}  ,\tilde{\mathcal{W}}^0 ,
  \frac{x_d}{\eps} \tilde{\mathcal{W}}^0 ) = 
\frac{x_d}{\eps} \tilde{\mathcal{W}}^0 \{  v^{0,\flat} . G_{1,1}^{\flat} (\rond{v}^0 ,\rond{w}^0 , x_d v^{0,\flat} ,x_d w^{0,\flat} ,\tilde{\mathcal{W}}^0) 
\\ + w^{0,\flat}   G_{1,2}^{\flat}   (\rond{v}^0  ,\rond{w}^0 , x_d  v^{0,\flat}  ,  x_d  w^{0,\flat}  , \tilde{\mathcal{W}}^0 )\} .
 \end{eqnarray*}
Because $G_{1,1}^{\flat}$ and $G_{1,2}^{\flat}$ are $C^\infty$, the function $ K_1^{\flat,2}$ is also $C^\infty$.

\end{enumerate}
\end{proof}

 We denote by 
\begin{eqnarray*}
 K_1^{\flat}   (\eps, x_d ,\rond{v}^0 , v^{0,\flat} , V^\eps , \eps V^\eps , \rond{w}^0 , w^{0,\flat}   , \tilde{\mathcal{W}}^0 
 ,\frac{x_d}{\eps} \tilde{\mathcal{W}}^0 , W^{\eps,\flat} ,\eps W^{\eps,\flat} ) := 
 \\  K_1^{\flat,2} (x_d ,\rond{v}^0 , v^{0,\flat}  , \rond{w}^0   , w^{0,\flat}   , \tilde{\mathcal{W}}^0 , \frac{x_d}{\eps} \tilde{\mathcal{W}}^0  ) 
 +  K_1^{\flat,1} (\eps , v^0 , V^\eps , \eps V^\eps , \mathcal{W}^0 ,  W^{\eps,\flat} ,\eps W^{\eps,\flat} )  
  \end{eqnarray*}
 Because  the first term does not involve $V^\eps$ and $W^{\eps,\flat}$ and
  the function $ K_1^{\flat,1}$ is affine with respect to its third variable $V^\eps$ and affine
with respect to its sixth variable $W^{\eps,\flat}$ with $\X_{\mathtt{v}^0} \,  v^0$ as factor in the leading coefficient, 
 the function  $ K_1^{\flat} $ is  affine with respect to its fifth variable $V^\eps$ and affine with respect to its eleventh argument 
 $W^{\eps,\flat}$  with
 $\X_{\mathtt{v}^0} \,  v^0$ as factor in the leading coefficient.
 Note that  as the profile $ \tilde{W}_0$ is rapidly decreasing in $X$, 
 the  profile $X \tilde{W}_0 $  is also rapidly decreasing. 
 Combine $(\ref{zoub})$, $(\ref{N})$ and $(\ref{zoub2})$ to find  $(\ref{sing})$.

\begin{step}[Use the ground state properties]

We are going to prove that the terms  $ K_1 (u^0  )$ and  $K_1 (\rond{v}^0  , \rond{w}^0 + \tilde{\mathcal{W}}^0 ) -
K_1 (\rond{u}^0  )$ are equal to $0$.

\end{step}

First because the ground state $u^0$ satisfies $(\ref{aller})$,  
\begin{eqnarray*}
 K_1 (u^0  )  =   S^\star  (u^0  )  (\X_{\mathtt{v}^0} \,    + \M^\star (u^0 ,\partial_x ))   v^0 =0 .
\end{eqnarray*}
Moreover, referring to  $(\ref{K1})$ and thanks to $(\ref{cond})$, we obtain:
\begin{eqnarray*}
 K_1 (\rond{v}^0 , \rond{w}^0 + \tilde{\mathcal{W}}^0) - K_1 (\rond{u}^0 )
= (  S^\star   (v^0 , \rond{w}^0  + \tilde{\mathcal{W}}^0) -   S^\star ( v^0 , \rond{w}^0 ) ) \X_{ \rond{\mathtt{v}}^{0}}  \rond{v}^{0} 
= 0  .
\end{eqnarray*}
\end{proof}

\section{WKB expansions}
\label{secform}

\subsection{Formal solutions}

We look for formal solutions  $(u^\eps )_\eps $ of $(\ref{aller})$ of the  form $ u^\eps =(v^0 + \eps V^\eps , W^\eps )$ 
where  $U^\eps := (V^\eps , W^\eps )$ is an expansion  
\begin{eqnarray}
\label{pk}
\sum_{n \geqslant 0} \eps^n \, \mathcal{U}^n (t,x, \frac{x_d}{\eps})
\end{eqnarray}
 where each  $ \mathcal{U}^n $ belongs to  $\P (T) $. 
 We rewrite  $ \mathcal{U}^n $ as  $ \mathcal{U}^n = ( \mathcal{V}^n ,  \mathcal{W}^n )$ 
 and we suppose  $  \un{\mathcal{W}}^0 := w^0$. 
 Let us explain what is meant by \textit{formal} \textit{solutions}.
 Plugging the expansion $(\ref{pk})$ into the system, using Taylor expansions and ordering the terms in powers of $\eps$, 
  we get a formal expansion in power series of $\eps$:
\begin{eqnarray}
\label{resu}
\sum_{n \geqslant -1} \eps^n \, \Phi^n (t,x, \frac{x_d}{\eps})
\end{eqnarray}
 where the  $(\Phi^n )_{n \geqslant -1}$ are in  $\P (T) $.
 We say that $(u^\eps )_\eps $ is a formal solution when all the resulting $\Phi^n $ are identically zero.
 The following theorem states that the system has formal solutions and that we can prescribe arbitrary initial values to the $(P_0
 \, \mathcal{U}^j |_{t=0} )_{j \in \N } $. 
  We introduce    $\Nrond_{init}  :=  H^\infty  (\R^d_+  , \mathcal{S}(\R_+))$ 
and the profile space  
\begin{eqnarray*}
 \P_{init} := \{ \mathcal{U} \in L^2 (\R^d_+ \times \R_+ ) \ \mathrm{such} \  \mathrm{that} \
 \mathrm{there} \  \mathrm{exist} \ \un{ \mathcal{U}} \in H^\infty (\R^d_+ ), \ \tilde{ \mathcal{U}} \in  \Nrond_{init} 
\\ \mathrm{such} \  \mathrm{that} \ \mathrm{for} \  \mathrm{all} \
  (x,X) \in \R^d_+ \times \R_+ ,\quad  \mathcal{U}  (x,X) = \un{ \mathcal{U}} (x) + \tilde{ \mathcal{U}} (x,X)  \}.
 \end{eqnarray*}

\begin{theo}
\label{try}

Assume that some profiles $(\mathcal{U}^j_{init} )_{j \in \N }$ such that $ (Id - P_0) \mathcal{U}^j_{init} =0$  and
$\mathcal{W}^0_{init} = w^0 |_{t=0}$ are given, 
then there exists a formal solution $(u^\eps) $  of $(\ref{aller})$ on $(0,T_0 )$ with some profiles  $(\mathcal{U}^j )_{j \in \N }$
in $\P (T_0 )$ such that $\un{\mathcal{W}}^0 := w^0 $ and that for all $j \in \N$, $P_0 \, \mathcal{U}^j |_{t=0}= \mathcal{U}^j_{init} $.
Moreover the profile $\mathcal{U}^0$ is polarized in the sense
 that $(Id - P_0) \tilde{\mathcal{U}} =0$ for all $(t,x) \in (0,T_0 ) \times \Omega$.

\end{theo}

Notice that Theorem $\ref{try}$ gives the existence of a formal solution till $T_0$.
At first sight, this may seems strange because Theorem $\ref{try}$ deals with large boundary layers and in some similar settings
large boundary layers have in general a nonlinear behavior (see for example paper \cite{MZ} of G.M\'etivier and K.Zumbrun about large
viscous boundary layers for noncharacteristic nonlinear hyperbolic problems and paper \cite{Gr1} of E.Grenier about large velocity
boundary layers for the Euler equations). 
It is also different from the result on entropy waves of C.Cheverry, O.Gu\`es and G.M\'etivier  \cite{cgm2}. 
In \cite{cgm2}, Theorem $3.4$ gives the local existence of formal solutions with large amplitude oscillations on the entropy.
More precisely Theorem $3.4$ of \cite{cgm2} gives the existence of the profile $\mathcal{U}^0$ for small time and we see on the
system $(3.36)$ of \cite{cgm2} that this is not because of a lack of analysis but because of an actual nonlinear effect.
At the opposite in the proof of Theorem $\ref{try}$ we will construct a formal solution thanks to linear profile problems.
This shows that a system can reveals additional transparency properties when looking at boundary layers instead of high frequency
oscillations.
The transparency property of \cite{moi3} supports this 
remark.\footnote{This phenomenon can also be observed when looking at boundary layers of "strong amplitude", following the terminology of
\cite{cgm1}.
Instead of nonlinear modulation equation as in \cite{cgm2}, we obtain linear profile equation.}

Since the profile $\mathcal{U}^0$ is polarized, we get $\tilde{\mathtt{V}}_\mathfrak{n}^0 = \tilde{\mathtt{P}}^0 = 0$ (see section
$\ref{setting}$ for the definition of $P_0$) as we have
mentioned it in the introduction. 
This means that the amplitude of normal velocity and pressure boundary layers is weaker than
the amplitude of tangential boundary layer. 
The fact that the profile  $\tilde{\mathtt{V}}_\mathfrak{n}^0 $ vanishes is a consequence of the conditions $(\ref{cond})$ and more
precisely of the fact that $\X_{\mathtt{v}^0} \,  \mathtt{p}^0 =0$
for $(t,x) \in (0,T_0) \times \Gamma$. Without this last condition, it is still possible to construct some nontrivial formal
solutions of the form $(\ref{pk})$ but with a nonvanishing profile $\tilde{\mathtt{V}}_\mathfrak{n}^0 $ (cf. section $\ref{coming}$). 
However such formal solutions can be unstable. 
We do not consider such formal solutions in this paper.
\\
\\ Next subsection is devoted to the proof of Theorem $\ref{try}$.
We outline here some key tools of the proof.
First we use the reduction of the system of section $3$.
We use in a essential way Proposition $\ref{sev}$.
For example because $K_1$ contains a factor $\eps$, the term $K^\eps$ is not singular with respect to $\eps$.
This is crucial because we look for formal solutions without boundary layers for $v^\eps$ but only with large boundary layers for
$V^\eps$.
In order to find the expansion  $(\ref{resu})$, we use several Taylor expansions with respect to $\eps$ and to $x_d$. 
The underlying idea is to obtain an expansion $(\ref{resu})$ inside which the profiles  $(\mathcal{U}^j )_{j \in \N }$ appear at the
highest possible order. 
We notably use the fact that $ \mathtt{v}_d^0$ and $\X_{\mathtt{v}^0} \, v^0$  vanish on the boundary and the fact that $K_1^\flat$ is affine with respect
to its eleventh argument with $\X_{\mathtt{v}^0} \, 
 v^0$ as factor in the leading coefficient (cf. Proposition $\ref{sev}$ and  Lemma  $\ref{psg}$).     
In order to find some  profiles  $(\mathcal{U}^j )_{j \in \N }$ for which the profiles  $(\Phi^j )_{j \in \N }$ 
of the resulting expansion  $(\ref{resu})$ are identically zero, 
we follow the strategy of \cite{Gu95}, \cite{moi}, \cite{moi2}, \cite{moi3} and exhibit a
sequence of profile problem. 
Theorem $\ref{casca}$ gives the existence of profiles $(\mathcal{U}^j )_{j \in \N }$ solving this sequence of problem.
Theorem $\ref{try}$ will be proved as a consequence of theorem $\ref{casca}$.
One key point in the proof of Theorem $\ref{casca}$ is that profile problems are linear, thanks to transparency properties.
One manifestation of this property is that $K_1^\flat$ is affine with respect
to its fifth argument. The profiles problems are linearized Euler equations for the regular part of the profiles and linear totally
characteristic transport equations for the boundary layer part of the profiles.

\subsection{Proof of Theorem $\ref{try}$}

We will proceed in two steps. 
First we explicit the expansion $(\ref{resu})$.
 Then we look for profiles  $(\mathcal{U}^j )_{j \in \N }$ for which the profiles  $(\Phi^j )_{j \in \N }$ 
of the resulting expansion  $(\ref{resu})$ are identically zero solving a sequence of profile problems.

\subsubsection{Expliciting the expansion $(\ref{resu})$} 

We introduce the operator $\H^0$ which maps a function $U \in H^\infty ((0,T) \times \R^d_+ )$ to the function 
$\H^0 (U, W^{\flat} ,  \tilde{\mathcal{W}}^0 ,  \frac{x_d}{\eps} \tilde{\mathcal{W}}^0 , \partial ) U \in H^\infty ((0,T) \times
\R^d_+ )$ where 
 \begin{eqnarray*}
  \H^0 (U, W^{\flat} ,  \tilde{\mathcal{W}}^0 ,  \frac{x_d}{\eps} \tilde{\mathcal{W}}^0 , \partial ) U:= 
  \mathfrak{L} (v^0 ,W, \partial ) U
 +  \begin{bmatrix}
 K^0 (V, \tilde{\mathcal{W}}^0 , \frac{x_d}{\eps} \tilde{\mathcal{W}}^0 ,    W^{\flat} )  
  \\ 0
 \end{bmatrix}
  \end{eqnarray*}
  with 
 \begin{eqnarray*}
 K^0 (V, \tilde{\mathcal{W}}_0 , \frac{x_d}{\eps} \tilde{\mathcal{W}}^0 , W^{\flat} ) := 
 K_{1}^{\flat} ( 0, x_d ,\rond{v}^0 , v^{0,\flat} , V , 0 , \rond{w}^0 , \rond{w}^{0,\flat}  ,  \tilde{\mathcal{W}}^0 , \frac{x_d}{\eps} \tilde{\mathcal{W}}^0 ,  
 W^{\flat} , 0) .
 \end{eqnarray*}
 Remember that the term $ K_{1}^{\flat}  $ is defined in Proposition $\ref{sev}$.
 As in $(\ref{singul})$, we denote
 $\tilde{\mathcal{W}}^0$ instead of $\tilde{\mathcal{W}}^0  (t,x,\frac{x_d}{\eps})$.
 Thanks to Prop $\ref{sev}$, the function $K^0 (V, \tilde{W} , X \tilde{W} ,W)$ is affine with respect
  to its first argument $V$ and affine with
 respect to its fourth argument $W$ with $\X_{\mathtt{v}^0} \,   v^0$ as factor in the leading coefficient. 

\begin{prop}

\begin{enumerate}[$(i)$]

\item We obtain 
 \begin{eqnarray*}
\un{\Phi}^{-1} =0,  \ \un{\Phi}^{0} = \H^0  (\un{  \mathcal{U}}^0 ,0,0 , \un{\mathcal{W}}^1 ,\partial )  \un{ \mathcal{U}^0 }, 
  \end{eqnarray*}
 and for $j \geqslant 1$,
 \begin{eqnarray*}
\un{\Phi}^{j} =  \H^0 (\un{  \mathcal{U}}^0 , 0,0,\un{ \mathcal{W}}^1 , \partial ) \un{  \mathcal{U}}^j -  \un{\mathcal{Q}}^j .
\end{eqnarray*}
The terms $\un{\mathcal{Q}}^j$ depend on  $(t,x)$, on the profiles $( \un{\mathcal{U}}^k )_{k  <j}$ and on their derivatives and affinely
on the profiles  $\un{\mathcal{U}}^{j}$ and $\un{\mathcal{W}}^{j+1}$. 
In particular, $ \un{\mathcal{Q}}^1_w = -  \un{\mathtt{V}}^0 . \nabla \un{\mathcal{W}}^0 $.

\item We obtain 
\begin{eqnarray*}
 \tilde{\Phi}^{-1} &=& \L_d \, \partial_X \tilde{\mathcal{U}}^0
\\  \tilde{\Phi}^{0} &=& \L_d \, \partial_X \tilde{\mathcal{U}}^1 
 + \H' (\mathcal{U}^0  , \partial ,  \partial_X ) \tilde{ \mathcal{U}}^0 ,
\\  \tilde{\Phi}^{j} &=& \L_d \, \partial_X \tilde{\mathcal{U}}^{j+1} 
+  \H' ( \mathcal{U}^0  ,  \partial ,  \partial_X ) \tilde{ \mathcal{U}}^j  
 -  \tilde{ \mathcal{Q}}^j \quad \mathrm{for} \ j \geqslant 1
\end{eqnarray*}
where

\begin{enumerate}[$\bullet$]

\item  the operator $ \H'  $ 
 maps a function $\mathcal{U} \in  \P(T) $ to
   $ \H' (\mathcal{U} ,\partial , \partial_X ) \tilde{\mathcal{U}} \in  \mathcal{N}(T) $ where
\begin{eqnarray*}
 \H' (\mathcal{U} ,\partial , \partial_X ) \tilde{\mathcal{U}} :=
 (X  \mathtt{v}^{0,\flat}_d +  \rond{\underline{\mathtt{V}}}_d +  \tilde{\mathtt{V}}_d ) S (v^0 ,\mathcal{W} )  \partial_X  \tilde{\mathcal{U}} 
 + \mathfrak{L} (v^0 , \mathcal{W} ,\partial)  \tilde{\mathcal{U} }
\\  + 
(S (v^0 ,\mathcal{W} ) - S (v^0 , \un{ \mathcal{W} })) \X_{\mathtt{v}^0}  \, \un{ \mathcal{U}}
+ 
\begin{bmatrix}
  K^0_I (  \mathcal{V} , \tilde{\mathcal{W}} , X \tilde{\mathcal{W}} )
\\ 0
\end{bmatrix} ,
 \end{eqnarray*}
where $K^0_I ( V , W , X W )$ is a $C^\infty$ function affine with respect to
its first argument $V$ such that   $K^0_I (  \mathcal{V}^0 , \tilde{\mathcal{W}}^0 , X \tilde{\mathcal{W}}^0 ) \in 
\mathcal{N}(T) $.

\item the terms $ \tilde{ \mathcal{Q}}^j $  depend of  $(t,x)$, of the profiles $( \mathcal{U}^k )_{k  <j}$, 
of their traces  and of their derivatives and
affinely of the profiles $\mathcal{U}^j $ and   $ \un{\mathcal{W}}^{j+1} $.

\end{enumerate}

\end{enumerate}
\end{prop}
\begin{proof}[Proof of $(i)$]
We introduce the operator $\H^\eps$ by the formula
\begin{eqnarray*}
 \H^\eps (U,\partial ) U:=  \mathfrak{L}( v^0 + \eps V ,W,\partial)   U + K^\eps .
\end{eqnarray*}
Remember that the term $K^\eps$ is defined in Proposition $\ref{Re1}$.
Assume that  $W $ is of the form $ W^{\eps} = \mathcal{W}^0 +\eps W^{\flat} $ where $ \mathcal{W}^0 \in \P $ 
with  $\un{\mathcal{W}}^0 = w^0$.
By a first order Taylor expansion in respect to $\eps$, using Proposition $\ref{sev}$, we see that  $\H^\eps $ can be put under the form:
\begin{eqnarray*}
 \H^\eps (U,\partial ) U = \H^0 (U, W^{\flat} , \tilde{\mathcal{W}}^0 ,  \frac{x_d}{\eps} \tilde{\mathcal{W}}^0 ,  \partial ) U
 + \sum_{j=1}^d B_{j}^{\eps} (v^0 , U) \eps \partial_j U + \eps \mathbf{M}^\eps (t,x,U ). 
 \end{eqnarray*}
 where the $B_{j}^{\eps}$ and $\mathbf{M}^\eps$ are $(d+2) \times (d+2) $ matrices depending in a $C^\infty$ way on their arguments
including $\eps$ up to $\eps=0$. They are  of the form 
 \begin{eqnarray*}
 B_{j}^{\eps} := 
\begin{bmatrix}
 B_{j}^{\star,\eps} & 0
\\ 0 & \mathtt{V}_j
 \end{bmatrix}
 , \quad
 \mathbf{M}^\eps :=
 \begin{bmatrix}
  \mathbf{M}^{\star,\eps} & 0
\\ 0 & 0
 \end{bmatrix}.
 \end{eqnarray*}
 Then we easily complete the proof of $(i)$.
\end{proof}
\begin{proof}[Proof of $(ii)$]
 We now look at the  $\tilde{\Phi}^{j}$.
Because the terms $\tilde{\mathcal{U}}^{j} (t,x,\frac{x_d}{\eps})$ contain a factor $\eps^{-1}$ in their argument, the normal derivative $\partial_d$ plays a
crucial role. 
Therefore, we define  
\begin{eqnarray}
\label{fd}
\X'_\mathtt{v} :=\X_\mathtt{v} - \mathtt{v}_d \partial_d \  \mathrm{and} \  \L'   (u, \partial )= \L  ( \partial) + S(u) \X'_\mathtt{v} .
\end{eqnarray}
Write $U := (V,W)$ to find that 
\begin{eqnarray}
 \L'   (u, \partial ) U = 
 \begin{bmatrix}
  \L^{\star,'}   (u, \partial ) V 
\\ \X'_\mathtt{v}  W
\end{bmatrix}
\quad \mathrm{where} \
 \L^{\star,'}   (u, \partial ) := \L^\star  ( \partial) + S^\star (u) \X'_\mathtt{v} .
\end{eqnarray}
Then
\begin{eqnarray*}
 \H^\eps (U,\partial ) U  &=&  (\mathtt{v}^0_d + \eps \mathtt{V}_d ) S^\eps (v^0 + \eps V ,W)  \partial_d U 
 + \L' (v^0 + \eps V ,W ,\partial) U   + K^\eps  ,
\\ \H^0 (U, W^{\flat}, \tilde{\mathcal{W}}^0 ,  \frac{x_d}{\eps} \tilde{\mathcal{W}}^0 ,  \partial ) U   &=&  
\mathtt{v}^0_d \, \S (v^0 ,W)  \partial_d U + \L' (v^0 , W,\partial) U  
 + K^0   (V, \tilde{\mathcal{W}}^0 ,  \frac{x_d}{\eps} \tilde{\mathcal{W}}^0 , W^{\flat}). 
\end{eqnarray*}
A key point when calculating the  $\tilde{\Phi}^{j}$ is that when a boundary layer profile  is in factor, 
we can do formally $"x_d =0"$. The
underlying idea is that for a $C^\infty$ scalar function $\Phi (x_d )$ and $\tilde{\mathcal{U}} \in  \mathcal{N}(T)$
using a Taylor expansion yields the existence of a $C^\infty$ scalar function $\Phi^\flat (x_d )$
 such that $\Phi (x_d ) = \Phi (0) + x_d \Phi^\flat (x_d ) $. Then
\begin{equation*}
\Phi (x_d ) \tilde{\mathcal{U}} (t,x,\frac{x_d}{\eps}) = \Phi (0) \tilde{\mathcal{U}} (t,x,\frac{x_d}{\eps}) 
+ \eps \Phi^\flat (x_d ) (X \tilde{\mathcal{U}} (t,x,X)) |_{X=\frac{x_d}{\eps}} .
\end{equation*}
Because $X  \tilde{\mathcal{U}} (t,x,X) $ still is in $\mathcal{N}(T)$, the second term can be put with the term of
order $\eps$. 
If $\Phi (0 )=0$, then you obtain that $\tilde{\mathcal{U}}$ does not appear in the resulting expansion at the order $0$.
We are going to apply successively this idea to the terms  
\begin{eqnarray*}
 \mathtt{v}^0_d \, S (v_0  ,W) \partial_d  \mathcal{U},  \ B_{d}^{\eps} (v_0 , \mathcal{U}) \eps \partial_d \mathcal{U} \ 
\mathrm{and} \ K^0 (\mathcal{V}^0 ,\tilde{\mathcal{W}}^0 ,X \tilde{\mathcal{W}}^0 ,\mathcal{W}^1) . 
\end{eqnarray*}
\begin{enumerate}[$\bullet$]

\item  Because  $\rond{\mathtt{v}}^0_{d} =0$, we have  $\mathtt{v}^0_d  = x_d  \mathtt{v}^{0,\flat}_d$.
Applying $ \mathtt{v}^0_d \, S (v^0  ,W) \partial_d$ to $ \tilde{\mathcal{U}}(t,x,\frac{x_d}{\eps})$ where 
  $\tilde{\mathcal{U}} \in \mathcal{N}(T) $ 
yields the term 
\begin{eqnarray*}
(\mathtt{v}^{0,\flat}_d \, S^0 (v^0  ,W) X \partial_X \tilde{\mathcal{U}}
 + \mathtt{v}^{0} \,  S^0 (v^0  ,W) \partial_d \tilde{\mathcal{U}} )|_{X=\frac{x_d}{\eps}} .
 \end{eqnarray*}
\item  We now look at the term  $B_{d}^{\eps} (v^0 ,\mathcal{U} )$.

\begin{lem}

There exist some matrices $B_{d}^{\eps,\flat}$ $C^\infty$ with respect to its arguments, including $\eps$ up to $0$, such that for
all $ \mathcal{U} \in \P(T)$, 
\begin{eqnarray}
\label{frg}
B_{d}^{\eps} (v^0 ,\mathcal{U} )  = \mathtt{V}_d \, S (v^0 ,\mathcal{W} ) 
+ \eps B_{d}^{\eps,\flat} (x_d , X, \rond{v}^0 , v^{0,\flat} ,\rond{\un{\mathcal{U}}} , \un{\mathcal{U}}^\flat ,\tilde{\mathcal{U}} ).
\end{eqnarray}
Moreover, the matrices  $B_{d}^{\eps,\flat}$ are affine with respect to $X$.

\end{lem}
\begin{proof}
We will proceed in three steps.
\\
\begin{enumerate}[$1$.]

\item  First  by a Taylor expansion there exist some matrices $B_{d}^{\eps,\flat,1}$, $C^\infty$ with respect to its arguments, 
including $\eps$ up to $0$, such that for all $\mathcal{U} \in \P(T)$, 
\begin{eqnarray*}
B_{d}^{\eps} (v^0 , \mathcal{U}) = B_{d}^{0} (v^0 , \mathcal{U}) 
+ \eps B_{d}^{\eps,\flat} (v^0 , \mathcal{U}  ),
\end{eqnarray*}
 where
\begin{eqnarray*}
 B_{d}^{0} (v^0 ,  \mathcal{U}) :=  \mathtt{V}_d \, S (v^0 ,\mathcal{W} ) + \mathtt{v}^0_d \, DS (v^0 ,\mathcal{W} ).( \mathcal{V},0) .
\end{eqnarray*}
\item  By a Taylor expansion, there exist some matrices $B_{d}^{0,\flat,1}$ and $B_{d}^{0,\flat,2}$,
 $C^\infty$ with respect to their arguments 
 \begin{eqnarray*}
 (v_1 , \mathcal{U}_1 ,v_2 , \mathcal{U}_2 ) \in \R^{d+1} \times \R^{d+2} \times \R^{d+1} \times \R^{d+2} ,
 \end{eqnarray*}
 such that for all $(v_1 , \mathcal{U}_1 ,v_2 ,\mathcal{U}_2 )$, 
 \begin{eqnarray}
 \label{tyh}
\quad  B_{d}^{0} (v_1 + v_2 , \mathcal{U}_1 + \mathcal{U}_2  ) =
 B_{d}^{0} ( v_1 , \mathcal{U}_1  ) + v_2 . B_{d}^{0,\flat,1} (v_1 ,\mathcal{U}_1 ,v_2 ,\mathcal{U}_2 )
\\  \nonumber + U_2 . B_{d}^{0,\flat,2} (v_1 ,\mathcal{U}_1 ,v_2 ,\mathcal{U}_2 ) .
 \end{eqnarray}
 We write  $ \mathtt{v}^0 = \rond{\mathtt{v}}^0 + x_d \, \mathtt{v}^{0,\flat}$ and for all  $\mathcal{U} \in \P (T)$, 
 $\un{\mathcal{U}} = \rond{\un{\mathcal{U}}} + x_d \, \un{\mathcal{U}}^\flat $.
 We apply  $(\ref{tyh})$ with 
 \begin{eqnarray*}
 (v_1 ,\mathcal{U}_1 ,v_2 ,\mathcal{U}_2  )= (\rond{v}^0 , \rond{\un{\mathcal{U}}}+ \tilde{\mathcal{U}} , x_d \, v^{0,\flat} ,  x_d \, \un{\mathcal{U}}^\flat ) ,
 \end{eqnarray*}
   and obtain for all $\mathcal{U} \in \P (T)$,
\begin{eqnarray*}
 B_d^0 (v^0 , \mathcal{U}) = B_d^0 (\rond{v}^0 ,\rond{\un{\mathcal{U}}} + \tilde{\mathcal{U}}) 
 +  x_d \, B_{d}^{0,\flat,2} (x_d ,\rond{v}^0 , v^{0,\flat} ,\rond{\un{\mathcal{U}}} , \un{\mathcal{U}}^\flat ,\tilde{\mathcal{U}} ),
\end{eqnarray*}
with 
\begin{eqnarray*}
 B_{d}^{0, \flat,2} (x_d ,\rond{v}^0 , v^{0,\flat} ,\rond{\un{\mathcal{U}}} , \un{U}^\flat ,\tilde{\mathcal{U}} )  &:= &
 v^{0,\flat}  . B_{d,1}^{0,\flat} (\rond{v}^0 ,\rond{\un{\mathcal{U}}} + \tilde{\mathcal{U}},  x_d \, v^{0,\flat} ,  x_d \, \un{U}^\flat )
\\ && + \un{U}^\flat . B_{d,2}^{0,\flat} (\rond{v}^0 ,\rond{\un{\mathcal{U}}} + \tilde{\mathcal{U}},  x_d \, v^{0,\flat} ,  x_d \, \un{U}^\flat ) .
 \end{eqnarray*}
\item   Moreover, for all $\mathcal{U} \in \P  (T) $, as $ \rond{\mathtt{v}}^0_d  =0$, 
 we have $ B_d^0 (\rond{v}^0 ,\mathcal{U}) =  \mathtt{V}_d S^0 (v^0 ,W)$.
As a consequence, we get $(\ref{frg})$ with 
\begin{eqnarray}
\label{NN}
 B_{d}^{\eps,\flat} (x_d , X, \rond{v}^0 , v^{0,\flat} ,\rond{\un{\mathcal{U}}} , \un{U}^\flat ,\tilde{\mathcal{U}} ):= B_{d}^{\eps,\flat} (v^0 , \mathcal{U}  )
 + X  B_{d}^{0, \flat,2} (x_d ,\rond{v}^0 , v^{0,\flat} ,\rond{\un{\mathcal{U}}} , \un{\mathcal{U}}^\flat ,\tilde{\mathcal{U}} ) 
\ \end{eqnarray}
which is affine with respect to $X$.

\end{enumerate}
\end{proof}

Write $\mathtt{V}_d = \rond{\underline{\mathtt{V}}}_d + x_d \mathtt{\underline{V}}_d^{\flat}$ and define
\begin{eqnarray*}
 B_{d}^{\eps,\flat,2} (x_d , X, \rond{v}^0 , v^{0,\flat} ,\rond{\un{\mathcal{U}}} , \un{\mathcal{U}}^\flat ,\tilde{\mathcal{U}} )
=  B_{d}^{\eps,\flat} (x_d , X, \rond{v}^0 , v^{0,\flat} ,\rond{\un{\mathcal{U}}} , \un{\mathcal{U}}^\flat ,\tilde{\mathcal{U}} )
+  \underline{\mathtt{V}}_d^{\flat} S (v^0 , \mathcal{W}) X  
\end{eqnarray*}
to see that by applying $B_{d}^{\eps} (v^0 ,\mathcal{U}) \eps \partial_d$ to $\tilde{\mathcal{U}}(t,x,\frac{x_d}{\eps})$ where 
  $\tilde{\mathcal{U}} \in \mathcal{N}(T) $,
 we get the term 
\begin{eqnarray*}
( (\rond{\underline{\mathtt{V}}}_d +  \tilde{\mathtt{V}}_d )  S (v^0 , \mathcal{W}) \partial_X \tilde{\mathcal{U}} 
+ \eps \{B_{d}^{\eps,\flat,2} (x_d , X, \rond{v}^0 , v^{0,\flat} ,\rond{\un{\mathcal{U}}} , \un{\mathcal{U}}^\flat ,\tilde{\mathcal{U}} ).\partial_X 
+  B_{d}^{\eps} (v^0 ,\mathcal{U}) \partial_d  \} \tilde{\mathcal{U}} ) |_{X=\frac{x_d}{\eps}}  .
\end{eqnarray*}
Moreover, since the matrices $B_{d}^{\eps,\flat}$ are affine with respect to $X$ (see $(\ref{NN})$), 
the two terms of the sum are in $\mathcal{N}(T)$.
\\
\item At first sight, the term 
 \begin{eqnarray*}
 T := K^0 (\mathcal{V}^0 ,\tilde{\mathcal{W}}^0 , X \tilde{\mathcal{W}}^0, \mathcal{W}^1 )
  -  K^0 ( \un{\mathcal{V}}^0 ,0 , 0,  \un{\mathcal{W}}^1 )
 \end{eqnarray*}
  must appear in 
$\tilde{\Phi}^{0}$ but we are going to redevelop this term.

\begin{lem}
\label{psg}

There exist two  $C^\infty$ functions $K^0_I$ and $K^0_{II}$  affine with respect to their  first argument $\mathcal{V}^0$ such that 
 \begin{eqnarray}
 \label{ty}
T  &:=& 
  K^0_I (  \mathcal{V}^0 , \tilde{\mathcal{W}}^0  , X \tilde{\mathcal{W}}^0 ) 
 + \eps K^0_{II} (\mathcal{V}^0 ,\tilde{\mathcal{W}}^0 , X \tilde{\mathcal{W}}^0 , \tilde{\mathcal{W}}^1 , X \tilde{\mathcal{W}}^1). 
\end{eqnarray}
\end{lem}

Moreover, the functions $K^0_I (  \mathcal{V}^0 , \tilde{\mathcal{W}}^0  , X \tilde{\mathcal{W}}^0 )$ and $ K^0_{II} (\mathcal{V}^0 ,\tilde{\mathcal{W}}^0 , X
\tilde{\mathcal{W}}^0 , \tilde{\mathcal{W}}^1 , X \tilde{\mathcal{W}}^1)$ are in $\mathcal{N}(T)$.
\begin{proof}
We will proceed in six steps.
\\
\begin{enumerate}[$1$.]
\item
We begin to split the term $T$ into three parts:  $T= T_1 +  T_2 +  T_3$ where
 \begin{eqnarray*}
 T_1 &:=& K^0 (\mathcal{V}^0 ,\tilde{\mathcal{W}}^0 , X \tilde{\mathcal{W}}^0, \mathcal{W}^1 ) - K^0 (\mathcal{V}^0
 ,\tilde{\mathcal{W}}^0 , X \tilde{\mathcal{W}}^0, \un{\mathcal{W}}^1 ),
\\  T_2 &:=& K^0 (\mathcal{V}^0 ,\tilde{\mathcal{W}}^0 , X \tilde{\mathcal{W}}^0, \un{\mathcal{W}}^1 ) 
-  K^0 (\mathcal{V}^0 ,0 , 0, \un{\mathcal{W}}^1 ),
\\  T_3 &:=&  K^0 (\mathcal{V}^0 ,0 , 0, \un{\mathcal{W}}^1 ) - K^0 ( \un{\mathcal{V}}^0 ,0 , 0, \un{\mathcal{W}}^1 ).
\end{eqnarray*}
\item  Because the term $ K^0 (V,\tilde{\mathcal{W}} ,  X \tilde{\mathcal{W}} , W)$ is affine with respect to its 
 first argument $V$ and affine  with respect to its fourth argument $W$ with $\X_{\mathtt{v}^0} \,   v^0$ as factor 
 in the leading coefficient, there exist some $C^\infty$ functions
 $K^0_1$ and  $K^0_2$  affine with respect to their first argument such that for all $(V,\tilde{\mathcal{W}},W) \in \R^{d+1} \times
 \R\times \R$,
 \begin{eqnarray}
 \label{L}
 K^0 (V ,\tilde{\mathcal{W}} , X \tilde{\mathcal{W}} , W ) &:=& K^0_1 (V ,\tilde{\mathcal{W}}, X \tilde{\mathcal{W}}  )
  + \X_{\mathtt{v}^0} \,   v^0 . K^0_2 (V ,\tilde{\mathcal{W}}, X \tilde{\mathcal{W}} ) .W .
   \end{eqnarray}
 \item To deal with the term $T_1$, we use twice the equality $(\ref{L})$, first with
  $(V ,\tilde{\mathcal{W}}, W)=( \mathcal{V}^0 ,\tilde{\mathcal{W}}^0 , \mathcal{W}^1 )$ and then with 
 $(V ,\tilde{\mathcal{W}}, W)=( \mathcal{V}^0 ,\tilde{\mathcal{W}}^0 ,  \un{\mathcal{W}}^1 )$ to find 
 \begin{eqnarray}
  \label{La}
 T_1 = \X_{\mathtt{v}^0} \,   v^0 .  K^0_2 (V ,\tilde{\mathcal{W}}, X \tilde{\mathcal{W}} ) .\tilde{\mathcal{W}}^1
 \end{eqnarray}
  \item To deal with the term $T_2$, we use twice the equality $(\ref{L})$, first with
  $(V ,\tilde{\mathcal{W}}, W)= ( \mathcal{V}^0 ,\tilde{\mathcal{W}}^0 ,  \un{\mathcal{W}}^1 )$ and then with 
 $(V ,\tilde{\mathcal{W}}, W)=( \mathcal{V}^0 ,0 , \un{ \mathcal{W}}^1 )$  to find 
 \begin{eqnarray}
  \label{Lb}
    T_2 &=&  K^0_1 ( \mathcal{V}^0 ,\tilde{\mathcal{W}}^0 ,  X \tilde{\mathcal{W}}^0 ) -  K^0_1 ( \mathcal{V}^0 ,0,0)  
\\ \nonumber &&+  \X_{\mathtt{v}^0} \,   v^0 . ( K^0_2 (  \mathcal{V}^0 ,\tilde{\mathcal{W}}^0 , X \tilde{\mathcal{W}}^0 ) 
   -   K^0_2 (  \mathcal{V}^0 , 0, 0 ))
   .\un{\mathcal{W}}^1 .
 \end{eqnarray}
 \item To deal with the term $T_3$,  we use twice the equality $(\ref{L})$, first with
  $(V ,\tilde{\mathcal{W}}, W)= (\mathcal{V}^0 ,0 ,  \un{\mathcal{W}}^1 )$ and then with 
 $(V ,\tilde{\mathcal{W}}, W)=( \un{\mathcal{V}}^0 ,0 ,  \un{\mathcal{W}}^1 )$
 \begin{eqnarray}
  \label{Lc}
  T_3 =  K^0_1 ( \mathcal{V}^0 ,0 , 0 ) -  K^0_1 (  \un{\mathcal{V}}^0 ,0,0) + 
 \X_{\mathtt{v}^0} \,   v^0 . ( K^0_2 (  \mathcal{V}^0 ,0 ,0) -   K^0_2 ( \un{\mathcal{V}}^0 , 0,0 ))  .\un{\mathcal{W}}^1 .
  \end{eqnarray}
   \item
 Using $(\ref{La})$, $(\ref{Lb})$, $(\ref{Lc})$, we get $(\ref{ty})$ with
 \begin{eqnarray*}
   K^0_I (  \mathcal{V}^0 , \tilde{\mathcal{W}}^0 ,  X \tilde{\mathcal{W}}^0 ) &:=&  K^0_1 ( \mathcal{V}^0 ,\tilde{\mathcal{W}}^0 ,  X \tilde{\mathcal{W}}^0 ) - 
    K^0_1 (  \un{\mathcal{V}}^0 ,0,0)
 \\  K^0_{II} (\mathcal{V}^0 ,\tilde{\mathcal{W}}^0 , X \tilde{\mathcal{W}}^0 , \un{\mathcal{W}}^1 , X \tilde{\mathcal{W}}^1) &:=&
    \frac{ \X_{\mathtt{v}^0} \,   v^0 }{x_d} X \{  K^0_2 (  \mathcal{V}^0 ,\tilde{\mathcal{W}}^0 ,  X \tilde{\mathcal{W}}^0 ). \tilde{\mathcal{W}}^1
 \\ &&+ ( K^0_2 (  \mathcal{V}^0 ,\tilde{\mathcal{W}}^0 ,  X\tilde{\mathcal{W}}^0 ) - K^0_2 (  \un{\mathcal{V}}^0 ,0,0)).
 \un{\mathcal{W}}^1 \}.
 \end{eqnarray*}

 We have used that $x_d = \eps X$.
 Because the term $\X_{\mathtt{v}^0} \,  v^0$ vanishes on the boundary, it contains a factor $x_d$. 
Therefore the functions $ K^0_I$ and $K^0_{II}$
 are $C^\infty$, affine with respect to their first argument. 
 Moreover, the functions 
 \begin{eqnarray*}
 K^0_I (  \mathcal{V}^0 , \tilde{\mathcal{W}}^0 ,  X \tilde{\mathcal{W}}^0 ) \quad \mathrm{and} \quad   K^0_{II} (\mathcal{V}^0 ,\tilde{\mathcal{W}}^0 , X
\tilde{\mathcal{W}}^0 , \tilde{\mathcal{W}}^1 , X \tilde{\mathcal{W}}^1)
 \end{eqnarray*}
 are in $\mathcal{N}(T)$.

 \end{enumerate}
\end{proof}

As a consequence, the function  $K^0_{II} (\mathcal{V}^0 ,\tilde{\mathcal{W}}^0 , X \tilde{\mathcal{W}}^0 , \un{\mathcal{W}}^1 , X
\tilde{\mathcal{W}}^1)$ can be put in $\tilde{\Phi}_{1}$.
Therefore the boundary layer profile $\mathcal{W}^1$ does not appear
in the contribution of the term  $T$  at the order $0$.
\\
\end{enumerate}
 Finally,  we get  the result and complete the proof of $(i)$.                                                                                                                            
\end{proof}

Write $ \mathcal{U} := (\mathcal{V} ,\mathcal{W} )$ to find that 
\begin{eqnarray}
\label{spil}
 \H' (X, \mathcal{U}  ,\partial , \partial_X ) \tilde{\mathcal{U}} :=
 \begin{bmatrix}
  \H^{\star,'} ( \mathcal{U}   ,\partial , \partial_X )  \tilde{\mathcal{V}}
 \\  (X  \mathtt{v}^{0,\flat}_d + \rond{\un{\mathtt{V}}}_d +  \tilde{\mathtt{V}}_d  ) \partial_X  \tilde{\mathcal{W}} 
 + \X_{\mathtt{v}^0}  \, \tilde{\mathcal{W}}
 \end{bmatrix}.
 \end{eqnarray}
 where
  \begin{eqnarray*}
 \H^{\star,'} (X,\mathcal{U} ,\partial , \partial_X )  \tilde{\mathcal{V}} :=
 (X \mathtt{v}^{0,\flat}_d + \rond{\un{\mathtt{V}}}_d +  \tilde{\mathtt{V}}_d  ) S^\star (v^0 , \mathcal{W})  \partial_X \tilde{\mathcal{V}} 
+  \mathfrak{L}^{\star} (v^0 , \mathcal{W} ,\partial)  \tilde{\mathcal{V}} 
\\ + (S^\star  (v^0 ,\mathcal{W} ) - S^\star  (v^0 , \un{ \mathcal{W} })) \X_{\mathtt{v}^0}  \, \un{ \mathcal{V}}
 +  K^0_I (  \mathcal{V} , \tilde{\mathcal{W}} , X  \tilde{\mathcal{W}}).
 \end{eqnarray*}
 
 \subsubsection{A sequence of profile problems}
 \label{tabl}

We split, for $j \geqslant 0$,  $ \un{\Phi}^{j}$  into
  \begin{eqnarray*}
  \un{\Phi}^{j} := (  \un{\Phi}^{j}_{v} ,  \un{\Phi}^{j}_{w} ) \in \R^{d+1} \times  \R .
  \end{eqnarray*}
   We define the problem 
\begin{eqnarray*}
(\mathcal{S}^j (T)): \quad
\left\{
\begin{array}{ccc}
&(Id -P_0 ) \til{\Phi}^{j-1} = P_0   \, \til{\Phi}^{j} = \un{\Phi}^{j}_{v} =  \un{\Phi}^{j+1}_{w}   =0
   \quad &\mathrm{where} \,  (t,x)  \in   (0,T) \times \Omega
\\    &\mathtt{V}^j_d =  0   \quad &\mathrm{where} \,  (t,x)   \in  (0,T) \times  \Gamma
\end{array}
\right.
\end{eqnarray*}
We illustrate our  strategy with the following table
\\
\\ \begin{tabular}{ccccc}
                    &  $(\mathcal{S}^{0})$                              & $(\mathcal{S}^{1})$                              & $(\mathcal{S}^{2})$           &   
		                 ...
\\                  &                                         &                                        &                            
        &  
\\    $\Phi^{-1}$   & $(Id - P_0 )    \tilde{\Phi}^{-1} $     &                                        &                            
        & 
\\   $\Phi^{0}$     & $ P_0  \, \tilde{\Phi}^0 ,\un{\Phi}^{0}_{v} $ & $(Id - P_0 ) \tilde{\Phi}^0$        &                            
        & 
\\    $\Phi_{1 }$   &   $ \un{\Phi}^{1}_{w} $                   & $ P_0  \, \tilde{\Phi}^1 ,\un{\Phi}^{1} $  &  $(Id - P_0 )   
\tilde{\Phi}^1 $   & 
\\ $\Phi^{2 }$      &                                         &    $ \un{\Phi}^{2}_{w} $                 &   $ P_0  \, \tilde{\Phi}^2
,\un{\Phi}^{2}_{v} $   &                        
\\   $\Phi^{3 }$ &         &           &  $ \un{\Phi}^{3}_{w} $  &                    
\\ ... & & & &
\end{tabular}
\\
\\ Each element of the left column contain  the  sum of the  corresponding line. When we solve successively 
the  problems $(\mathcal{S}^j (T))$ the unknown profiles are for all $j \geqslant 0$, $\mathcal{V}^j $, $ \un{\mathcal{W}}^{j+1} $
and $\tilde{\mathcal{W}}^{j}$ when solving the  problem  $(\mathcal{S}^j (T))$. 
  Notice that $ \un{\Phi}^{0}_{w} =  \X_{\mathtt{v}^0} \,  w^0  =0 $ because the ground state $(v^0 ,\mathcal{W}^0 )$ 
 is solution of the Euler system.  
  Our goal in this section is to prove the following theorem:

 \begin{theo}
 \label{casca}
 
 There exist profiles $(\mathcal{U}^j )_{j \geqslant 0}$ in $\mathcal{P} (T_0 )$  which verify $\un{\mathcal{W}}^0 :=
 w^0 $, $(\mathcal{S}^j (T_0))_{j \geqslant 0} $ and
  $ P_0 \mathcal{U}^j |_{t=0} =  \mathcal{U}^j_{init}$ for all $j \geqslant 0$.
 Moreover the profile $\mathcal{U}^0$ is polarized in the sense
 that $(Id - P_0) \tilde{\mathcal{U}} =0$ for all $(t,x) \in (0,T_0 ) \times \Omega$.

\end{theo}

Theorem $\ref{try}$ is a consequence of theorem $\ref{casca}$.
\begin{proof}
We will proceed in two steps first studying the problem $(\mathcal{S}^0 (T_0))$ then for all $ j \geqslant 1$, 
the problems  $(\mathcal{S}^j (T_0) )$.
The problem $(\mathcal{S}^0 (T))$ splits into two sub-problems 
\begin{eqnarray}
\label{1P}
\left\{
\begin{array}{ccc}
& \un{\Phi}^{0}_{v} = \un{\Phi}^{1}_{w} = 0  &\quad \mathrm{where} \,  (t,x)  \in  (0,T) \times  \Omega
\\  & \un{\mathtt{v}}^0_d =  0              &\quad \mathrm{where} \,  (t,x)   \in   (0,T) \times \Gamma
\end{array}
\right.
\end{eqnarray}
and   
\begin{eqnarray}
\label{2P}
(Id - P_0 ) \tilde{\Phi}^{-1} =  P_0  \, \tilde{\Phi}^0 =0  \quad \mathrm{where} \,  (t,x)  \in   (0,T) \times \Omega
\end{eqnarray}
Let us begin to look at $(\ref{1P})$: we have  
\begin{eqnarray*}
 \un{\Phi}^{0}_{v} &=&  \mathfrak{L}^\star (u^0 ,\partial)    \un{\mathcal{V}}^0 +  K^0  ( \un{\mathcal{V}}^0 ,0 ,0, \un{\mathcal{W}}^1 ) ,
 \\  \un{\Phi}^{1}_{w} &=&   \X_{\mathtt{v}^0} \,  \un{\mathcal{W}}^1 +  \mathtt{\un{V}}^0 . \nabla \un{\mathcal{W}}^0 .
\end{eqnarray*}
In  $(\ref{1P})$, the determinations of $\un{\mathcal{V}}^0$ and $\un{\mathcal{W}}^1$ are coupled. 
The existence of solutions is given by the following lemma.

 \begin{lem}
\label{noncon}

There exist some profiles $ (\un{\mathcal{V}}^0 , \un{\mathcal{W}}^1)$ 
in $H^\infty ( (0,T_0 ) \times \Omega)$ solution of:
\begin{eqnarray}
\label{elair} \qquad 
\left\{
\begin{array}{cc}
\mathfrak{L} (u^0 ,\partial) 
    \begin{bmatrix}
     \un{\mathcal{V}}^0
      \\ \un{\mathcal{W}}^1
   \end{bmatrix}
   +  
    \begin{bmatrix}
   K^0 ( \un{\mathcal{V}}^0  ,0 ,0, \un{\mathcal{W}}^1 ) 
    \\ \un{\mathtt{V}}^0 . \nabla \un{\mathcal{W}}^0 
   \end{bmatrix}
    = 0  
      &\quad \mathrm{where} \,  (t,x) \in      (0,T_0 ) \times \Omega
\\  \underline{\mathtt{V}}^0_d =0  &\quad \mathrm{where} \,  (t,x) \in  (0,T_0 ) \times \Gamma
\end{array}
\right.
\end{eqnarray}
 such that 
 \begin{eqnarray*}
  P_0^\star \, \un{\mathcal{V}}^0 |_{t=0} =  \un{\mathcal{V}}^0_{init},  \ 
 \un{\mathcal{W}}^1 |_{t=0} = \un{\mathcal{W}}^1_{init} .
 \end{eqnarray*}
 \end{lem}
\begin{proof}
We will proceed in three steps.

\begin{enumerate}[$(i)$]

\item 
It is important to see that these equations are linear because the function  $K^0 $ is affine 
with respect to its first and third argument.

\item 
Some compatibility conditions on the corner $\{t=x_d =0 \}$ are necessary to the existence of solutions in $H^\infty$.
The existence of initial compatible initial data $(\hat{\un{\mathcal{V}}}^0_{init} , \hat{\un{\mathcal{W}}}^1_{init} )$ are given in 
Proposition $3.10$ of  \cite{moi}. Moreover Proposition $3.10$ of \cite{moi} insures that we can prescribe arbitrarily  $P_0 \, 
\hat{\un{\mathcal{V}}}^0_{init} $ and $ \hat{\un{\mathcal{W}}}^1_{init} $. It is therefore possible to choose $P_0 \, 
\hat{\un{\mathcal{V}}}^0_{init} := \un{\mathcal{V}}^0_{init}$ and $\hat{\un{\mathcal{W}}}^1_{init} = \un{\mathcal{W}}^1_{init}$.

\item
Then  Th\'eor\`eme $3.10$ of  \cite{moi} yields the existence of a solution $(\un{\mathcal{V}}^0 , \un{\mathcal{W}}^1)$ in $H^\infty
((0,T_0 ) \times \Omega)$.

\end{enumerate}
\end{proof}
\begin{rem}

Referring to Proposition $\ref{sev}$, we can  compute the term  $ K^0 ( \un{\mathcal{V}}^0  ,0 ,0, \un{\mathcal{W}}^1 ) $ and see
that the Prob. $(\ref{elair})$ can be restated as a boundary value problem for the linearized Euler system. 
To see this by another way, plug $v^\eps = v^0 + \eps \un{\mathcal{V}}^0$ and 
$w^\eps = w^0 + \eps \un{\mathcal{W}}^0$ in $(\ref{po})$.

\end{rem}

 Let us attack  $(\ref{2P})$. 
 To solve it we use the following result:

 \begin{lem}
 \label{lem2}
 
 There is a unique  profile  $\til{U}^0$  in a function of  $\Nrond (T_0)$ solution of
\begin{eqnarray}
\label{lem2a}
&(Id-P_0 )  \til{\mathcal{U}}^0 =0 \quad &\mathrm{where} \,  (t,x,X) \in (0,T_0) \times \Omega_ \times \R_+ ,
\\ \label{lem2b}
& P_0  \,  \H' (X ,v_0  , \mathcal{W}^0 , \partial , \partial_X ) \til{\mathcal{U}}^0 = 0 
  \quad &\mathrm{where} \,  (t,x,X) \in  (0,T_0) \times \Omega  \times \R_+ ,
  \\ \label{lem2c}
&  P_0  \, \til{\mathcal{U}}^0 = \til{\mathcal{U}}^0_{init}  \quad &\mathrm{where} \,  t=0 .
\end{eqnarray}
\end{lem}
\begin{proof}
If we solve directly the problem, because the system is \`a-priori nonlinear, we can only claim the local existence of a solution $
\til{\mathcal{U}}^0$. A more acute analysis is possible using transparency properties.
In order to exploit them, we are going to proceed in two steps first looking for $ \til{\mathcal{W}}^0$ 
then looking  for the tangential velocity $\mathtt{V}_\mathfrak{t}$.

\begin{enumerate}[$(i)$]

\item Referring to $(\ref{spil})$, and thanks to the equations $(\ref{lem2a})$ and $(\ref{lem2b})$,  
 we find  for $\til{\mathcal{W}}^0$ the equation
 \begin{eqnarray*}
(X  \mathtt{v}^{0,\flat}_d + \rond{\un{\mathtt{V}}}_d +  \tilde{\mathtt{V}}_d ) \partial_X  \,  \tilde{\mathcal{W}}^0 
 + \X_{\mathtt{v}^0}  \,  \tilde{\mathcal{W}}^0 =0 . 
\end{eqnarray*}
Thanks to $(\ref{lem2a})$ and  $(\ref{elair})$,  we have 
 \begin{eqnarray}
 \label{prune}
 \til{\mathtt{V}}^0_d =0 \quad \mathrm{and} \  \rond{\un{\mathtt{V}}}_d^0 =0 .
 \end{eqnarray}
 Therefore we get
 \begin{eqnarray}
 \label{clem}  \quad
 \mathtt{v}^{0,\flat}_d  \, X \partial_X  \tilde{\mathcal{W}}^0 
 + \X_{\mathtt{v}^0}  \,  \tilde{\mathcal{W}}^0 =0 . 
\end{eqnarray}
This is a linear transport equation for which the  boundary  is totally characteristic.
As a consequence, there is none compatibility condition at the corner $\{ t= x_d =0 \}$. 
Therefore there exists one (and only one) solution $\til{\mathcal{W}}^0$ in $\P(T_0 )$ of $(\ref{clem})$ such that 
$\til{\mathcal{W}}^0 |_{t=0} = \til{\mathcal{W}}^0_{init} $.

\item 
We now look for $P_0^\star \, \til{\mathcal{V}}^0$ i.e. for $\tilde{\mathtt{V}}^0_\mathfrak{t}$.
 Referring to $(\ref{spil})$, using once more $(\ref{prune})$, we get 
 \begin{eqnarray}
  \label{M1}
 \H^{\star,'} (X ,\mathcal{U}  ,\partial , \partial_X ) \tilde{\mathcal{V}}  = 
 \mathtt{v}^{0,\flat}_d \, S^\star (v^0 , \mathcal{W} )  X \partial_X \tilde{\mathcal{V}}  +  
 \mathtt{v}^{0}_d  \,S^\star  (v^0 , \mathcal{W})  \partial_d  \tilde{\mathcal{V}}
 + \L^{\star,'} (v^0 , \mathcal{W} ,\partial)  \tilde{\mathcal{V}} 
\\ \nonumber +(S^\star  (v^0 ,\mathcal{W}^0 ) - S^\star  (v^0 , \un{ \mathcal{W} })) \X_{\mathtt{v}^0} \, \un{ \mathcal{U}}
 +   K^0_I (  \mathcal{V} , \tilde{\mathcal{W}} , X \tilde{\mathcal{W}} ) .
 \end{eqnarray}
 Referring to section $(\ref{setting})$,  we get 
 \begin{eqnarray}
 \label{M3}
  P_0^{\star}  \, S^\star  (v^0 , \un{ \mathcal{W} }) P_0^{\star} = \rho( \mathtt{p}^0 , \mathcal{W}^0 ) Id_{d-1}
 \quad \mathrm{and} \  P_0^{\star}  \,  \L^{\star} P_0^{\star} = 0, 
  \end{eqnarray}
We denote by $K^0_{I,\mathfrak{t}}$ the $d-1$ first components of $K^0_{I}$.
 Thanks to  $(\ref{M1})$ and $(\ref{M3})$, we deduce from the equation $(\ref{lem2b})$ the following equation for the tangential
 velocity $\tilde{\mathtt{V}}^0_\mathfrak{t}$:
 \begin{eqnarray}
 \label{poli} \quad
\\ \nonumber  \rho( \mathtt{p}^0 ,\mathcal{W}^0)  (\mathtt{v}^{0,\flat}_d X \partial_X 
 + \X_{\mathtt{v}^0})  \, \tilde{\mathtt{V}}^0_\mathfrak{t}
  + K^0_{I,\mathfrak{t}} (  \mathcal{V}^0 , \tilde{\mathcal{W}}^0 , X \tilde{\mathcal{W}}^0 ) 
 =  (  \rho (\mathtt{p}^0 ,\mathcal{W}^0 ) 
-  \rho ( \mathtt{p}^0 , \un{ \mathcal{W}}^0)) \X_{\mathtt{v}^0}  \, \un{ \mathtt{V}}^0_\mathfrak{t}  .
 \end{eqnarray}
 We face a totally characteristic initial boundary value problem and once more time, there is none compatibility conditions. 
 The system is first order symmetric hyperbolic.
Because the function $K^0_{I,\mathfrak{t}}$ is affine with respect to its first argument, this equation $(\ref{poli})$ is linear.
Therefore there exists one (and only one) $\tilde{\mathcal{V}}_0 $ solution of  the equation $(\ref{poli})$ 
on  $(0,T_0) \times \Omega  \times \R_+ $ with $P_0^\star \, \tilde{\mathcal{V}}^0 |_{t=0} = \tilde{\mathcal{V}}^0_{init} $.

\end{enumerate}
\end{proof}

\begin{rem}

Referring to Proposition $(\ref{sev})$, we can compute 
 \begin{eqnarray*}
K^0_{I,\mathfrak{t}} (  \mathcal{V}^0 , \tilde{\mathcal{W}}^0 , X \tilde{\mathcal{W}}^0 ) :=  
(  \rho (\mathtt{p}^0 ,\mathcal{W}^0 ) -  \rho ( \mathtt{p}^0 , \un{ \mathcal{W}}^0))
 \frac{\X_{\mathtt{v}^0}  \,  \mathtt{v}^0_\mathfrak{t} }{x_d } X 
\\ + (  \rho (\mathtt{p}^0 ,\mathcal{W}^0 ) \mathtt{V}^0 - \rho (\mathtt{p}^0 , \un{\mathcal{W}}^0 ) \un{\mathtt{V}}^0 ).  \nabla 
\mathtt{v}^0_\mathfrak{t} .
 \end{eqnarray*}

\end{rem}

 For all $ j \geqslant 1$, the problem  $(\mathcal{S}^j (T_0) )$ splits into several problems. 
First, we solve:
\begin{eqnarray}
\label{depol}
(Id-P_0 )  \tilde{\Phi}^{j-1} =0 \quad \mathrm{where} \, (t,x,X) \in  (0,T_0) \times \Omega \times \R_+
\end{eqnarray}
then 
\begin{eqnarray}
\label{pol}  
 P_0  \,  \tilde{\Phi}^{j} =0   \quad \mathrm{where} \, (t,x,X) \in  (0,T_0) \times  \Omega \times \R_+
\end{eqnarray}
and 
 \begin{eqnarray}
\label{rft}
\left\{
\begin{array}{ccc}
&\un{\Phi}^{j}_{v}  =   \un{\Phi}^{j+1}_{w}  =   0 \quad &\mathrm{where} \,  (t,x) \in  (0,T_0) \times  \Omega
\\    &\un{\mathtt{V}}^j_d  =0   \quad &\mathrm{where} \,  (t,x) \in  (0,T_0) \times \Gamma
\end{array} 
\right. 
\end{eqnarray}
To do so, we will proceed in three steps.

\begin{enumerate}[$(i)$]

\item 
We solve  $(\ref{depol})$ and defining $(Id-P_0) \til{\mathcal{U}}^j$ as the unique solution in  $\Nrond (T)$ of
\begin{eqnarray*}
\L_d   \, \partial_X (Id-P_0) \tilde{\mathcal{U}}^{j}= 
- (Id-P_0 ) ( \H' (X, \mathcal{U}^0 ,   \partial ,  \partial_X ) \tilde{ \mathcal{U}}^{j-1} 
 -  \tilde{ \mathcal{Q}}^{j-1} ).
\end{eqnarray*}
We stress that $\partial_X $ is an automorphism of  $\Nrond (T)$.
Notice that we cannot prescribe arbitrarily $(Id-P_0) \tilde{\mathcal{U}}^{j}$ at $t=0$.

\item
The problem  $(\ref{rft})$ is equivalent to a linearized Euler problem as in Lemma $\ref{noncon}$ for
$(\un{\mathcal{V}}^{j},\un{\mathcal{W}}^{j+1})$.
We obtain easily the existence of profiles $ \un{\mathcal{V}}^j$, $\un{\mathcal{W}}^{j+1}$ 
 in $H^\infty ((0,T_0) \times \Omega)$ solution of this problem with 
  \begin{eqnarray*}
   P_0^\star  \,  \un{  \mathcal{V}}^j |_{t=0} =  \un{\mathcal{V}}^j_{init} , \  
 \un{  \mathcal{W}}^{j+1} |_{t=0} = \un{  \mathcal{W}}^{j+1}_{init} .
 \end{eqnarray*}    
\item We  define $P_0 \, \til{\mathcal{U}^j}$ as the unique solution in  $\Nrond (T_0)$ of 
the totally characteristic linear hyperbolic initial boundary value:
\begin{eqnarray*}
P_0  \, \H' (X, \mathcal{U}^0 , \partial, \partial_X )  P_0  \,   \tilde{ \mathcal{U}}^j 
 =  P_0 ( \tilde{ \mathcal{Q}}^{j}   
 - (Id-P_0 ) \H' (X, \mathcal{U}^0 ,  \partial , \partial_X ) \tilde{ \mathcal{Q}}^j )
\end{eqnarray*}    
with $ P_0  \, \tilde{ \mathcal{U}}^j |_{t=0} =  \tilde{ \mathcal{U}}^j_{init}$.
Thus $(\ref{pol})$ is proved. 

\end{enumerate}
\end{proof}

\section{Stability}
\label{stab}

In this section, we are interested in the existence and the propagation of exact solutions of $(\ref{aller})$ asymptotic to the
formal solutions constructed in the previous section. 
Thus we assume that formal solutions $(u^\eps_{formal})_\eps $ of
$(\ref{aller})$ on $(0,T_0 ) \times \Omega $ of the  form 
\begin{eqnarray*}
 u^\eps_{formal} =(v^0 + \eps V^\eps_{formal} , W^\eps_{formal} ) \quad
\mathrm{where} \   U^\eps_{formal} := (V^\eps_{formal} , W^\eps_{formal} )
\end{eqnarray*}    
 is an expansion  
\begin{eqnarray*}
\sum_{n \geqslant 0} \eps^n \, \mathcal{U}^n (t,x, \frac{x_d}{\eps})
\ \mathrm{with} \ 
\mathcal{U}^n = (\mathcal{V}^n , \mathcal{W}^n ) \in \P (T_0 ), \  \underline{\mathcal{W}}^0 = w^0 
\end{eqnarray*}    
are given.
We obtain approximate solutions $u^\eps_a = ( v^\eps_a =   v^0 + \eps V^\eps_a , W^\eps_a )$ 
of the system $(\ref{aller})$, choosing for $n \in \N$, 
\begin{eqnarray*}
V^\eps_a  (t,x) &:=& \sum_{k=0}^n \eps^{k} \, \mathcal{V}^k (t,x, \frac{x_d}{\eps}),
\\ W^\eps_a  (t,x) &:=& \sum_{k=0}^n \eps^{k} \, \mathcal{W}^k (t,x, \frac{x_d}{\eps}) + \eps^{n+1} \, \un{\mathcal{W}}^{n+1} (t,x) .
\end{eqnarray*}    
We  denote $Z_0 := \partial_t$, $Z_i := \partial_i$ for $1 \leq i \leq n-1$ and $Z_d := h(x_d ) \partial_d$, where $h$ is a
bounded and $C^\infty$ function on $\R$ such that $h(x_d) \neq 0$ for $x_d  \neq 0$, $h(x_d) =x_d$ when $0\leqslant x_d \leqslant 1$
and $h(x_d) = 1$ when $ x_d \geqslant 2$. 
The family $(Z_i )_{0 \leqslant i \leqslant d}$ generates the algebra of $C^\infty$ tangent vector fields to $\Gamma$.
For all $l \in \N$, we  denote by $Z^l$ the collection of the operators of the form  $Z_0^{\alpha_0} ... Z_d^{\alpha_d}$ where
$\alpha_0 ,... , \alpha_d $ are in $\N$ and satisfy $|\alpha|:= \alpha_0 +... + \alpha_d =l$. 
To simplify, we introduce the notations 
\begin{eqnarray*}
L^2 (T) :=  L^2 ((0,T) \times \R_+^d) \quad \mathrm{and} \quad  L^\infty (T)  :=  L^\infty ((0,T) \times \R_+^d ).
\end{eqnarray*}
 The family $(U^\eps_a )_{\eps \in ]0,1]}$ verifies the following estimates:
\begin{eqnarray}
\label{linf}
\mathrm{sup}_{0 < \eps \leqslant 1} || (\eps \partial_d)^k \, Z^l \, U^\eps_a ||_{L^\infty (T_0)  } < \infty .
\end{eqnarray}
For $m \in \N$ and $T >0$, we denote by $\mathbf{E}^m (T)$ the set 
\begin{eqnarray*}
\mathbf{E}^m (T) := \{   (u^\eps)_\eps   \in L^2 ( T ) /  \quad 
 \mathrm{sup}_{0 < \eps \leqslant 1}  
 \sum_{0 \leqslant 2k+l \leqslant m}  ||  (\eps \partial_d )^k \,  Z^l \,  u ||_{L^2 ( T)}  < \infty  \} .
\end{eqnarray*}
To be clear about our notations, let us stress that in the previous inequality, the sum  
\begin{eqnarray*}
\sum_{0 \leqslant 2k+l \leqslant m}  ||  (\eps \partial_d )^k  \, Z^l \,  u ||_{L^2 ( T)}
\end{eqnarray*}
 stands for
\begin{eqnarray*}
\sum_{0 \leqslant 2k+|\alpha| \leqslant m}  ||  (\eps \partial_d )^k  \, Z^\alpha \,  u ||_{L^2 ( T)}
\end{eqnarray*}
Let us explain why we use the set $\mathbf{E}^m (T)$. 
We begin with a brief review about smooth solutions of characteristic hyperbolic initial boundary value problem.
We referred to the work of O.Gu\`es \cite{Gu90}.\footnote{An analysis in $L^2$ was achieved by J.Rauch \cite{R}.}
First the boundary matrix of the system $(\ref{po2})$ is $A_d (u) := S(u) \mathtt{v}_d + \mathbf{L}_d $.
When $x_d =0$, $\mathtt{v}_d =0$ and the rank of $A_d (u) = \mathbf{L}_d $ is constant.
This suggests to use an extra-derivative namely $x_d \partial_d$ to the tangential derivatives $\partial_t, ...,\partial_{d-1}$.
This yields the notion of conormal regularity.
Then, handle normal derivatives for characteristic problem needs carefulness.
An example by \cite{mo} shows that in general there are not $H^s$ estimates (unlike the noncharacteristic case). 
We can extirpate $A_d \partial_d u$ from the equation but because the matrix $A_d$ is not invertible, 
this does not provide estimates
for the components in the kernel of $A_d$.
However these components satisfy a transport equation (cf. \cite{rare}, \cite{mefr}, \cite{Gu90})
with a source term which contains two conormal derivatives of all the
components.   
An iteration yields the idea to use one normal derivative for two conormal derivatives.
Thus in \cite{Gu90}, O.Gu\`es uses the spaces
\begin{eqnarray*}
E^m (T) := \{   u  \in L^2 ( T ) /  \quad 
 \sum_{0 \leqslant 2k+l \leqslant m}  || \partial_d^k  \, Z^l  \, u ||_{L^2 ( T )}  < \infty  \} .
\end{eqnarray*}
Here we face a singular perturbation problem. 
More precisely, we look at boundary layers which corresponds to variations
in $\frac{x_d}{\eps}$. 
That is why we introduce some more adapted sets $\mathbf{E}^m (T)$ with the derivatives $\eps \partial_d $ instead of $\partial_d $.
This idea of using some derivatives with $\eps$ in factor for some singular perturbations problems was used in
\cite{Gu93},
\cite{Gu92}, \cite{cgm2},... 
Here, this idea is applied to anisotropic Sobolev spaces.
Let us mention one technical point.
In \cite{Gu90}, O.Gu\`es uses a reduction of the system.
For the Euler system, this corresponds to choose the thermodynamic variables $p$, $\mathtt{v}$, $s$.

Looking at the table of subsection $\ref{tabl}$,  we see that  
\begin{eqnarray}
\label{E}
\H^\eps (U^\eps_a , \partial) U^\eps_a = \eps^{n+ \frac{1}{2}} 
 \begin{bmatrix}
 R^\eps_v
\\ \eps R^\eps_W
\end{bmatrix} 
\quad \mathrm{for} \ \mathrm{all} \ (t,x) \in (0,T_0) \times \Omega ,
\end{eqnarray}  
with
 \begin{eqnarray*}
 R^\eps_v &:=& (Id - P_0^* ) \tilde{\Phi}^n_v  + \eps \un{\Phi}^{n+1}_v + ...
\\ R^\eps_W &:=&  \tilde{\Phi}^n_w + \eps \un{\Phi}^{n+2}_w + ...
\end{eqnarray*}  
  Notice that if $\tilde{\Phi} \in \mathcal{N}(T)$ then the family $(\Phi^\eps )_{\eps}$ defined by
  \begin{eqnarray*}
 \Phi^\eps (t,x) := \eps^{-\frac{1}{2}} \,  \tilde{\Phi} (t,x,\frac{x_d}{\eps}) \quad \mathrm{for} \ \mathrm{all} \ (t,x) \in (0,T) \times \Omega ,
 \end{eqnarray*}  
 is in  $\mathbf{E}^m (T) $, for all $m \in \N$.
 Thus the family $( \eps^{-\frac{1}{2}} \,  R^\eps )_{\eps}$ defined by
 \begin{eqnarray*}
 R^\eps (t,x) :=  \begin{bmatrix}
\eps R^\eps_v
\\ R^\eps_{W}
\end{bmatrix}
\quad \mathrm{for} \ \mathrm{all} \ (t,x) \in (0,T) \times \Omega ,
\end{eqnarray*}  
 is in  $\in \mathbf{E}^m (T) $, for all $m \in \N$. 
 
 The system of equations $(\ref{E})$ is equivalent to the two equations 
  \begin{eqnarray}
  \label{E1} \quad &\mathfrak{L}^\star (u^\eps_a ,\partial)  V^\eps_a 
  + \frac{1}{\eps} K_1 (v^0 , \partial v^0 , u^\eps_a )
= \eps^{n+ \frac{1}{2}} R^\eps_v 
 \quad &\mathrm{for} \ \mathrm{all} \ (t,x) \in (0,T_0) \times \Omega,
  \\ \label{E2}
 & \X_{\mathtt{v}^\eps_a } W^\eps_a = \eps^{n+ \frac{3}{2}} R^\eps_W \quad &\mathrm{for} \ \mathrm{all} \ (t,x) \in (0,T_0) \times \Omega.
  \end{eqnarray}  
 Multiplying the equation $(\ref{E1})$ by $\eps$, using that $u^0$ satisfies the equation 
 \begin{eqnarray*}
 \mathfrak{L} (u^0 ,\partial)  (\partial_x ) v^0 =0 ,
 \end{eqnarray*}  
 we obtain that the family $(u^\eps_a )_\eps $ satisfy
\begin{eqnarray*}
\mathfrak{L} (u^\eps_a ,\partial) u^\eps_a 
=
\eps^{n+ \frac{1}{2}}
\begin{bmatrix}
\eps R^\eps_v
\\ \eps R^\eps_{W}                                                                                                                   
                                                      
\end{bmatrix} \quad \mathrm{when} \ x_d > 0 ,
\\ \mathtt{v}^\eps_{a,d} = 0  \quad \mathrm{when} \ x_d = 0 .
\end{eqnarray*}    
We look for solutions $(u^\eps )_\eps \in H^{\infty} $ of the problem $(\ref{aller})$ of the form 
\begin{eqnarray}
\label{fofo}
u^\eps := (v^\eps, W^\eps ) \ \mathrm{with} \ v^\eps = v^\eps_a + \eps^{M+1} \, 
V^\eps_R , \ W^\eps = W^\eps_a + \eps^{M} \, W^\eps_R .
\end{eqnarray}    
 We denote $U^\eps_R := (V^\eps_R  , W^\eps_R )$.
 We begin with a result of propagation.

\begin{theo}
\label{propa}

Let $m >\frac{d}{2}+5 $, $n \geqslant 1$, $M \in ]\frac{1}{2}, n+ \frac{1}{2} [$ and  $T \in ]0,T_0[$. 
We assume that we have a family of exact solutions $(u^\eps )_\eps \in H^{\infty} ((0,T) \times \Omega )$ of the
problem $(\ref{aller})$ of the form $(\ref{fofo})$ where the family $(U^\eps_R )_{\eps}$ is in the set $\mathbf{E}^m (T)$.
Then there is $\eps_0 \in ]0,1]$ such that for all $\eps \in ]0,\eps_0 ]$, the solution $u^\eps$ can be extended in a solution in
$H^\infty ((0,T_0 ) \times \Omega )$ of the problem $(\ref{aller})$ and is of the form $(\ref{fofo})$ where we have extended
$(U^\eps_R )_{\eps}$ in a family of $\mathbf{E}^m (T_0)$.

\end{theo}

The assumption $n+ \frac{1}{2} > M$ insures that the $u^\eps_a$ are accurate enough approximate solutions and the assumption $M >
\frac{1}{2}$ that the $u^\eps_a$ are close to $u^\eps$. 
Thanks to these assumptions,  Theorem $\ref{propa}$ claims that it is possible to extend the $u^\eps$, for $\eps \in ]0,\eps_0 ]$,
till $T_0$ i.e. till the ground state exists.
Notice that we can extend the result of Theorem $\ref{propa}$ to the more general case $n \in \N$, $M \geqslant \frac{1}{2}$.
To do so, consider  $\hat{n} \geqslant M - \frac{1}{2}$ and 
approximate solutions of order  $\hat{n}$:
\begin{eqnarray*}
\hat{u}^\eps_a = ( \hat{v}^\eps_a , \hat{w}^\eps_a ) =  ( v^0 + \eps \hat{V}^\eps_a , \hat{W}^\eps_a )
\end{eqnarray*}    
of the system $(\ref{aller})$, choosing
\begin{eqnarray*}
\hat{V}^\eps_a  (t,x) &:=&  \sum_{k=0}^n \eps^{k} \, \mathcal{V}^k (t,x, \frac{x_d}{\eps}),
\\ \hat{W}^\eps_a  (t,x) &:=& \sum_{k=0}^n \eps^{k} \, \mathcal{V}^k (t,x, \frac{x_d}{\eps}) + \eps^{n+1} \, \un{\mathcal{W}}^{n+1} (t,x) .
\end{eqnarray*} 
We define the family  $(\hat{U}^\eps_R )_{\eps}$ by 
\begin{eqnarray*}
\hat{U}^\eps_R := U^\eps_R + \eps^{-M} \, (U^\eps - \hat{U}^\eps_a ) \in  \mathbf{E}^m (T) .
\end{eqnarray*} 
On $(0,T)$, we get $U^\eps =\hat{U}^\eps_a + \eps^{M} \hat{U}^\eps_R $. 
Applying Theorem $\ref{propa}$, we obtain some extensions of the $(\hat{U}^\eps_R )_{\eps}$ in a family of $\mathbf{E}^m (T_0)$. 
Then we extend $( U^\eps_R  )_{\eps}$ in a family of $\mathbf{E}^m (T_0)$ , setting 
\begin{eqnarray*}
U^\eps_R = \hat{U}^\eps_R  - \eps^{-M} \, (U^\eps - \hat{U}^\eps_a ).
\end{eqnarray*} 
In this paper, we consider a ground state $u^0$ in $H^\infty ((0,T_0) \times \Omega )$  and formal solutions with $H^\infty$
regularity. It could also be possible to extend to ground states of high but finite regularity. 

We also give a result of existence.

\begin{theo}
 \label{di}

Let $m >\frac{d}{2}+4 $, $n \geqslant 1$, $M \in ]\frac{1}{2}, n+ \frac{1}{2} [$.
There exists $\eps_0 \in ]0,1]$ such that for all $\eps \in ]0,\eps_0 ]$, 
there exists a solution $u^\eps \in H^\infty ((0,T_0 ) \times \Omega )$ of
$(\ref{aller})$ of the form $(\ref{fofo})$ where $(U^\eps_R)_{\eps}$ is in $\mathbf{E}^m (T_0)$.

\end{theo}

It is important to understand that Theorem $\ref{propa}$ and $\ref{di}$ yields exact solutions till $T_0$. 
Our method is based on the sets $\mathbf{E}^m$ and on some estimates uniform with respect to $\eps$. 
Theorem $\ref{int}$ given in the introduction is a consequence of Theorem $\ref{try}$ and Theorem  $\ref{di}$. 

Because the problem $(\ref{aller})$ comes from a system of conservation laws and $\mathbf{\lambda}(u,\xi)$ 
is an eigenvalue with constant
multiplicity (cf. section $2$), it is possible to obtain $L^\infty$ estimates, even for $d \geqslant 1$. 
We refer to  papers \cite{m2}, \cite{sonique} of G.M\'etivier, paper \cite{RR} of J.Rauch and M.Reed and paper \cite{Gu90}.
Therefore we can weaken the regularity of the solution and prove a propagation result for some solutions admitting only
one normal derivative in $L^2$. 
We introduce the sets
\begin{eqnarray*}
\mathbf{A}^m (T) := \{   (u^\eps)_\eps   \in L^2 ( T) /  \quad 
 \mathrm{sup}_{0 < \eps \leqslant 1}  
 \sum_{0 \leqslant l \leqslant m}  ||  Z^l  \, u ||_{L^2 ( T )} 
 +  \sum_{0 \leqslant l \leqslant m-2}  ||  \eps \partial_d   \, Z^l  \, u ||_{L^2 ( T )} < \infty  \} .
\end{eqnarray*}
We will also use some norms built on $L^\infty$.
We denote by $Z^\eps $ the collection of the derivatives  $Z_0 $,...,$Z_d$ and $ \eps \partial_d $. 
Because the boundary is characteristic, we will need not only the Lipschitz norms but higher order $L^\infty$ control, 
as O.Gu\`es in
\cite{Gu90} and G.M\'etivier in \cite{sonique}.
We denote by $L^\infty (T)$ the space  $L^\infty (T) =L^\infty ((0,T) \times \R_+^d )$.
We introduce the norms
\begin{eqnarray*}
||u||_{k,T }  &:=& \sum_{0 \leqslant k \leqslant m} ||  Z^k u ||_{L^\infty (T )} ,
\\ ||u||^*_{\eps,T}  &:=& ||u||_{0,T} + ||Z^\eps  u||_{1,T} ,
\\  ||u||_{\eps,Lip,T}  &:=&  ||u||_{0,T} +   ||Z^\eps u ||_{0,T} .
\end{eqnarray*}
Remember that, by abuse of notation, we denote for example  $||Z^\eps u ||_{0,T}$
 for 
\begin{eqnarray*}
 \sum_{0 \leqslant i  \leqslant d} ||Z_i u ||_{0,T} + ||  \eps \partial_d u ||_{0,T} .
 \end{eqnarray*}
We introduce the sets
\begin{eqnarray*}
\mathbf{\Lambda}^m (T) := \{  (u^\eps)_\eps   \in L^\infty ( T)  /  \quad 
 \mathrm{sup}_{0 < \eps \leqslant 1} \, ||u||^*_{\eps,T} < \infty  \} . 
\end{eqnarray*}
It is also possible to tackle the limit case $M = \frac{1}{2} $, but we can prove the propagation only till $T_1 \in ]T,T_0]$.
We incorporate this limit case in the following Theorem.

\begin{theo}
\label{propa2}

Let $m > \frac{d}{2}+5$,  $M \geqslant \frac{1}{2} $, $T \in ]0,T_0[$. 
We assume that we have a family of exact solutions $(u^\eps )_\eps $ of the
problem $(\ref{aller})$ on $(0,T)$ of the form $(\ref{fofo})$ 
where $(U^\eps_R)_{\eps}$ is in $\mathbf{A}^m (T) \cap \mathbf{\Lambda}^m (T) $.
Then there exists $T_1 \in ]T,T_0]$ and there is $\eps_0 \in ]0,1]$ such that for all $\eps \in ]0,\eps_0 ]$,
 the solution $u^\eps$ can be extend in a solution on
$(0,T_1 )$ of the problem $(\ref{aller})$ and is of the form $(\ref{fofo})$ where we have extended
$(U^\eps_R)_{\eps}$ in a function of $\mathbf{A}^m (T_1) \cap \mathbf{\Lambda}^m (T_1)$.
Moreover if $M> \frac{1}{2} $, we can take $T_1 = T_0$.  

\end{theo}

It could be also possible to treat the limit case $M = \frac{1}{2} $ with the sets $\mathbf{E}^m (T_0)$ and incorporate a result of
propagation till $T_1 \in ]T,T_0]$ in Theorem $\ref{propa}$. We did not do so for sake of clarity
The rest of this section is devoted to the proof of Theorem $\ref{propa}$ and Theorem  $\ref{di}$.
The proof of Theorem  $\ref{di}$ needs carefulness about the existence of compatible initial data.      
Subsection  $\ref{bobo}$ is devoted to this question.
In subsection $\ref{Reduction}$ we perform a reduction in a problem for $(U^\eps_R )_{\eps}$ which are the real unknown.
We obtain (cf. Prop. $\ref{Reduction2}$) that $U^\eps_R$ satisfies a quasi-linear symmetric hyperbolic boundary value problem.
As for the originating Euler problem, the boundary is conservative and characteristic of constant multiplicity.
At first look, this system is singular with respect to  $\eps$ because of a factor $\eps^{-1}$ in the equation of $V^\eps_R$.
However a further analysis reveals that in fact the singular term contains $\frac{x_d}{\eps} W^\eps_R $ as factor.
This will be a key point in order to surmount the apparent singularity. 
In order obtain existence of  $(U^\eps_R )_\eps$ till $T_0$, we will use a family of iterative schemes. 
Thus we will supply in subsection $\ref{Linesti}$ linear estimates which are the core the proof. 
We will successively perform $L^2$ estimates, conormal estimates and normal estimates.
Several difficulties occur and are melt.
First, the boundary is characteristic. 
As we have explained it above, to tackle this problem we get inspired by the paper of O.Gu\`es \cite{Gu90}. 
We adapt the method of O.Gu\`es substituting the derivative $\eps \partial_d $ to the derivative $\partial_d $ 
in order to obtain uniform estimates. 
Moreover, we use estimates of $\frac{x_d}{\eps} W^\eps_R $ in order to surmount the apparent singularity.
In subsection $\ref{Itsch}$, 

The proof of Theorem $\ref{propa2}$ is not detailed as it obtained from Th\'eor\`eme $3$ of \cite{Gu90} in the same way Theorem
$\ref{propa}$ and Theorem  $\ref{di}$ are inspired from Th\'eor\`eme $1$ and $2$ of \cite{Gu90}.
The discussion about $M-\frac{1}{2}$ appears in the proof in subsection $\ref{Itsch}$ when we use the Sobolev embedding Lemma 
$\ref{sobo}$. Some minor modifications at this step allow to tackle the limit case $M =\frac{1}{2}$.

\subsection{Reduction}
\label{Reduction}

Because the $(U^\eps_R )_\eps$ are the real unknown, we begin with a reduction.
We will denote
\begin{eqnarray*}
W^{\eps,\flat}_a  (t,x) := \sum_{k=1}^n \eps^{k-1} \, \mathcal{W}^k (t,x, \frac{x_d}{\eps}) 
+ \eps^{n} \, \un{\mathcal{W}}^{n+1} (t,x) .
\end{eqnarray*}
\begin{prop}
\label{Reduction2}

There are some $(d+1) \times (d+1)$ matrices $J^\alpha$, 
 some functions $J^{\beta,1}$, $J^{\beta,2}$
 with values in $\R^{d+1}$, such that for all $\eps \in ]0,1]$,
 a function $u^\eps = (v^\eps, W^\eps )$ of the form $(\ref{fofo})$ verifies $(\ref{aller})$ 
 if and only if $U^\eps_R := (V^\eps_R  , W^\eps_R )$ verifies 
\begin{eqnarray}
\label{poi}
 \\ \nonumber (\mathfrak{L}^* (u^\eps ,\partial) + J^\alpha ).V^\eps_R
+   W^\eps_R .(J^{\beta,1} +  \frac{x_d}{\eps} .J^{\beta,2}) = - \eps^{n+\frac{1}{2}-M} R^\eps_v \quad &\mathrm{where} \, (t,x) \in (0,T) \times \Omega
\\ \label{poi2} \X_{\mathtt{v}^\eps }  \, W^\eps_R + \eps \mathtt{V}^\eps_R . \nabla W^\eps_a = -\eps^{n+\frac{3}{2}-M}  R^\eps_W 
&\quad \mathrm{where} \, (t,x) \in (0,T) \times \Omega
\\ \label{poi3} \mathtt{V}^\eps_{R,d} = 0  &\quad \mathrm{where} \, (t,x) \in (0,T) \times \Gamma .
\end{eqnarray}
Moreover, the matrices $J^\alpha$ and the functions  $J^{\beta,1}$, $J^{\beta,2}$ depend in a $C^\infty$ way of 
\begin{eqnarray*}
(\eps,V^\eps_a , Z_\eps V^\eps_a ,W^{\eps,\flat}_a, \eps^{M} \, U^\eps_R ). 
\end{eqnarray*}
\end{prop}
\begin{proof}
We only sketch the if part.
The converse is left to the reader.
Thus we assume that $(u^\eps)$ satisfies $(\ref{aller})$. 
According to Proposition $\ref{Re1}$, the function $U^\eps$ verifies $(\ref{Re2})$ 
which is equivalent to the two coupled equations $(\ref{Re3})$. 
\\ 
\\ We begin to show $(\ref{poi})$.
Because the family of functions $(u^\eps )_\eps$ is of the form $(\ref{fofo})$, for each $\epsilon$, 
the function $W^\eps$ is of the form $(\ref{qz})$ where $ \mathcal{W}^0 \in \P(T)$
 is such that $\un{\mathcal{W}}^0 = w^0$ and $ W^{\eps,\flat} =  W^{\eps,\flat}_a + \eps^{M-1}  \, W^{\eps}_R $.
Thus we can use Proposition $\ref{sev}$. 
We will denote 
\begin{eqnarray*}
\label{rr}
K_1^\flat (V^\eps, \eps V^\eps , W^{\eps,\flat} ,\eps W^{\eps,\flat}) 
\end{eqnarray*}
instead of
\begin{eqnarray*}
K_{1}^{\flat} (\eps, x_d ,\rond{v}^0 , v^{0,\flat} , V^\eps , \eps V^\eps , \rond{w}^0 , w^{0,\flat} ,  \tilde{\mathcal{W}}^0
 ,\frac{x_d}{\eps} \tilde{\mathcal{W}}^0 , W^{\eps,\flat} , \eps W^{\eps,\flat}) 
 \end{eqnarray*}
because these arguments are the ones which
are important at this step of the analysis.
Let us introduce our strategy. 
We would like to write $(\ref{rr})$
as a perturbation of 
\begin{eqnarray*}
K_1^\flat (V^\eps_a , \eps V^\eps_a , W^{\eps,\flat}_a ,\eps W^{\eps,\flat}_a )
\end{eqnarray*}
 of order $\eps^M$. 
But because $W^{\eps,\flat} := W^{\eps,\flat}_a + \eps^{M-1}  \, W^{\eps}_R $ a naive Taylor expansion fails to give the desired
result.
The trick lies in the way $K_1^\flat $ depends of $W^{\eps,\flat}$ (cf. Proposition $\ref{sev}$).
This allows to factorize by $\eps^M$ and $ \frac{x_d}{\eps}$.

\begin{lem}
\label{Redlem1}

There are  some $(d+1) \times (d+1)$ matrices $K_1^{\flat,\alpha}$ and  some $C^\infty$ 
functions $K_1^{\flat,\beta,1}$, $K_1^{\flat,\beta,2}$  with values in $\R^{d+1}$ such that 
\begin{eqnarray}
\label{KK}
\quad K_1^\flat (V^\eps, \eps V^\eps , W^{\eps,\flat} ,\eps W^{\eps,\flat}) &=&
 K_1^\flat (V^\eps_a , \eps V^\eps_a  , W^{\eps,\flat}_a  ,\eps W^{\eps,\flat}_a )
\\ \nonumber &&+ \eps^{M} \,   K_1^{\flat,\alpha} (\eps, V^\eps_a ,W^{\eps,\flat}_a , \eps^{M}  \, U^\eps_R) .V^\eps_R
\\ \nonumber &&+ \eps^{M} \,  W^\eps_R . K_1^{\flat,\beta,1} (\eps, V^\eps_a ,W^{\eps,\flat}_a , \eps^{M}  \, U^\eps_R) 
\\ \nonumber &&+ \eps^{M} \, \frac{x_d}{\eps} W^\eps_R . K_1^{\flat,\beta,2} (\eps, V^\eps_a , W^{\eps,\flat}_a ). 
\end{eqnarray}
\end{lem}
\begin{proof}
We will proceed in three steps. 

\begin{enumerate}

\item \label{pl1} By a first order Taylor development, there exist some $C^\infty$ $(d+1) \times (d+1)$ matrices $K_1^{\flat,\alpha}$ 
and some  $C^\infty$ functions $K_1^{\flat,\beta,1}$ with values in $\R^{d+1}$ such
that 
\begin{eqnarray*}
K_1^\flat (V^\eps, \eps V^\eps , W^{\eps,\flat} ,\eps W^{\eps,\flat}) &=&
 K_1^\flat (V^\eps_a , \eps V^\eps_a  , W^{\eps,\flat}  ,\eps W^{\eps,\flat}_a )
\\ &&+ \eps^{M} \, K_1^{\flat,\alpha} (\eps, V^\eps_a ,W^{\eps,\flat}_a , \eps^{M}  \, U^\eps_R) .V^\eps_R 
\\ &&+ \eps^{M} \, W^\eps_R . K_1^{\flat,\beta,1} (\eps, V^\eps_a  , W^{\eps,\flat}_a , \eps^{M}  \, U^\eps_R). 
\end{eqnarray*}
\item \label{pl2} According to Proposition $\ref{sev}$, the function 
\begin{eqnarray*}
K_1^\flat (V^\eps_a , \eps V^\eps_a  , W^{\eps,\flat} ,\eps
W^{\eps,\flat}_a)
\end{eqnarray*}
 is affine with respect to $  W^{\eps,\flat} $ with $\X_{\mathrm{v}^0} v^0$ as a factor in the leading term.
This means that  there is a $C^\infty$ function $K_1^{\flat,\eta}$ with values in $\R^{d+1}$ such that
\begin{eqnarray*}
K_1^\flat (V^\eps_a , \eps V^\eps_a  , W^{\eps,\flat} , \eps W^{\eps,\flat}_a) &=&
K_1^\flat (V^\eps_a , \eps V^\eps_a  , W^{\eps,\flat}_a , \eps W^{\eps,\flat}_a) 
\\ &&+ \eps^{M-1} \, \X_{\mathrm{v}^0}  \, v^0  .W^\eps_R . K_1^{\flat,\eta}
(\eps,V^\eps_a  ,  W^{\eps,\flat}_a ) .
\end{eqnarray*}
\item Thanks to the conditions $(\ref{cond})$, we get $ \X_{\mathrm{v}^0} v^0 =0$ when $x_d =0$.
We define
\begin{eqnarray*}
 K_1^{\flat,\beta,2} (\eps,V^\eps_a , W^{\eps,\flat}_a ) :=  \frac{1}{x_d} \X_{\mathrm{v}^0}  \, v^0 .K_1^{\flat,\eta}
(\eps,V^\eps_a  ,  W^{\eps,\flat}_a ) .
\end{eqnarray*}
The function $K_1^{\flat,\beta,2}$ is  $C^\infty$ and thanks to $(\ref{pl1})$ and $(\ref{pl2})$ 
we get $(\ref{KK})$.

\end{enumerate}
\end{proof}

We now look at the term $S^* (u^\eps ) \X_{\mathtt{v}^\eps} \,  V^\eps $.

\begin{lem}
\label{Redlem2}

There are some $(d+1) \times (d+1)$ matrices $S^{*,\alpha}$ and some functions  $S^{*,\beta}$ with values in
$\R^{d+1}$ such that 
\begin{eqnarray*}
S^* (u^\eps ) \X_{\mathrm{v}^\eps} V^\eps = 
S^* (u^\eps_a ) \X_{\mathrm{v}^\eps_a} V^\eps_a + \eps^{M} \, S^* (u^\eps )  \X_{\mathrm{v}^\eps}  \, V^\eps_R
            + \eps^{M} S^{*,\alpha} . V^\eps_R + \eps^{M} \, W^\eps_R . S^{*,\beta} .
\end{eqnarray*}
Moreover, the matrices $S^{*,\alpha}$ and $S^{*,\beta}$ depend in a $C^\infty$ way of 
\begin{eqnarray}
\label{RR}
(\eps, V^\eps_a , Z_\eps  V^\eps_a ,  \eps^{M} \, U^\eps_R )  .
 \end{eqnarray}
\end{lem}
\begin{proof}
We proceed in three steps.

\begin{enumerate}

\item We apply  $(\ref{zou})$  with $(v_1 , w_1 ,v_2 , w_2 ) = (u^\eps_a , \eps^{M} \,  ( \eps V^\eps_R ,  W^\eps_R  ))$ and obtain
\begin{eqnarray}
\label{tyj}
S^* (u^\eps ) = S^* (u^\eps_a ) &+& \eps^{M+1} \, V^\eps_R . S^{*,\flat}_1 (u^\eps_a , \eps^{M} \,  ( \eps V^\eps_R ,  W^\eps_R  )) 
\\    \nonumber                          &+&  \eps^{M} \, W^\eps_R . S^{*,\flat}_2 (u^\eps_a , \eps^{M} \, ( \eps V^\eps_R ,  W^\eps_R  )) .
\end{eqnarray}
\item We get 
\begin{eqnarray*}
 \X_{\mathrm{v}^\eps} V^\eps = \X_{\mathrm{v}^\eps_a} V^\eps_a + \eps^{M} \, \X_{\mathrm{v}^\eps} V^\eps_R 
 + \eps^{M} \, \mathtt{V}^\eps_R .(\eps \nabla  V^\eps_a  ).
\end{eqnarray*}
\item We introduce the $(d+1) \times (d+1)$ matrices $S^{*,\alpha}$ such that for all $V \in \R^{d+1}$,
\begin{eqnarray*}
S^{*,\alpha} V := \eps V . S^{*,\flat}_1 . ( \X_{\mathrm{v}^\eps_a}  \, V^\eps_a  
+  \eps^{M} \, \mathtt{V}^\eps_R . (\eps \nabla  V^\eps_a  ). 
\end{eqnarray*}
Notice that $\X_{\mathrm{v}^\eps_a} V^\eps_a $ is not singular and can be expressed thanks to $V^\eps_a $, $Z_\eps V^\eps_a$.
Thus the matrices $S^{*,\alpha} $ are $C^\infty$ with respect to $(\ref{RR})$.

We also introduce the functions 
 \begin{eqnarray*}
S^{*,\beta}  :=  S^{*,\flat}_2 . (\X_{\mathrm{v}^\eps_a}  \, V^\eps_a  +  \eps^{M} \, \mathtt{V}^\eps_R. (\eps \nabla  V^\eps_a  ))
\end{eqnarray*}
which take values in $\R^{d+1}$ and are $C^\infty$ with respect to $(\ref{RR})$.
Thanks to $1$ and $2$, we get $(\ref{tyj})$. 

\end{enumerate}
\end{proof}

We define  the matrices $(d+1) \times (d+1)$ matrices $J^\alpha$
\begin{eqnarray*}
J^\alpha (\eps,V^\eps_a , Z_\eps V^\eps_a , W^{\eps,\flat}_a , \eps^{M} \, U^\eps_R ) &:=& 
 K_1^{\flat,\alpha} (\eps, V^\eps_a ,W^{\eps,\flat}_a , \eps^{M} \, U^\eps_R) 
\\  &&+ S^{*,\alpha} (\eps, V^\eps_a , Z_\eps  V^\eps_a , \eps^{M} \, U^\eps_R  ) ,
\end{eqnarray*}
and, for $i \in \{1,2\}$, the $\R^{d+1}$-valued function
\begin{eqnarray*}
 J^{\beta,i} (\eps,V^\eps_a , Z_\eps V^\eps_a , W^{\eps,\flat}_a , \eps^{M} \, U^\eps_R ) &:=&
  K_1^{\flat,\beta,i} (\eps,V^\eps_a ,  W^{\eps,\flat}_a , \eps^{M} \, U^\eps_R ) 
\\ &&+ S^{*,\beta} (\eps, V^\eps_a , Z_\eps  V^\eps_a , \eps^{M} \, U^\eps_R ) .
 \end{eqnarray*}
Thanks to Lemma $\ref{Redlem1}$ and $\ref{Redlem2}$ and Equation $(\ref{E1})$, we show that $V^\eps_R$  satisfies $(\ref{poi})$.
\\
\\ As
\begin{eqnarray*}
 \X_{\mathrm{v}^\eps}  \, W^\eps = \X_{\mathrm{v}^\eps_a}  \, W^\eps_a + \eps^{M}  \, \X_{\mathrm{v}^\eps}  \, W^\eps_R 
 + \eps^{M}  \, \mathtt{V}^\eps_R .(\eps \nabla  W^\eps_a  ),
\end{eqnarray*}
 and using Equation $(\ref{E2})$,  we get $(\ref{poi2})$.
We easily show $(\ref{poi3})$ and complete the proof.
\end{proof}
In the formulas $(\ref{poi})-(\ref{poi2})-(\ref{poi3})$, there is a singular factor $\eps^{-1}$ 
which appears in the term $\frac{x_d}{\eps} W^\eps_R J^{\beta,2}$. 
One idea would be to try to obtain estimates for $W^\eps_R$ ponderated by some $\eps$. 
The difficulty lies in the fact that the equation $(\ref{poi3})$ of $W^\eps_R$ involves $V^\eps_R$ 
 in return by the term $\eps  \mathtt{V}^\eps_R  .\nabla  W^\eps_a $.
Fortunately, because of the special form of $\nabla  W^\eps_a $, we will see that it is possible to find good estimates of
$\frac{x_d}{\eps} W^\eps_R$. This allows to overcome the false singularity.

\subsection{Initial data}
\label{bobo}

In order to obtain smooth solutions of $(\ref{poi})-(\ref{poi2})-(\ref{poi3})$, some compatibility conditions for the initial data
$(U^\eps_{R,init})_{\eps \in ]0,1]}$ are necessary. 
At the order $0$, this reads 
\begin{eqnarray*}
\mathcal{R}^{0,\eps} : \  \mathtt{V}^\eps_{R,init,d} |_{x_d =0}= 0.
\end{eqnarray*}
We now explain what are the compatibility conditions at order $j \geqslant 1$. 
From the equation $(\ref{poi})$, we can extirpate $\partial_t \mathtt{V}^\eps_{R,init,d} $ in function of spatial derivatives and so by
restriction its trace $\partial_t \mathtt{V}^\eps_{R,init,d} |_{ t=x_d =0 }$ on the corner $\{  t=x_d =0 \}$.
More precisely, there exists a $C^\infty$ function $\H_1$ such that 
\begin{eqnarray*}
\partial_t \mathtt{V}^\eps |_{ t=x_d =0 } = \H_1 ((\partial_x^l U^\eps_{R,init} |_{x_d =0 } )_{l \leqslant 1} ) - \partial_t
R^\eps_{\mathtt{V}_d}  |_{t=x_d =0} .
\end{eqnarray*}
in fact, the function $\H_1$ also depends in a $C^\infty$ way of $t,x,v^0$ and its derivatives, the profiles $\mathcal{U}^j$ and
their $\eps$-conormal derivatives.
We purposely not specify these arguments for sake of clarity.
by iteration, we can also express the time derivatives $\partial_t^j \mathtt{V}^\eps_{R,d}$ for $j \geqslant 2$ by the equation
$(\ref{poi})$. Therefore there exist some $C^\infty$ functions $\H^j$ such that 
\begin{eqnarray*}
\partial_t^j \mathtt{V}^\eps_R |_{ t=x_d =0 } = \H_j ((\partial_x^l U^\eps_{R,init} |_{x_d =0 } )_{l \leqslant j} ) - \partial_t^j
R^\eps_{\mathtt{V}_d}  |_{t=x_d =0} .
\end{eqnarray*}
For $j \geqslant 1$, we define the $j$th order compatibility condition
\begin{eqnarray*}
(\mathcal{R}^{j,\eps}): \ \H_j ((\partial_x^l U^\eps_{R,init} |_{x_d =0 } )_{l \leqslant j} ) = \partial_t^j
R^\eps_{\mathtt{V}_d}  |_{t=x_d =0} .
\end{eqnarray*}
Thus for $\eps \in ]0,1]$, if for $T >0$, $U^\eps_{R} \in H^\infty ( (0,T) \times \Omega) $ satisfies the equation
$(\ref{poi})-(\ref{poi2})-(\ref{poi3})$ and $U^\eps_{R} |_{t=0} = U^\eps_{R,init}$, then $U^\eps_{R} $  satisfies the compatibility
condition $(\mathcal{R}^{j,\eps})_{j \in \N}$. 
Next proposition will show that there exist some families $(U^\eps_{R,init})_{\eps \in ]0,1]}$ which satisfy the  compatibility
condition $(\mathcal{R}^{j,\eps})_{j \in \N}$ and some estimates uniform with respect to $\eps$.
Moreover we can prescribe arbitrarily the components $(\mathtt{V}^\eps_{R,init,y} , \mathtt{S}^\eps_{R,init})_{\eps}$ 
among the families
which satisfy convenient uniform estimates.

\begin{prop}
\label{exidi}

Let $(\mathtt{V}^\eps_{R,init,y} , \mathtt{S}^\eps_{R,init})_{\eps}$ be a family bounded in $H^\infty (\Omega)$. 
Then there exists a family
$(\mathtt{V}^\eps_{R,init,d} , \mathtt{P}^\eps_{init})_{\eps}$ bounded in $H^\infty (\Omega)$ such that for all $\eps \in ]0,1]$, for all $j\in
\N$, the compatibility condition $\mathcal{R}^{\eps,j}$ is verified.

\end{prop}

The underlying reason to the possibility  to choose arbitrarily the tangential velocity and the entropy is that these components are
the characteristic ones. Indeed we will see in the proof below the crucial role played by the normal derivatives.
We need for uniform estimates in order to find a family of solutions $(U^\eps_R )_\eps$ in $E^\infty ((0,T_0) \times \Omega )$.
We succeed in this goal obtaining in Proposition $\ref{exidi}$ a family $(U^\eps_{R,init})_\eps$ which even do not contain
singularities in $\frac{x_d}{\eps}$.
We refer to \cite{moi}, \cite{moi2} for other examples of existence of initial data compatible at all order with uniform existence in
a setting of boundary layer theory.
We will use in the proof of Proposition $\ref{exidi}$, as in the references above, a Borel lemma.
\begin{proof}
We will proceed in three steps.
\begin{enumerate}
\item We begin to analyze more accurately the compatibility conditions $ \mathcal{R}^{\eps,j}$ and more particularly 
the way the functions $(\H_j )_{j \in \N^*}$ involve the normal derivatives of $\mathtt{v}^\eps_{R,init,d}$ and $p^\eps_{init}$.
Indeed there are some $C^\infty$ functions $(\tilde{\H}_j )_{j \in \N^*}$
such that for all $j \in \N^*$,
\begin{eqnarray*}
\H^{2j} ((\partial_x^l U^\eps_{R,init} |_{x_d =0 } )_{l \leqslant 2j} ) = 
\tilde{\H}^{2j} ((\partial_x^\alpha U^\eps_{R,init} |_{x_d =0 } )_{\alpha \in I_j} ) + (\frac{1}{\alpha^\eps})^j \partial_d^{2j}
\mathtt{v}^\eps_{R,init,d} ,
\\  \H^{2j+1} ((\partial_x^l U^\eps_{R,init} |_{x_d =0 } )_{l \leqslant 2j+1} ) = 
\tilde{\H}^{2j} ((\partial_x^\alpha U^\eps_{R,init} |_{x_d =0 } )_{\alpha \in I'_j} ) + (\frac{1}{\rho^\eps})^j \partial_d^{2j}
p^\eps_{R,init} .
\end{eqnarray*}
\item For all $j \in \N$, there exists a family of functions $( \rond{\mathtt{V}}^{\eps,j}_{R,init,d} , \rond{P}^{\eps,j}_{R,init}
)_\eps$ bounded in $H^\infty ((0,T_0) \times \Gamma )$ such that for all $j \in \N^*$
\begin{eqnarray*}
\tilde{\H}^{2j} ((\partial_x^{\alpha'} \rond{U}^{\eps,\alpha_d}_{R,init}  )_{\alpha \in I_j} )
 + (\frac{1}{\alpha^\eps})^j  \rond{\mathtt{V}}^{\eps,2j}_{R,init,d} =0,
\\  \tilde{\H}^{2j+1} ((\partial_x^{\alpha'} \rond{U}^{\eps,\alpha_d}_{R,init}  )_{\alpha \in I'_j} ) 
+ (\frac{1}{\rho^\eps})^j  \rond{\mathtt{P}}^{\eps,2j+1}_{R,init} =0.
\end{eqnarray*}
\item We use a Lemma by E.Borel in order to end the proof (cf. \cite{moi}).
\end{enumerate}
\end{proof}
\subsection{Linear estimate}
\label{Linesti}

We begin to look at the following linear problem:
\begin{eqnarray}
 \label{AL} \qquad (\mathfrak{L}^\star ( \underline{u}^\eps ,\partial)
+ \underline{J}^\alpha ).V^\eps_R +  W^\eps_R .(\underline{J}^{\beta,1} + 
\frac{x_d}{\eps} \underline{J}^{\beta,2})
 = - R^\eps_v
 \quad \mathrm{when} \ x_d >0 ,
\\ \label{BL} \X_{\underline{\mathtt{v}}^\eps } W^\eps_R + \eps \mathtt{V}^\eps_R . \nabla W^\eps_a = -\eps R^\eps_W 
\quad \mathrm{when} \ x_d > 0,
\\ \label{CL} \mathtt{V}^\eps_{R,d} = 0 \quad \mathrm{when} \ x_d =0,
\end{eqnarray}
where $\underline{J}^\alpha $ denotes the $(d+1) \times (d+1)$ matrix 
\begin{eqnarray*}
\underline{J}^\alpha &:=& J^\alpha (\eps, V^\eps_a , Z_\eps V^\eps_a , W^{\eps,\flat}_a, \eps^{M} \, \underline{U}^\eps_R  ),
\end{eqnarray*}
and for $i \in \{1,2\}$,  $\underline{J}^{\beta,i} $ denotes the $\R^{d+1}$-valued function
\begin{eqnarray*}
\underline{J}^{\beta,i} &:=& J^{\beta,i} (\eps, V^\eps_a , Z_\eps V^\eps_a , W^{\eps,\flat}_a, \eps^M \, \underline{U}^\eps_R  ) .
\end{eqnarray*}
  The family $(\underline{u}^\eps )_{\eps} := (\underline{v}^\eps , \underline{W}^\eps )_{\eps} $ is a given family of the form 
  \begin{eqnarray}
  \label{linform}
    \underline{v}^\eps := v^\eps_a + \eps^{M+1} \,  \underline{V}^\eps_R ,  
   \quad \underline{W}^\eps  := W^\eps_a + \eps^{M} \,  \underline{W}^\eps_R 
   \end{eqnarray}
   with
    \begin{eqnarray*}
   \underline{U}^\eps_R := (\underline{V}^\eps_R , \underline{W}^\eps_R ) \in \mathbf{E}^m (T_0  ). 
  \end{eqnarray*}
  We introduce the classic spaces of conormal distributions 
   \begin{eqnarray*}
  H^{0,m} (T) := \{ u \in L^2 (T ) / \ 0 \leqslant k  \leqslant m, \ Z^k u \in L^2 ( T ) \} .
   \end{eqnarray*}
  These spaces will be endowed by the following weighted norms:
\begin{eqnarray}
\label{defi}
|u|_{m,\lambda,T} &:=& \sum_{0 \leqslant k \leqslant m} \lambda^{m-k}  \, ||e^{-\lambda T} Z^k u ||_{L^2 (T)} .
\end{eqnarray}
Because we face a characteristic boundary problem with functions in $\frac{x_d}{\eps}$, we will also use the following  norms:
\begin{eqnarray}
\label{defiN}
 |u|^\mathcal{N}_{\eps, m,\lambda,T} &:=& |u|_{m,\lambda,T}  +  |(\eps \partial_d )u|_{m,\lambda,T}
\\ \label{defiE}  |u|^\mathbf{E}_{\eps, m,\lambda,T} &:=& 
\sum_{0 \leqslant 2k+l \leqslant m} \lambda^{m-2k-l}  \, |Z^l (\eps \partial_d )^k |_{0,\lambda,T}
\end{eqnarray}
A link between the $L^2$-type norm: $|.|^\mathbf{E}_{\eps,m,\lambda,T}$ and the $L^\infty$-type norm: $||.||^*_{\eps,T}$ is given by
the following Sobolev Embedding lemma:

\begin{lem}[Sobolev Embedding]
\label{sobo}

Let  $m > \frac{n}{2}+5$. There is  $c >0$ such that for all $T \in [0,T_0]$, 
  for all $u \in C_0^\infty ((0,T) \times \Omega )$,
\begin{eqnarray*}
\sqrt{\eps}  ||u||^*_{\eps,T} \leqslant c T e^{\lambda T }  \, |u|^\mathbf{E}_{\eps,m,\lambda,T} .
\end{eqnarray*}

\end{lem} 
\begin{proof}
 Let us introduce the family of functions $(\tilde{u}^\eps )_{\eps \in ]0,1]}$ 
by 
\begin{eqnarray*}
\tilde{u}^\eps (t,x ):= u(t,y,\eps x_d ) \quad \mathrm{for} \ \mathrm{all} \ (\eps,t,x )  \in ]0,1] \times [0,T_1]  \times 
\Omega. 
\end{eqnarray*}
We introduce $ ||u||^*_{T}  := ||u||_{0,T} + ||  u||_{1,T}$.
We get
\begin{eqnarray*}
|| \tilde{u}^\eps ||^*_T := ||u ||^*_{\eps,T} 
 \ \mathrm{and} \   || \tilde{u}^\eps ||^E_{m,\lambda,T} := \eps^{-\frac{1}{2}} \, ||u ||^\mathbf{E}_{\eps,m,\lambda, T} .
 \end{eqnarray*}
  However Lemma $II.1.1$ of \cite{Gu90} shows that there exists $c >0$ such that for all $u \in E^m$, 
  \begin{eqnarray*}
||u||^*_{T} \leq cTe^{\lambda T } |u|^E_{m,\lambda,T}  .
\end{eqnarray*}
Using this for the family $(\tilde{u}^\eps)_\eps $ leads to Lemma $\ref{sobo}$.
\end{proof}
We will use the following version of the Gagliardo-Nirenberg estimates (cf. \cite{Gu90}).
\begin{lem}[Gagliardo-Nirenberg]
\label{Ga1}
There exists $C>0$ such that for all $T \in ]0,T_0]$, for all $u \in L^\infty (T) \cap H^{0,m} (T)$, 
for all $l$ and $m$ such that $l  \leqslant k  \leqslant m$ and for all real $\lambda \geqslant 1$
\begin{eqnarray*}
\lambda^{k-l} |e^{-\frac{2}{p}\lambda t}  \, Z^l u |_{L^p ((0,T) \times \Omega )} \leqslant C ||u||^{1-\frac{k}{m}}_{0,T}
|u|^{\frac{k}{m}}_{m,\lambda,T} ,
\end{eqnarray*}
where $p:= \frac{2m}{k}$.
\end{lem} 
We will also use 
\begin{lem}[Gagliardo-Nirenberg]
\label{Ga2}
Let $m \in \N$ be even.
There exists $C>0$ such that for all $T \in ]0,T_0]$, for all $u \in L^\infty (T) \cap E (T)
$, for all $l$ and $m$ such that $l  \leqslant k  \leqslant m$, for all $\eps \in ]0,1]$ and for all real $\lambda \geqslant 1$,
\begin{eqnarray*}
\lambda^{k-l} |e^{-\frac{2}{p}\lambda t} (\eps \partial_d )^k Z^l u |_{L^p ((0,T) \times \Omega )} \leqslant 
C ||u||^{1-\frac{k}{m}}_{0,T} (|u|^\mathbf{E}_{\eps,m,\lambda,T})^{\frac{k}{m}},
\end{eqnarray*}
where $p:= \frac{2m}{k}$.
\end{lem} 
\begin{proof}
It is a special case of $(Ap-II-3)$ given in \cite{Gu90}, p.$643$.
\end{proof}
Lemma $\ref{Ga1}$ and $\ref{Ga2}$ imply the following Moser's type inequalities.
\begin{lem}[Moser Inequalities]
\label{Moser}
Let $F \in C^\infty  $ such that   $F(0) =0$, and  $\mathfrak{R} >0$. There exists a real  $l$ such that 
 for all $T \in [0,T_0]$, if  
 \begin{eqnarray*}
 g \in L^\infty (T ) \cap  H^{0,m} (T) 
 \quad (\mathrm{resp}. \ L^\infty (T ) \cap E^m (T ) \ \mathrm{and} \ m \  \mathrm{even} ),
 \end{eqnarray*}
   verifies  $||g||_{0,T} \leqslant  \mathfrak{R}$,
  we have, for all  $\lambda  \geqslant 1$, for all  $\eps \in ]0,1]$, we have: 
\begin{eqnarray*}
| F(g) |_{m,\lambda,T} \leqslant l | g |_{m,\lambda,T} ,\quad (\mathrm{resp}. \ | F(g) |^\mathbf{E}_{\eps, m,\lambda,T} 
\leqslant l | g |^\mathbf{E}_{\eps, m,\lambda,T} ).
\end{eqnarray*}
\end{lem} 
\begin{lem}
\label{Moser2}
For all  $\underline{m} \in \N$, there exists a real $c$ such that for all  $T$ verifying  $0 \leqslant T \leqslant T_0$, 
and  $\alpha, \beta \in L^\infty (T )$, we have:
\begin{enumerate}[$(1)$]
\item \label{mosera}
If  $\alpha, \beta \in H^{0,m} (T)$ and if  $j,k,l \in \N$ are such that 
 $0 \leqslant j+l \leqslant k \leqslant \underline{m} $:
\begin{eqnarray*}
\lambda^{\underline{m}-k} |Z^j \alpha \ Z^l \beta |_{\eps,0,\lambda,T} 
\leqslant c( ||\alpha ||_{0,T} \ |\beta |_{ \underline{m},\lambda,T} +  |\alpha |_{\underline{m},\lambda,T} \ ||\beta ||_{ 0,T })
\\ |\alpha \ \beta |_{\underline{m},\lambda,T} \leqslant c ( ||\alpha ||_{0,T} \  | \beta |_{ \underline{m},\lambda,T}  
 +  |\alpha |_ {\underline{m} ,\lambda,T} \   ||\beta ||_{ 0,T })
\end{eqnarray*}
\item If  $m$ is even, if $\alpha$, $\beta \in \mathbf{E}^m (T) ,$ and $ k,k',l,l',p  \in \N$ are 
 such  that  $2(k+k')+l+l' \leqslant p \leqslant \underline{m}$:
\begin{eqnarray*}
\lambda^{\underline{m}-p} |(\partial_n^k Z^l \alpha )(\partial_n^{k'} Z^{l'} \alpha ) |_{0,\lambda,T}
 &\leqslant& c (|| \alpha ||_{0,\lambda,T} \  |\beta|^\mathbf{E}_{\eps,\underline{m},\lambda,T} 
 + | \alpha |^\mathbf{E}_{\eps,\underline{m},\lambda,T} \  || \beta ||_{0,\lambda,T} ) ,
\\ | \alpha \  \beta |^\mathbf{E}_{\eps,\underline{m},\lambda,T} 
&\leqslant& c ( || \alpha  ||_{0,T}  \ |\beta|^\mathbf{E}_{\eps,\underline{m},\lambda,T}
 + | \alpha |^\mathbf{E}_{\eps,\underline{m},\lambda,T} \  || \beta ||_{0,T} ).
\end{eqnarray*}
\end{enumerate}
\end{lem}
We combine Lemma $\ref{Moser}$ and Lemma $\ref{Moser2}$ to find the following corollary.
\begin{coro}
Let  $F \in C^\infty (\R^k ,\R )$ such that  $F(0)=0$ and  $\mathfrak{R} >0$. 
 There exists  $C>0$ such that, for all  $T \in [0,T_0 ]$, 
for all  function $g \in Lip((0,T) \times \Omega ) \cap \mathcal{N} (T )$, $g: \Omega_T \rightarrow \R^k $
 which verify  $||g||_{Lip(T)} \leqslant \mathfrak{R}$ then  $F(g) \in  \mathcal{N} (\Omega_T )$ and we have:
\begin{eqnarray*}
|F(g)|^\mathcal{N}_{\eps,m,\lambda,T} \leqslant C |g|^\mathcal{N}_{\eps,m,\lambda,T} .
\end{eqnarray*}
\end{coro}
\begin{proof}
Referring to the definition of $|.|^\mathcal{N}_{\eps,m,\lambda,T}$ (see $(\ref{defiN})$),
 we get
\begin{eqnarray*}
|F(g)|^\mathcal{N}_{\eps,m,\lambda,T} &:=& |F(g)|_{m,\lambda,T} + | \eps  \partial_d F(g)|_{m-2,\lambda,T} .
\end{eqnarray*}
We use Lemma $\ref{Moser}$ to bound $|F(g)|_{m,\lambda,T}$.
Referring to the definition of $|.|_{m,\lambda,T}$ (see  $(\ref{defi})$),
 we get
\begin{eqnarray*}
 | \eps  \partial_d F(g)|_{m-2,\lambda,T}    & =&  
  \sum_{0 \leqslant k \leqslant m-2} \lambda^{m-2-k} ||e^{-\lambda T}  \, Z^k \eps \partial_d F(g)||_{L^2 (T)} 
\\     & =&  \sum_{0 \leqslant k= k_1 + k_2  \leqslant m-2} \lambda^{m-2-k} \mathfrak{N}_{k_1 , k_2 } 
\end{eqnarray*}
where 
\begin{eqnarray*}
 \mathfrak{N}_{k_1 , k_2 } := ||e^{-\lambda T}   \, Z^{k_1}  D_g F(g) .  Z^{k_1} \eps \partial_d g ||_{L^2 (T )} .
\end{eqnarray*}
Then we use Lemma $\ref{Moser2}$ $(\ref{mosera})$. 
\end{proof}
We will prove in subsection $\ref{L2esti}$ and $\ref{hoe}$ the following linear $\eps$-uniform estimates 
which will be useful in subsection  $\ref{Itsch}$. 
\begin{theo}
\label{pok}
Let $m$ be an even integer larger than $2$. Let $\mathfrak{R}$ be a strictly positive real. 
Then there is $\lambda_0 > 1$ such that if $\eps^M ||
\underline{U}_R^\eps ||^*_{\eps,T_0} \leqslant \mathfrak{R}$ then 
for all $\lambda \geqslant \lambda_0$, for all $\eps \in ]0,\eps_0]$,
\begin{eqnarray*}
|U_R^\eps |^\mathbf{E}_{\eps,m,\lambda,T_0 } \leqslant
 \lambda^{-1} \lambda_0 .( 1+ \eps^M  \, || U_R^\eps ||^*_{\eps,T_0} . | \underline{U}_R^\eps |^\mathbf{E}_{\eps,m,\lambda,T_0}  ).
\end{eqnarray*}
\end{theo}

We will proceed in several steps. 
In subsection $\ref{L2esti}$, we will perform $L^2$ estimates proceeding in three steps, estimating $W^\eps_R$ then $\frac{x_d}{\eps}
W^\eps_R$ and finally $V^\eps_R$. 
higher order estimates are more delicate to obtain.
In subsection $\ref{hoe}$, we will give preliminary results with some commutators estimates (cf. Proposition $\ref{lapin}$ and with
some  estimates of $\mathtt{v}^\eps_d . \partial_d U^\eps_R $ (cf. Lemma $\ref{carotte}$). We will begin with conormal estimates
(subsection $\ref{conorm}$) then looking for normal estimates (subsection $\ref{norm}$). 
As in the $L^2$ estimate, we will look for estimates about $W^\eps_R $, $\frac{x_d}{\eps}
W^\eps_R$ and finally $V^\eps_R$.
Furthermore, in the normal estimate, we will distinguish the estimates about tangential velocity $\mathtt{v}^\eps_{R,y}$ which
corresponds to characteristic components and the estimates about normal velocity $\mathtt{v}^\eps_{R,d}$ which
corresponds to a noncharacteristic component. 
\subsection{$L^2$ estimate}
\label{L2esti}
We will proceed in three steps estimating $W^\eps_R$ then $\frac{x_d}{\eps}
W^\eps_R$ and finally $V^\eps_R$.  
According to section $5$, $\nabla W^\eps_a$ is of the form $\beta_1 (\eps,t,x) +
\eps^{-1} \beta_2 (\eps,t,x,\frac{x_d}{\eps})$ where $\beta_1$ and $\beta_2$ are
$C^\infty_b$ with respect to  their arguments including $\eps$ up to $0$. 
Moreover  $\beta_2$ is rapidly decreasing  with respect to  its last argument.
\subsubsection{Estimate of $W^\eps_R$}
\label{estiW}
The aim of this subsection is to prove 
\begin{prop}
Let $\mathfrak{R} >0$. There exists $\lambda_0 > 1$ such that if $\eps^M \, || \underline{V}^\eps_R ||^*_{\eps,T_0} \leqslant \mathfrak{R}$
 then for all $\lambda \geqslant \lambda_0$, 
 \begin{eqnarray*}
 | W^\eps_R |_{0,\lambda,T_0} \leqslant \frac{\lambda}{\lambda_0} (1+ | V^\eps_R |_{0,\lambda,T_0}).
 \end{eqnarray*}
\end{prop}
\begin{proof}
We rewrite the equation $(\ref{BL})$ as  $\X_{\underline{\mathtt{v}^\eps}} W^\eps_R =f^\eps_I$ 
where $ f^\eps_I :=  (\eps \beta_1 + \beta_2 ) V^\eps_R -\eps R^\eps_W $. 
By a classic $L^2$ estimate, we obtain that for $\mathfrak{R} >0$, there exists $\lambda_0 \geqslant 1$ such that if
$||\underline{\mathtt{v}}^\eps ||_{Lip} \leqslant  \mathfrak{R}$ then
 for all $\lambda \geqslant \lambda_0$, 
 \begin{eqnarray*}
 | W^\eps_R |_{0,\lambda,T_0} \leqslant \frac{\lambda_0}{\lambda}
|f^\eps_I  |_{0,\lambda,T_0}.
 \end{eqnarray*}
Using $(\ref{linform})$
and estimating $| f^\eps_I  |_{0,\lambda,T_0}$ we
end the proof.
\end{proof}
\subsubsection{Estimate of $\frac{x_d}{\eps} W^\eps_R$}
\label{estiXW}
We will exploit the special form of the equation $(\ref{BL})$ thanks to a second estimate 
which concerns $\frac{x_d}{\eps} W^\eps_R$.
\begin{prop}
Let $\mathfrak{R}  >0$. 
There exists $\lambda_0 > 1$ such that if $\eps^M || \underline{V}^\eps ||_{\eps,T_0}^* \leqslant \mathfrak{R}$ 
then for all $\lambda \geqslant \lambda_0$, 
$ |\frac{x_d}{\eps} W^\eps_R |_{0,\lambda,T_0} \leqslant \frac{\lambda}{\lambda_0} (1+ | V^\eps_R |_{0,\lambda,T_0})$.
\end{prop}
\begin{proof}
We will proceed in two steps.
First we will look for an equation for $\frac{x_d}{\eps} W^\eps_R$ and then we will use a $L^2$
estimate for this equation.
\begin{enumerate}[$1$.]  
\item We begin to calculate 
\begin{eqnarray}
\label{aq}
\X_{\underline{\mathtt{v}}^\eps}  \, ( \frac{x_d}{\eps} W^\eps_R ) = \frac{x_d}{\eps} \X_{\underline{\mathtt{v}}^\eps}  \, W^\eps_R +
\underline{\mathtt{v}}_d^\eps . \frac{1}{\eps} W^\eps_R .
\end{eqnarray}
Because $\underline{\mathtt{v}}_d^\eps =0 $ when $x_d =0$, there exist some $C^\infty$ functions $
\underline{\mathtt{v}}^{\eps,\flat}$ such that $\underline{\mathtt{v}}_d^\eps = x_d  \underline{\mathtt{v}}^{\eps,\flat}$.
Using the equation $(\ref{BL})$, we deduce from the equation $(\ref{aq})$ that $\frac{x_d}{\eps} W^\eps_R $ verifies the equation 
\begin{eqnarray}
\label{aq2}
\X_{\underline{\mathtt{v}}^\eps}  \, ( \frac{x_d}{\eps} W^\eps_R ) -  \underline{\mathtt{v}}^{\eps,\flat} ( \frac{x_d}{\eps} W^\eps_R ) 
= f^\eps_{II} ,
\end{eqnarray}
where 
\begin{eqnarray}
f^\eps_{II} &:=& \frac{x_d}{\eps} (- \eps \mathtt{V}^\eps_R . \nabla W^\eps_a -  \eps R^\eps_W ) ,
\\          &:=& -x_d \beta_1 V^\eps_R - \frac{x_d}{\eps} \beta_2 V^\eps_R - x_d R^\eps_W  .
\end{eqnarray}
\item By a classic $L^2$ estimate, we obtain that for $\mathfrak{R} >0$, there exists $ \lambda_0 \geqslant 1$ such that if 
\begin{eqnarray*}
||\underline{\mathtt{v}}^{\eps} ||_{Lip ((0,T_0 ) \times \Omega )} 
+ || \underline{\mathtt{v}}_d^{\eps} ||_{L^\infty (T_0  )}  \leqslant R
\end{eqnarray*}
 then for all $\lambda  \geqslant  \lambda_0$ ,
\begin{eqnarray*}
| \frac{x_d}{\eps} W^\eps_R |_{0,\lambda,T_0} \leqslant \frac{\lambda_0}{\lambda} | f^\eps_{II} |_{0,\lambda,T_0} .
\end{eqnarray*}
Because of the form of $\underline{\mathtt{v}}^\eps $ (see $(\ref{linform})$), the family $(||
\underline{\mathtt{v}}^\eps ||_{Lip} + || \underline{\mathtt{v}}^{\eps,\flat} ||_{L^\infty} )_\eps $ 
is bounded when the family $(\eps^{M} || \underline{V}^\eps_R  ||^*_{\eps,T_0}  )_\eps$ is bounded. 
Moreover because $\beta_2$ is rapidly decreasing with respect to its last argument,  
the family $(\frac{x_d}{\eps} \beta_2 (\eps,t,x,\frac{x_d}{\eps}))_\eps $ is
bounded in $L^\infty$. This ends the proof.
\end{enumerate}
\end{proof}
\subsubsection{Estimate of $ V^\eps_R$}
\label{estiV}
We will prove the following result:
\begin{lem}
Let $\mathfrak{R} >0$. There exists $ \lambda_0 > 0 $ such that 
if $\eps^{M} || \underline{U}^\eps_R ||_{\eps,T_0}^* \leqslant  \mathfrak{R}$ 
then for all $\lambda \geqslant \lambda_0$, $ | V^\eps_R |_{0,\lambda,T_0} \leqslant  \frac{ \lambda_0}{\lambda}$.
\end{lem}
\begin{proof}
For each $\eps \in ]0,1]$, the function $V^\eps_R$ is solution of the following boundary value problem: 
\begin{eqnarray}
\label{Sx}
\left\{
\begin{array}{ccc}
&(  \mathfrak{L}^\star (\underline{u}^\eps  ,\partial) + \underline{J}^\alpha  ).V^\eps_R =
f^\eps_{III} \quad &\mathrm{when} \ x_d > 0 ,
\\ &\mathtt{V}^\eps_{R,d} = 0  &\quad \mathrm{when} \ x_d = 0 ,
\end{array}
\right.
\end{eqnarray}
where 
\begin{eqnarray*}
f^\eps_{III}  :=  - W^\eps_R \underline{J}^{\beta,1}  - \frac{x_d}{\eps}  W^\eps_R \underline{J}^{\beta,2} - R^\eps_W .
\end{eqnarray*}
It is a symmetric hyperbolic system for which the boundary is conservative and characteristic of constant multiplicity. 
Thanks to \cite{R}, the following $L^2$ estimate holds: for $\mathfrak{R}> 0$, there exists $\lambda > 0$, such that if  
\begin{eqnarray*}
||\underline{\mathtt{v}}^\eps ||_{Lip} + || \underline{J}^\alpha ||_{L^\infty} 
+ || \underline{J}^\beta ||_{L^\infty} \leqslant \mathfrak{R}
\end{eqnarray*}
then for all $\lambda \geqslant \lambda_0$, 
\begin{eqnarray*}
|V^\eps_R |_{0,\lambda,T_0} \leqslant  \frac{\lambda}{\lambda_0} |f^\eps_{III} |_{0,\lambda,T_0}  .
\end{eqnarray*}
Using the form of $\underline{\mathtt{v}}^\eps $, $\underline{J}^\alpha$, $\underline{J}^\beta$ and 
estimating  $|f^\eps_{III} |_{0,\lambda,T_0}$ thanks to subsection $\ref{estiW}$ and $\ref{estiXW}$ yields the conclusion.
\end{proof}
\subsection{Higher order estimates} 
\label{hoe}
We look for estimates uniform with respect to $0<\eps \leqslant 1$ of 
\begin{eqnarray*}
\lambda^{m-2k-l} |(\eps \partial_d )^k Z^l U^\eps_R |_{0,\lambda,T_0 } \quad \mathrm{for} \  2k+l \leqslant m .
\end{eqnarray*}
In this subsection, we give crucial preliminary technical results. 
We group them in two kinds.
In a first time, we look for commutators estimates. 
In a second time, we look for estimates which will be useful when estimating source terms.
\subsubsection{Commutators estimates} 
We begin with the estimates of the commutators $[ \underline{J} , (\eps \partial_d )^k Z^l ] \varphi$.
Because we will proceed in two steps, estimating first conormal estimates then normal estimates, we give specific estimates
of the commutators with the derivatives $Z^l$ for $l \leqslant m$ (i.e. in the limit case $k=0$).
\begin{prop}
\label{pomme}
Let $m$ be an integer. 
There is an increasing function $C:\R_+ \rightarrow \R_+^*$ such that for all $\eps \in ]0,1]$, 
if $\eps^M  \, || \underline{U}^\eps_R ||_{0,T}
\leqslant  \mathfrak{R}$ then for all $\varphi \in C^\infty_0 ((0,T) \times \overline{\Omega} )$,
 for all $k,l$ such that $2k+l \leqslant m$, for all $\lambda
\geqslant \lambda_0$,
\begin{enumerate}[$a)$]
\item \label{pomme1}
if $m$ is even, 
\begin{eqnarray*}
\lambda^{m-2k-l} |[ \underline{J} (\underline{u}^\eps ) , (\eps \partial_d )^k  Z^l ] \varphi|_{0,\lambda,T} 
\leqslant C(\mathfrak{R}) ( |\varphi |^\mathbf{E}_{\eps, m,\lambda,T} 
+ \eps^M || \varphi ||_{0,T} .|\underline{U}^\eps_R |^\mathbf{E}_{\eps,m,\lambda,T} ) , 
\end{eqnarray*}
 \item \label{pomme2}
 for any $m$, 
\begin{eqnarray*}
\lambda^{m-l} |[ \underline{J}  (\underline{u}^\eps ) ,  Z^l ] \varphi |_{0,\lambda,T}
 \leqslant C(\mathfrak{R}) ( |\varphi |_{m,\lambda,T} 
+ \eps^M || \varphi ||_{0,T} .|\underline{U}^\eps_R |_{m,\lambda,T} ) , 
\end{eqnarray*}
\end{enumerate}
where $\underline{J}$ stands for $\underline{J}^\alpha $, $\underline{J}^{\beta,1}$ or $\underline{J}^{\beta,2}$.
\end{prop}
\begin{proof}[Proof of $(\ref{pomme2})$]
We introduce a notation: if $l \in \N, l \neq 0$ and $\phi$ a vector-valued function which components are
denoted  $\phi_i$, we denote $Z^{<l>} \phi$ the collection of the terms
of the form $Z^{l_1}  \phi_{i_1} ... Z^{l_r}  \phi_{i_r} $
 where $1 \leqslant  r  \leqslant l , \ l_1 + ... l_r = l $
and $l_i \geqslant 1$ and if $l =0$, $Z^{<l>} \phi =0$.
The commutator $[ \underline{J} ,  Z^l ] \varphi$ is a linear combination of terms of the form 
\begin{eqnarray*}
\Phi (\eps, V^\eps_a , Z_\eps V^\eps_a , W^{\eps,\flat}_a , \eps^M \, \underline{U}^\eps_R ). \mathfrak{P}
\end{eqnarray*}
 with 
 \begin{eqnarray*}
 \mathfrak{P} := Z^{<l_1>} ( V_{a}^{\eps} ) . Z^{<l_2>} ( Z_\eps  V_{a}^{\eps} ) . Z^{<l_3>} (  W^{\eps,\flat}_a  ) .
Z^{<l_4>} (\eps^M \, \underline{U}^\eps_R ).  Z^{l_5}  \phi ,
\end{eqnarray*}
and  $l_1 + l_2 + l_3 + l_4 + l_5= l $.
 Thanks to estimates $(\ref{linf})$ and Lemma $\ref{Moser}$, we get
 \begin{eqnarray*}
\lambda^{m-l} |[ \underline{J}  (\underline{u}^\eps ) ,  Z^l ] \varphi |_{0,\lambda,T}
 \leqslant C(\mathfrak{R}) \lambda^{m-l} 
 | Z^{<l_4>} (\eps^M \, \underline{U}^\eps_R ).  Z^{l_5}  \phi |_{m,\lambda,T}  .
 \end{eqnarray*}
We end the proof thanks to Lemma $\ref{Moser2}$ $(1)$.
\end{proof}
\begin{proof}[Proof of $(\ref{pomme1})$]
We introduce a notation: if $l \in \N, l \neq 0$ and $\phi$ a vector-valued function which components are
denoted  $\phi_i$, we denote $(\eps \partial_d)^{<k>} \phi$ the collection of the terms
of the form $ (\eps \partial_d)^{k_1}  \phi_{i_1} ... (\eps \partial_d)^{k_r}  \phi_{i_r} $
 where $1 \leqslant  r  \leqslant l , \ k_1 + ... + k_r = l $
and $k_i \geqslant 1$ and if $l =0$, $(\eps \partial_d)^{<k>} \phi =0$.
The commutator $[ \underline{J} ,  (\eps \partial_d )^k Z^l ] \varphi$ is a linear combination of terms of the form 
\begin{eqnarray*}
\Phi (\eps, V^\eps_a , Z_\eps V^\eps_a , W^{\eps,\flat}_a , \eps^M \, \underline{U}^\eps_R ). \tilde{\mathfrak{P}}
\end{eqnarray*}
with 
 \begin{eqnarray*}
\tilde{\mathfrak{P}} := 
(\eps \partial_d)^{<k_1>} Z^{<l_1>}  ( V_{a}^{\eps} ) . (\eps \partial_d)^{<k_2>} Z^{<l_2>} ( Z_\eps  V_{a}^{\eps} ) 
 (\eps \partial_d)^{<k_3>} Z^{<l_3>} (  W^{\eps,\flat}_a  ) .
\\ .(\eps \partial_d)^{<k_4>}  Z^{<l_4>}   (\eps^M \, \underline{U}^\eps_R ).  (\eps \partial_d)^{<k_5>} Z^{l_5} \underline{J} \phi ,
\end{eqnarray*}
 with $l_1 + l_2 + l_3 + l_4 + l_5= l $, $k_1 + k_2 + k_3 + k_4 + k_5 = k $.
 Thanks to estimates $(\ref{linf})$ and Lemma $(\ref{Moser})$, we get
 \begin{eqnarray*}
\lambda^{m-2k-l} |[ \underline{J}  (\underline{u}^\eps ) ,  Z^l ] \varphi |_{0,\lambda,T}
 \leqslant C(\mathfrak{R}) \lambda^{m-l} 
 | Z^{<l_4>} (\eps^M \, \underline{U}^\eps_R ).  Z^{l_5} \underline{J} \phi |_{m,\lambda,T}  .
 \end{eqnarray*}
We end the proof thanks to Lemma $(\ref{Moser2})$ $(2)$.
\end{proof}
 A main point would be to estimate some commutators involving $S(\underline{u}^\eps )
X_{\underline{\mathtt{v}}^\eps}$.   
We will estimate separately $S(\underline{u}^\eps ) \underline{\mathtt{v}}^\eps_d \partial_d$ and the one $S(\underline{u}^\eps )
X'_{\underline{\mathtt{v}}^\eps}$ (see $(\ref{fd})$ for the definition of $X'_{\underline{\mathtt{v}}^\eps}$).
Moreover because we will proceed in two steps, estimating first conormal estimates then normal estimates, we give specific estimates
of the commutators with the derivatives $Z^l$ for $l \leqslant m$ (i.e. in the limit case $k=0$).
\begin{prop}
\label{lapin}
Let $m$ be an integer. 
There is an increasing function $C:\R_+ \rightarrow \R_+^*$ such that for all $\eps \in ]0,1]$, 
if $\eps^M  \, || \underline{U}^\eps_R ||^*_{\eps,T}
\leqslant  \mathfrak{R}$ then for all $\varphi \in C^\infty_0 ((0,T) \times \overline{\Omega} )$,
 for all $k,l$ such that $2k+l \leqslant m$, for all $\lambda
\geqslant \lambda_0$,
\begin{enumerate}[$a)$]
\item if $m$ is even, 
\begin{eqnarray*}
\lambda^{m-2k-l} |[S(\underline{u}^\eps ).\underline{\mathtt{v}}^\eps_d \partial_d , (\eps \partial_d )^k  Z^l ] \varphi
|_{0,\lambda,T} \leqslant C(\mathfrak{R}) ( |\varphi |^\mathbf{E}_{\eps, m,\lambda,T} 
+ \eps^M || \varphi ||^*_{\eps,T} .|\underline{U}^\eps_R |^\mathbf{E}_{\eps,
m,\lambda,T} ) , 
\end{eqnarray*}
\item for any $m$, 
\begin{eqnarray*}
\lambda^{m-l}  |[S(\underline{u}^\eps ).\underline{\mathtt{v}}^\eps_d \partial_d ,  Z^l ] \varphi |_{0,\lambda,T} \leqslant 
C(\mathfrak{R}) (|\varphi|^\mathcal{N}_{\eps, m,\lambda,T} 
+ \eps^M || \varphi ||^*_{\eps,T} .|\underline{U}^\eps_R |^\mathcal{N}_{\eps,m,\lambda,T}),
\end{eqnarray*}
\item if $m$ is even, 
\begin{eqnarray*}
\lambda^{m-2k-l} |[S(\underline{u}^\eps ). \X'_{\underline{\mathtt{v}}^\eps} , (\eps \partial_d )^k  Z^l ] \varphi
|_{0,\lambda,T} \leqslant C(\mathfrak{R})  ( |\varphi |^\mathbf{E}_{\eps, m,\lambda,T} 
+ \eps^M ||\varphi||_{1,T} .|\underline{U}^\eps_R |^\mathbf{E}_{\eps,
m,\lambda,T} ) , 
\end{eqnarray*}
\item for any $m$, 
\begin{eqnarray*}
\lambda^{m-l}  |[S(\underline{u}^\eps ). \X'_{\underline{\mathtt{v}}^\eps} ,  Z^l ] \varphi |_{0,\lambda,T} \leqslant 
  C(\mathfrak{R})  (|\varphi|_{m,\lambda,T} 
  + \eps^M || \varphi ||_{1,T} .|\underline{U}^\eps_R |_{m,\lambda,T}).
\end{eqnarray*}
\end{enumerate}
\end{prop}
We begin with the proof of $c)$ and $d)$ which are simpler.
Because $X'_{\underline{\mathtt{v}}} := \partial_t + \sum_{i=1}^{d-1} \mathtt{v}_i \partial_i $, the estimates  $c)$ and $d)$ are
consequences of the following lemma.
\begin{lem}
Let $m$ be an integer.
Let $F$ be a $C^\infty$ function.
There is an increasing function $C:\R_+ \rightarrow \R_+^*$ such that for all $\eps \in ]0,1]$, 
if $\eps^M  \, || \underline{U}^\eps_R ||^*_{\eps,T}
\leqslant  \mathfrak{R}$ then for all $\varphi \in C^\infty_0 ((0,T) \times \overline{\Omega} )$,
 for all $k,l$ such that $2k+l \leqslant m$, for all $\lambda
\geqslant \lambda_0$,
\begin{enumerate}[$(i)$]
\item \label{re1} if $m$ is even, 
\begin{eqnarray*}
\lambda^{m-2k-l} |[ F(\underline{u}^\eps ). Z , (\eps \partial_d )^k  Z^l ] \varphi|_{0,\lambda,T} \leqslant 
C(\mathfrak{R})  ( |\varphi |^\mathbf{E}_{\eps, m,\lambda,T}  
+ \eps^M ||\varphi||_{1,T} .|\underline{U}^\eps_R |^\mathbf{E}_{\eps,m,\lambda,T} ) , 
\end{eqnarray*}
\item \label{re2} for any $m$, 
\begin{eqnarray*}
\lambda^{m-l}  |[ F(\underline{u}^\eps ). Z ,  Z^l ] \varphi |_{0,\lambda,T} \leqslant 
  C(\mathfrak{R})  (|\varphi|_{m,\lambda,T} 
  + \eps^M || \varphi ||_{1,T} .|\underline{U}^\eps_R |_{m,\lambda,T}).
\end{eqnarray*}
\end{enumerate}
\end{lem}
\begin{proof}[Proof of $(\ref{re2})$]
The commutator $[ F(\underline{u}^\eps ). Z ,  Z^l ] \varphi $ is a linear combination of terms of the form
\begin{eqnarray}
\label{bb}
Z^{l_1} F(\underline{u}^\eps ). Z^{l_2} \varphi \quad \mathrm{when} \ l_1 + l_2 = l +1, \ l_1 \geqslant 1, \ l_2  \geqslant 1.
\end{eqnarray}
We introduce $l_1' := l_1 -1$, $l_2' := l'_2 -1$.
Thus the term $(\ref{bb}$) reads 
\begin{eqnarray*}
Z^{l'_1} (Z F(\underline{u}^\eps )). Z^{l'_2} (Z \varphi ).
\end{eqnarray*}
We apply Lemma $\ref{Moser2}$ $(1)$ with 
\begin{eqnarray*}
\alpha := Z F(\underline{u}^\eps ) \ ,  \beta:= Z \varphi \ , \underline{m} := m-1   
\end{eqnarray*}
and obtain the estimate
\begin{eqnarray*}
\lambda^{m-l}  |[ F(\underline{u}^\eps ). Z ,  Z^l ] \varphi |_{0,\lambda,T} \leqslant 
C ( ||F(\underline{u}^\eps ) ||_{1,T} .  |\varphi |_{m,\lambda,T}
  +  |F(\underline{u}^\eps ) |_{m,\lambda,T} . || \varphi  ||_{1,T} ).
\end{eqnarray*}
We use that $(u^\eps )_\eps$ is of the form $(\ref{linform})$ and apply Lemma $\ref{Moser}$ to complete the proof.
\end{proof}
\begin{proof}[Proof of $(\ref{re1})$]
We assume that $k \geqslant 1$. 
The case $k=0$ corresponds to $(\ref{re2})$.
The commutator $[ F(\underline{u}^\eps ). Z , (\eps \partial_d )^k  Z^l ] \varphi$ is a linear combination of terms of the form
\begin{eqnarray}
\label{bb1}
&(\eps \partial_d )^k  Z^{l_1} F(\underline{u}^\eps ). Z^{l_2} \varphi
 \quad &\mathrm{with} \ l_1 + l_2 = l +1, \ l_1 \geqslant 1, \ l_2  \geqslant 1,
\\ \label{bb2}
&(\eps \partial_d )^{k_1}  F(\underline{u}^\eps ). (\eps \partial_d )^{k_2} Z^{l+1} \varphi 
\quad &\mathrm{with} \ k_1 + k_2 = k, \ k_1  \geqslant 1 .
\end{eqnarray}
\begin{enumerate}[$\bullet$]
\item We begin to deal with the term $(\ref{bb1})$.
We introduce $k' := k -1$, $l_2' := l'_2 -1$.
We commute $\eps \partial_d$ with $Z^{l_1}$.
Thus it suffices to deal  with terms of the form
\begin{eqnarray*}
 (\eps \partial_d )^{k'}  Z^{l_1} (\eps \partial_d  F(\underline{u}^\eps )). Z^{l'_2} (Z \varphi ) .
\end{eqnarray*}
We apply Lemma $\ref{Moser2}$ $(1)$ with 
\begin{eqnarray}
\label{jus}
\alpha :=  \eps \partial_d F(\underline{u}^\eps ) \ ,  \beta:= Z \varphi \ , \underline{m} := m-2 .   
\end{eqnarray}
We obtain the estimate
\begin{eqnarray*}
\lambda^{m-l}  |[ F(\underline{u}^\eps ). Z ,  Z^l ] \varphi |_{0,\lambda,T} \leqslant 
C ( ||F(\underline{u}^\eps ) ||_{\eps,Lip,T} .  |\varphi |^\mathbf{E}_{m,\lambda,T}
  +  |F(\underline{u}^\eps ) |^\mathbf{E}_{m,\lambda,T} . || \varphi  ||_{1,T} ).
\end{eqnarray*}
We end the estimate as in the proof of $(\ref{re2})$.
\item
To deal with the term $(\ref{bb2}$) and complete the proof, we directly apply  Lemma $\ref{Moser2}$ $(1)$ to $(\ref{jus})$ and
proceed in the same way.
\end{enumerate}
\end{proof}
We tackle the Proof of $b)$.
\begin{proof}[Proof of $b)$]
The commutator $[S(\underline{u}^\eps).\underline{\mathtt{v}}^\eps_d \partial_d ,  Z^l ] \varphi $
is a linear combination of some terms of the form 
\begin{eqnarray*}
 Z^{l_1} S(\underline{u}^\eps) . Z^{l_2} \underline{\mathtt{v}}^\eps_d . Z^{l_3} \partial_d  \varphi , 
\quad \mathrm{with} \  l_1 + l_2 + l_3  = l, \  l_1 + l_2 \geqslant 1 .
\end{eqnarray*}
Because $\underline{\mathtt{v}}^\eps$ is of the form $\underline{\mathtt{v}}^\eps =  \mathtt{v}^0  + \eps
\underline{\mathtt{V}}^\eps  $, it suffices to deal with the terms 
\begin{eqnarray}
\label{ff1}  
Z^{l_1} S(\underline{u}^\eps) . Z^{l_2} \mathtt{v}^0_d . Z^{l_3} \partial_d  \varphi ,
\\ \label{ff2} 
Z^{l_1} S(\underline{u}^\eps) . Z^{l_2} \underline{\mathtt{V}}^\eps_d . Z^{l_3} (\eps \partial_d ) \varphi .
\end{eqnarray}
\begin{step2} 
We deal with the terms of the form $(\ref{ff1})$.
\end{step2}
Using that $ \mathtt{v}^0_d = x_d \mathtt{v}^{0,\flat}_d$ and commuting $x_d$ with $Z^{l_3}$, 
we have to deal with the terms of the form
\begin{eqnarray*}
Z^{l_1} S(\underline{u}^\eps) . Z^{l_2} \mathtt{v}^0_d . Z^{l'_3} Z \varphi ,
\quad \mathrm{with} \  l_3' \leqslant l_3 .
\end{eqnarray*}
We use $L^\infty$ estimate of $\mathtt{v}^0$.
Thus we have to control 
\begin{eqnarray}
\label{mm}
\lambda^{m-l} |Z^{l_1} S(\underline{u}^\eps) . Z^{l'_3} Z \varphi |_{0,\lambda,T} .
\end{eqnarray}
We distinguish two cases:
\begin{enumerate}[$\bullet$]
\item if $l_1 \leqslant 1$, we control $Z^{l_1} S(\underline{u}^\eps)$ in $L^\infty$ norm by Lemma
$\ref{Moser}$ and  the term $(\ref{mm})$ is bounded by $c. |\varphi |_{m,\lambda,T} $. 

\item if $l_1 \geqslant 2$, we apply Lemma $\ref{Moser2}$ with 
\begin{eqnarray*}
\alpha := Z^2 S(\underline{u}^\eps) ,\ \beta :=  Z \varphi  ,\ \underline{m} := m-2.  
\end{eqnarray*}
\end{enumerate}
\begin{step2} 
We deal with the terms of the form $(\ref{ff2})$.
\end{step2}
We distinguish two cases:
\begin{enumerate}[$\bullet$]

\item if $l_1 \leqslant 1$, we control $Z^{l_1} S(\underline{u}^\eps )$ in $L^\infty$ norm and then we
deal with 
\begin{eqnarray*}
\lambda^{m-l} |Z^{l_2} \underline{\mathtt{V}}^\eps_d .  Z^{l_3} (\eps \partial_d ) \varphi |_{0,\lambda,T} .
\end{eqnarray*}
with $l_2 + l_3 \leqslant l$. We distinguish two sub cases:
\begin{enumerate}[$(a)$]
\item if $l_2 \geqslant 2$, we introduce $l_2' := l_2 -2$.
Because $\underline{\mathtt{V}}^\eps_d = \underline{\mathtt{V}}^\eps_{a,d} + \eps^M \, \underline{\mathtt{V}}^\eps_{R,d}$,
we study the terms 
\begin{eqnarray}
\label{AAA1}
\lambda^{m-l} | Z^{l'_2} (Z^2 \mathtt{V}^\eps_{a,d} ) . Z^{l_3} (\eps  \partial_d ) \varphi |_{0,\lambda,T} ,
\\ \label{AAA2}
\eps^M \lambda^{m-l} | Z^{l'_2} (Z^2 \underline{\mathtt{V}}^\eps_{R,d} ) . Z^{l_3} (\eps  \partial_d ) \varphi |_{0,\lambda,T} .
\end{eqnarray}
We use the $L^\infty$ estimates of $(Z^2 \mathtt{V}^\eps_{a,d} )$ to control the term $(\ref{AAA1})$.
To control the term $(\ref{AAA2})$, we apply Lemma $\ref{Moser2}$ with 
\begin{eqnarray*}
\alpha := Z^2 \underline{\mathtt{V}}^\eps_{R,d} ,\ \beta :=  \eps \partial_d  \varphi  ,\ \underline{m} := m-2.  
\end{eqnarray*}
\item if $l_2 \leqslant 1$, we use the equality $\underline{\mathtt{V}}^\eps_{d} = x_d \underline{\mathtt{V}}^{\eps,\flat}_{d}$.
 We commute the factor $x_d$ with the derivative $ Z^{l_3}$. 
 Thus, to control the term $(\ref{AAA2})$, we study terms of the form
\begin{eqnarray*} 
\lambda^{m-l} |  Z^{l_2} (\eps \underline{\mathtt{V}}^{\eps,\flat}_{d} ) .  Z^{l'_3} (x_d  \partial_d ) \varphi |_{0,\lambda,T} ,
\end{eqnarray*}
with $l'_3 \leqslant l_3$.
We control $Z^{l_2}  (\eps \underline{\mathtt{V}}^{\eps,\flat}_{d})$ in $L^\infty$ norm by 
$ || \underline{\mathtt{V}}^\eps_{R,d} ||^*_{\eps,T}$. 
\end{enumerate}
\item if $l_1 \geqslant 2$, we apply twice  the Moser inequality of Lemma $\ref{Moser2}$ with  
\begin{eqnarray*} 
\alpha := Z^2 S(\underline{u}^\eps ) ,\ \beta := \underline{\mathtt{V}}^{\eps}_{d} ,\ \gamma := \eps \partial_d \varphi ,\ \underline{m} := m-2. 
\end{eqnarray*}
\end{enumerate}
\end{proof}
We now give the proof of $a)$. 
\begin{proof}[Proof of $a)$]
To explain in a clearer way our method, we will deal with the commutator
\begin{eqnarray}
\label{jj}
[\underline{\mathtt{v}}^\eps_d \partial_d , (\eps \partial_d )^k Z^l ] \varphi 
\end{eqnarray}
instead of 
\begin{eqnarray*}
[S(\underline{u}^\eps).\underline{\mathtt{v}}^\eps_d \partial_d , (\eps \partial_d )^k Z^l ] \varphi .
\end{eqnarray*}
This avoids heavy notations and does not change the mathematical analysis.
The commutator $(\ref{jj})$ is a linear combination of some
terms of the form 
\begin{eqnarray}
\label{tf1} &(\eps \partial_d )^k Z^{l_1} \underline{\mathtt{v}}^\eps_d . Z^{l_2} \partial_d  \varphi , 
 \quad &\mathrm{with} \ \mathrm{l_1 \neq 0} \ , l_1 + l_2 = l,
\\ \label{tf2} &(\eps \partial_d )^{k_1} \underline{\mathtt{v}}^\eps_d .  (\eps \partial_d )^{k_2} Z^{l} \partial_d  \varphi ,
  \quad &\mathrm{with} \ \mathrm{k_1 \neq 0} \ ,  k_1 + k_2 = k.
\end{eqnarray}
We look at the terms of the form $(\ref{tf1})$. 
Because $\underline{\mathtt{v}}^\eps$ is of the form $\underline{\mathtt{v}}^\eps =  \mathtt{v}^0  + \eps
\underline{\mathtt{V}}^\eps  $, it suffices to deal with the terms 
\begin{eqnarray}
\label{rf1} & (\eps \partial_d )^k Z^{l_1} \mathtt{v}^0_d . Z^{l_2} \partial_d  \varphi , \quad  &\mathrm{with} \ \mathrm{l_1 \neq 0} \ , l_1 + l_2 = l,
\\ \label{rf2} &(\eps \partial_d )^k Z^{l_1}  \underline{\mathtt{V}}^\eps_d .  Z^{l_2}  (\eps \partial_d ) \varphi 
 \quad  &\mathrm{with} \ \mathrm{l_1 \neq 0} \ , l_1 + l_2 = l .
\end{eqnarray}
We begin with the terms of the form $(\ref{rf1})$. 
Using that $ \mathtt{v}^0_d = x_d \mathtt{v}^{0,\flat}_d$ and commuting $x_d$ with $Z^{l_2}$, 
we have to deal with the terms of the form
\begin{eqnarray*}
(\eps \partial_d )^k Z^{l_1} \mathtt{v}^{0,\flat}_d . Z^{l'_2} (\eps \partial) \varphi 
\quad \mathrm{with} \ l_2' \leqslant l_2 .
\end{eqnarray*}
Using $L^\infty$ estimate of $\mathtt{v}^0$ yields the result.
We now have a look for the terms of the form $(\ref{rf2})$. 
We distinguish three cases:

\begin{enumerate}[$\bullet$]

\item if $k \neq 0$ then we write $k=k'+1$ and we commute $\eps \partial_d$ with $Z^{l_1}$. 
This yields some terms of the form 
\begin{eqnarray*}
(\eps \partial_d )^{k'} Z^{l'}  (\eps \partial_d \underline{\mathtt{V}}^\eps_d ) Z^{l_2} \eps \partial_d \varphi 
\end{eqnarray*}
with $l' \leqslant l_1$. 
We use that $\underline{\mathtt{V}}^\eps_d = \underline{\mathtt{V}}^\eps_{a,d} + \eps^M \, \underline{\mathtt{V}}^\eps_{R,d}$ 
and study the terms 
\begin{eqnarray}
\label{Ai}
\lambda^{m-2k-l} |   (\eps \partial_d )^{k'} Z^{l'}  (\eps \partial_d \underline{\mathtt{V}}^\eps_{a,d} ) Z^{l_2} \eps \partial_d
\varphi  |_{0,\lambda,T} ,
\\  \label{Bi}
\lambda^{m-2k-l} | \eps^M \, (\eps \partial_d )^{k'} Z^{l'}  (\eps \partial_d \underline{\mathtt{V}}^\eps_{R,d} ) Z^{l_2} \eps
\partial_d \varphi  |_{0,\lambda,T} .
\end{eqnarray}
Using the $L^\infty$ estimates $(\ref{linf})$, we get
$ (\ref{Ai}) \leqslant C  | \varphi |^\mathbf{E}_{\eps,m,\lambda,T} $.
We apply Lemma $\ref{Moser2}$ $(2)$ to 
\begin{eqnarray*}
\alpha := \eps \partial_d \mathtt{V}^\eps_{d,R} , \ \beta := \eps \partial_d  \varphi , \ \underline{m} =m-2
\end{eqnarray*}
 to bound $(\ref{Bi})$ by 
\begin{eqnarray}
\label{AM}
C \eps^M (|\varphi|^\mathbf{E}_{\eps,m,\lambda,T} + ||\varphi ||_{\eps,Lip,T} .|\underline{U}^\eps_R |^\mathbf{E}_{\eps,m,\lambda,T}). 
\end{eqnarray}
\item if $k=0$ and $l_1 \geqslant 2$ then with $l_1 := l_1' + 2$, the term $(\ref{rf2})$ can be rewritten 
\begin{eqnarray*}
(\eps \partial_d )^{k} Z^{l'_1  }  (Z^{2} \underline{\mathtt{V}}^\eps_d ). Z^{l_2} (\eps  \partial_d ) \varphi .
\end{eqnarray*}
We use that $\underline{\mathtt{V}}^\eps_d = \mathtt{V}^\eps_{a,d} + \eps^M \underline{\mathtt{V}}^\eps_{R,d}$ and study
the terms 
\begin{eqnarray}
\label{Aj}
\lambda^{m-2k-l} |(\eps \partial_d )^{k} Z^{l'_1}  (Z^{2} \mathtt{V}^\eps_{a,d} ). Z^{l_2} (\eps  \partial_d ) \varphi |_{0,\lambda,T},
\\  \label{Bj}
\lambda^{m-2k-l} |\eps^M  (\eps \partial_d )^{k} Z^{l'_1}  (Z^{2} \underline{\mathtt{V}}^\eps_{R,d} ). Z^{l_2} (\eps  \partial_d )
\varphi |_{0,\lambda,T}.
\end{eqnarray}
We use the $L^\infty$ estimates $(\ref{linf})$ to bound the term  $(\ref{Aj})$. 
We apply Lemma $\ref{Moser2}$ with 
\begin{eqnarray*}
\alpha := Z^2 \underline{V}^\eps_{R,d} ,\ \beta := \eps \partial_d \varphi ,\ \underline{m} = m-2
\end{eqnarray*}
 to bound $(\ref{Bj})$ by $(\ref{AM})$.
 %                 \\  
\item if $k=0$ and $l_1 =1$, using that $\underline{\mathtt{V}}^\eps_d = x_d \underline{\mathtt{V}}^{\eps,\flat}_d $ 
and commuting $x_d$ with
$Z^{l_2}$ yield a sum of terms of the form 
\begin{eqnarray*}
Z^{l_1} \eps \underline{\mathtt{V}}^{\eps,\flat}_d .Z^{l_2} (x_d  \partial_d \varphi )  .
\end{eqnarray*}
We bound it by $|| \eps \underline{\mathtt{V}}^{\eps,\flat}_d  ||_{1,T} . |\varphi |_{m,\lambda,T}$ and conclude with $|| \eps
\underline{\mathtt{V}}^{\eps,\flat}_d  ||_{1,T} \leqslant || \underline{\mathtt{V}}^{\eps}_d ||^*_{\eps,T} $. 
We can deal with the terms $(\ref{tf2})$ with the same methods.
\end{enumerate}
 For the proof of $b)$, proceed in the same way and substitute Lemma $\ref{Moser2}$ $(2)$ by Lemma $\ref{Moser2}$ $(1)$.
\end{proof}
\subsubsection{Source term estimates}

We now give two lemmas which will be useful when estimating source terms.

\begin{lem}
\label{poire}
Let $m$ be an integer. 
There is an increasing function $C:\R_+ \rightarrow \R_+^*$ such that  for all $\eps \in ]0,1]$, 
if $\eps^M \, || \underline{U}^\eps_R ||_{0,T}
\leqslant  \mathfrak{R}$ then for all $\varphi \in C^\infty_0 ((0,T) \times \overline{\Omega} )$,
 for all $\lambda \geqslant \lambda_0$,

\begin{enumerate}[$a)$]

\item for any $m$, 
\begin{eqnarray*}
|  \underline{J}. \varphi  |_{m,\lambda,T}
\leqslant C(\mathfrak{R}) (|\varphi  |_{ m,\lambda,T} 
+ \eps^M ||\varphi  ||_{0,T} .|\underline{U}^\eps_R |_{m,\lambda,T}),
\end{eqnarray*}
\item for $m$ even,
\begin{eqnarray*}
|  \underline{J}. \varphi   |^\mathbf{E}_{\eps,m,\lambda,T}
\leqslant  C(\mathfrak{R}) (  |\varphi |^\mathbf{E}_{\eps,m,\lambda,T} 
+  \eps^M || \varphi ||_{0,T} . | \underline{U}^\eps_R |^\mathbf{E}_{\eps,m,\lambda,T} ).
\end{eqnarray*}
\end{enumerate}
where $\underline{J}$ stands for $\underline{J}^\alpha $, $\underline{J}^{\beta,1}$ or $\underline{J}^{\beta,2}$.
\end{lem}

The proof of Lemma $\ref{poire}$ is straightforward and mainly lies on  Lemma $\ref{Moser2}$. 
It is left to the reader.

\begin{lem}
\label{carotte}
Let $m$ be an integer. 
There is an increasing function $C:\R_+ \rightarrow \R_+^*$ such that  for all $\eps \in ]0,1]$, 
if $\eps^M \, || \underline{U}^\eps_R ||^*_{\eps,T}
\leqslant  \mathfrak{R}$ then  for all $\lambda \geqslant \lambda_0$,
\begin{enumerate}[$a)$]
\item for any $m$, 
\begin{eqnarray*}
| \underline{\mathtt{v}}^\eps_d \partial_d  U^\eps_R |_{m-1,\lambda,T}
\leqslant C(\mathfrak{R}) (|U^\eps_R  |^\mathcal{N}_{\eps, m,\lambda,T} 
+ \eps^M || U^\eps_R  ||^*_{\eps,T} .|\underline{\mathtt{v}}^\eps_R |^\mathcal{N}_{\eps,m,\lambda,T}),
\end{eqnarray*}
\item for $m$ even,
\begin{eqnarray*}
| \underline{\mathtt{v}}^\eps_d \partial_d  U^\eps_R |^\mathbf{E}_{\eps,m-1,\lambda,T}
\leqslant  C(\mathfrak{R}) (  |U^\eps_R |^\mathbf{E}_{\eps,m,\lambda,T} 
+  \eps^M || U^\eps_R ||_{\eps,Lip} . | \underline{\mathtt{v}}^\eps_R |^\mathbf{E}_{\eps,m,\lambda,T} ).
\end{eqnarray*}
\end{enumerate}
\end{lem}

\begin{proof}[Proof of $a)$]
We will proceed in four steps.
\begin{enumerate}[$(i)$]
\item 
Because $ \underline{\mathtt{v}}^\eps_d =  \underline{\mathtt{v}}^0_d + \eps \underline{\mathtt{V}}^\eps_d $ (cf. $(\ref{linform})$)
 and $\mathtt{v}^0_d = x_d \mathtt{v}^{0,\flat}_d$, we get  
\begin{eqnarray*}
|\underline{\mathtt{v}}^\eps_d .\partial_d U^\eps_R |_{m-1,\lambda,T}
\leqslant C | U^\eps_R |_{m,\lambda,T} 
+ |\underline{\mathtt{V}}^\eps_d .(\eps \partial_d ) U^\eps_R |_{m-1,\lambda,T} .
\end{eqnarray*}
\item 
In order to estimate $ |\underline{\mathtt{V}}^\eps_d .(\eps \partial_d ) U^\eps_R |_{m-1,\lambda,T}$, 
we will control the terms of the form
\begin{eqnarray}
\label{cordf}
\lambda^{m-1-l} |(  Z^{l_1} \underline{\mathtt{V}}^\eps_d ). 
 Z^{l_2} (\eps \partial_d ) U^\eps_R |_{0,\lambda,T} ,
\end{eqnarray}
where $l_1, l_2 \in \N, \ l_1 + l_2 =l \leqslant m-1 $.
Because $\underline{V}^\eps = V^\eps_a + \eps^M \, \underline{V}^\eps_R$, 
it suffices to deal with the terms of the form 
\begin{eqnarray}
\label{cordfa}
\lambda^{m-1-l} |{ Z^{l_1} \mathtt{V}^\eps_{a,d} }. 
 Z^{l_2} (\eps \partial_d ) U^\eps_R |_{0,\lambda,T} ,
\\ \label{cordfb}
\lambda^{m-1-l} |{  Z^{l_1} \underline{\mathtt{V}}^\eps_{R,d} }. 
Z^{l_2} (\eps \partial_d ) U^\eps_R |_{0,\lambda,T} .
\end{eqnarray}
\item 
We begin with the term $(\ref{cordfa})$. We distinguish two cases:
\\
\begin{enumerate}[$\bullet$]
\item If $l_1 > 0$, then  $l_2 \leqslant m$ we use the $L^\infty$ estimates $(\ref{linf})$
 and bound the term $(\ref{cordfa})$ by 
\begin{eqnarray*}
 c.| (\eps \partial_d ) U^\eps_R |_{m-2,\lambda,T}  \leqslant  c.| U^\eps_R |_{m,\lambda,T} .
\end{eqnarray*}
\item If $l_1 =0$, then we write $\mathtt{V}^\eps_{a,d}  = x_d . \underline{\mathtt{V}}^{\eps,\flat}_{a,d} $ 
and commute $x_d$ with $ Z^{l_2}$ to find that $(\ref{cordfa})$ is a sum of terms of the form
\begin{eqnarray}
\label{cordfa3}
\lambda^{m-1-l} |\eps \underline{\mathtt{V}}^{\eps,\flat}_{a,d}  Z^{l'_2} U^\eps_R |_{0,\lambda,T} .
\end{eqnarray}
with $l'_2 \leqslant l_2 + 1$. 
Thus 
\begin{eqnarray*}
(\ref{cordfa3}) \leqslant C || \eps \underline{\mathtt{V}}^{\eps,\flat}_{a,d} ||_\infty . | U^\eps_R |_{m,\lambda,T }
\leqslant C  || \mathtt{V}^{\eps}_{a,d} ||_{\eps,T}^*  . | U^\eps_R |_{m,\lambda,T } . 
 \end{eqnarray*}
    \end{enumerate}
\item 

We now look at the term $(\ref{cordfb})$. We distinguish three cases:
\\
\begin{enumerate}[$\bullet$]
\item If $l_1 \neq 0$, we apply Lemma $\ref{Moser2}$ $(1)$ to 
\begin{eqnarray*}
\alpha := Z \underline{\mathtt{V}}^\eps_{R,d} , \
 \beta := (\eps \partial_d ) U^\eps_R , \ \underline{m} := m-2. 
 \end{eqnarray*}
The term $(\ref{cordfb})$ is bounded by 
\begin{eqnarray*}
\label{cordf1}
 c. \{ || \underline{\mathtt{V}}^\eps_{R,d} ||_{1,T} .| U^\eps_R |^\mathcal{N}_{\eps,m,\lambda,T} 
+ | \underline{\mathtt{V}}^\eps_{R,d} |_{m,\lambda,T} . || U^\eps_R ||_{\eps,Lip}  \}.
\end{eqnarray*}
\item If $l_1 =0$, we write $\underline{\mathtt{V}}^\eps_{R,d}  = x_d . \underline{\mathtt{V}}^{\eps,\flat}_{R,d} $ 
and commute $x_d$ with $ Z^{l_2}$ to find that $(\ref{cordf})$ is a sum of terms of the form
\begin{eqnarray}
\label{cordf3}
\lambda^{m-1-l} |\eps \underline{\mathtt{V}}^{\eps,\flat}_{R,d} Z^{l'_2} U^\eps_R |_{0,\lambda,T} .
\end{eqnarray}
with  $l'_2 \leqslant l_2 + 1$. 
Thus the term $(\ref{cordf3})$ is bounded by 
$ || \eps \underline{\mathtt{V}}^{\eps,\flat}_{R,d} || _\infty . |  U^\eps_R |_{ m,\lambda,T }$.   
    \end{enumerate}
\end{enumerate}
\end{proof}
The proof of $a)$ follows the same path and is slightly more complicated.
\begin{proof}[Proof of $b)$]
We will proceed in four steps.
\begin{enumerate}[$(i)$]
\item 
Because $ \underline{\mathtt{v}}^\eps_d =  \underline{\mathtt{v}}^0_d + \eps \underline{\mathtt{V}}^\eps_d $ (cf. $(\ref{linform})$)
 and $\mathtt{v}^0_d = x_d \mathtt{v}^{0,\flat}_d$, we get 
\begin{eqnarray*}
|\underline{\mathtt{v}}^\eps_d .\partial_d U^\eps_R |^\mathbf{E}_{m-1,\lambda,T}
\leqslant C | U^\eps_R |^\mathbf{E}_{m,\lambda,T} 
+ |\underline{\mathtt{V}}^\eps_d .(\eps \partial_d ) U^\eps_R |^\mathbf{E}_{m-1,\lambda,T} .
\end{eqnarray*}
\item 
In order to estimate $ |\underline{\mathtt{V}}^\eps_d .(\eps \partial_d ) U^\eps_R |^\mathbf{E}_{\eps,m-1,\lambda,T}$, 
we will control the terms of the form
\begin{eqnarray}
\label{rdf}
\lambda^{m-1-2k-l}  \, |{ (\eps \partial_d )^{k_1} Z^{l_1} \underline{\mathtt{V}}^\eps_d }. 
(\eps \partial_d )^{k_2} Z^{l_2} (\eps \partial_d ) U^\eps_R |_{0,\lambda,T} ,
\end{eqnarray}
where 
\begin{eqnarray*}
k_i, l_i \in \N, \ k_1 + k_2 = k, \ l_1 + l_2 =l, \ 2k+l \leqslant m-1 .
\end{eqnarray*}
Because $\underline{V}^\eps = V^\eps_a + \eps^M \, \underline{V}^\eps_R$, 
it suffices to deal with the terms of the form 
\begin{eqnarray}
\label{rdfa}
\lambda^{m-1-2k-l}  \, |{ (\eps \partial_d )^{k_1} Z^{l_1} \mathtt{V}^\eps_{a,d} }. 
(\eps \partial_d )^{k_2} Z^{l_2} (\eps \partial_d ) U^\eps_R |_{0,\lambda,T} ,
\\ \label{rdfb}
\lambda^{m-1-2k-l}  \, |{ (\eps \partial_d )^{k_1} Z^{l_1} \underline{\mathtt{V}}^\eps_{R,d} }. 
(\eps \partial_d )^{k_2} Z^{l_2} (\eps \partial_d ) U^\eps_R |_{0,\lambda,T} .
\end{eqnarray}
\item 
We begin with the term $(\ref{rdfa})$. We distinguish two cases:
\\
\begin{enumerate}[$\bullet$]
\item If $l_1 + k_1 > 0$, then  $2k_2 + l_2 \leqslant m$ we use the $L^\infty$ estimates $(\ref{linf})$
 and bound the term $(\ref{rdfa})$ by 
\begin{eqnarray*}
 c.| (\eps \partial_d ) U^\eps_R |^\mathbf{E}_{m-2,\lambda,T}  \leqslant  c.| U^\eps_R |^\mathbf{E}_{m,\lambda,T} .
\end{eqnarray*}
\item If $l_1 = k_1=0$, then we write $\mathtt{V}^\eps_{a,d}  = x_d . \underline{\mathtt{V}}^{\eps,\flat}_{a,d} $ 
and commute $x_d$ with $(\eps \partial_d )^{k_2} Z^{l_2}$ to find that $(\ref{rdfa})$ is a sum of terms of the form
\begin{eqnarray}
\label{rdfa3}
\lambda^{m-1-2k-l} \,  |\eps \underline{\mathtt{V}}^{\eps,\flat}_{a,d} (\eps \partial_d )^{k'_2} Z^{l'_2} U^\eps_R |_{0,\lambda,T} .
\end{eqnarray}
with $k'_2 \leqslant k_2$, $l'_2 \leqslant l_2 + 1$. 
Thus 
\begin{eqnarray*}
(\ref{rdfa3}) \leqslant C || \eps \underline{\mathtt{V}}^{\eps,\flat}_{a,d} ||_\infty . |  U^\eps_R |^\mathbf{E}_{m,\lambda,T }
\leqslant C  || \mathtt{V}^{\eps}_{a,d} ||_{\eps,T}^*  . | U^\eps_R |^\mathbf{E}_{m,\lambda,T } . 
 \end{eqnarray*}
    \end{enumerate}
\item 
We now look at the term $(\ref{rdfb})$. We distinguish three cases:
\\
\begin{enumerate}[$\bullet$]

\item If $l_1 \neq 0$, we apply Lemma $\ref{Moser2}$ $(2)$ to 
\begin{eqnarray*}
\alpha := Z \underline{\mathtt{V}}^\eps_{R,d} , \
 \beta := (\eps \partial_d ) U^\eps_R , \ \underline{m} := m-2. 
 \end{eqnarray*}
The term $(\ref{rdfb})$ is bounded by 
\begin{eqnarray*}
\label{rdf1}
 c. \{ || \underline{\mathtt{V}}^\eps_{R,d} ||_{1,T} .| U^\eps_R |^\mathbf{E}_{\eps,m,\lambda,T} 
+ | \underline{\mathtt{V}}^\eps_{R,d} |^\mathbf{E}_{\eps,m,\lambda,T} . || U^\eps_R ||_{\eps,Lip}  \}.
\end{eqnarray*}
\item If $l_1 =0$ and $ k_1 \neq 0$, we apply Lemma $\ref{Moser2}$ $(2)$ to 
\begin{eqnarray*}
\alpha := \eps \partial_d \underline{\mathtt{V}}^\eps_{R,d} , \  \beta := \eps \partial_d U^\eps_R , \  \underline{m} := m-2. 
\end{eqnarray*}
The term $(\ref{rdfb})$ is bounded by 
\begin{eqnarray*}
\label{rdf2}
\lambda^{-1} c. \{ || \underline{\mathtt{V}}^\eps_{R,d} ||_{\eps,Lip,T} .| U^\eps_R |^\mathbf{E}_{\eps,m,\lambda,T} 
+ | \underline{\mathtt{V}}^\eps_{R,d} |^\mathbf{E}_{\eps,m,\lambda,T} . || U^\eps_R ||_{\eps,Lip} \}.
\end{eqnarray*}
\item If $l_1 = k_1=0$, we write $\underline{\mathtt{V}}^\eps_{R,d}  = x_d . \underline{\mathtt{V}}^{\eps,\flat}_{R,d} $ 
and commute $x_d$ with $(\eps \partial_d )^{k_2} Z^{l_2}$ to find that $(\ref{rdf})$ is a sum of terms of the form
\begin{eqnarray}
\label{rdf3}
\lambda^{m-1-2k-l} |\eps \underline{\mathtt{V}}^{\eps,\flat}_{R,d} (\eps \partial_d )^{k'_2} Z^{l'_2} U^\eps_R |_{0,\lambda,T} .
\end{eqnarray}
with $k'_2 \leqslant k_2$, $l'_2 \leqslant l_2 + 1$. 
Thus  the term $(\ref{rdf3})$ is bounded by 
$ || \eps \underline{\mathtt{V}}^{\eps,\flat}_{R,d} || _\infty . |  U^\eps_R |^\mathbf{E}_{m,\lambda,T }$.   
    \end{enumerate}
\end{enumerate}
\end{proof}
We now attack the conormal estimates of $U^\eps_R$.
\subsection{Conormal estimates}
\label{conorm}
We look for estimates of $\lambda^{m-l} | Z^l U^\eps_R |_{0,\lambda,T_0}$.
As for $L^2$ estimates, we will proceed in three steps estimating 
\begin{eqnarray*}
\lambda^{m-l} | Z^l W^\eps_R |_{0,\lambda,T_0 }, \ \lambda^{m-l} | Z^l
(\frac{x_d}{\eps} W^\eps_R )|_{0,\lambda,T_0 } \ \mathrm{and} \ \lambda^{m-l} | Z^l V^\eps_R |_{0,\lambda,T_0 } .
\end{eqnarray*}
\subsubsection{Estimate of $W^\eps_R$}
\label{coW}
We apply the derivative $Z^l$ to the equation $\X_{\underline{v}^{\eps}} W^\eps_R = f^\eps_I $.
We get $ \X_{\underline{v}^{\eps}} (Z^l W^\eps_R ) = \tilde{f}^\eps_I $, where $\tilde{f}^\eps_I  := [\X_{\underline{\mathtt{v}}^{\eps}}
,Z^l ] W^\eps_R + Z^l f^\eps_I $.
 Then we estimate $\tilde{f}^\eps_I $ thanks to Lemma $\ref{Moser2}$ 
  and use a classic $L^2$ estimate.
 \subsubsection{Estimate of $\frac{x_d}{\eps}  W^\eps_R$}
\label{coXW}
 Because the control of the commutator $[ \underline{\mathtt{v}}^{\eps,\flat}_d \frac{x_d}{\eps} ,Z^l ] W^\eps_R $ seems difficult,
 we write 
 $\underline{\mathtt{v}}^{\eps,\flat}_d  \frac{x_d}{\eps} W^\eps_R  =
  \underline{\mathtt{v}}^{0}_d \frac{x_d}{\eps} W^\eps_R + \underline{\mathtt{V}}_d^\eps. W^\eps_R  $. 
 Thus   $\frac{x_d}{\eps}  W^\eps_R$ verifies the equation
 \begin{eqnarray*}
\label{coco}
 \X_{\underline{v}^{\eps}} (\frac{x_d}{\eps}  W^\eps_R ) = 
 \underline{\mathtt{v}}^{0,\flat}_d \frac{x_d}{\eps}  W^\eps_R +  \underline{\mathtt{V}}_d^\eps .
 W^\eps_R + f^{\eps}_{II} .
\end{eqnarray*}
We apply the derivative $Z^l$ to this equation and find for $Z^l (\frac{x_d}{\eps}  W^\eps_R )$ the equation
\begin{eqnarray*}
( \X_{\underline{\mathtt{v}}^{\eps}} -  \underline{\mathtt{v}}^{0,\flat}_d ) Z^l (\frac{x_d}{\eps}  W^\eps_R ) = \tilde{f}^{\eps}_{II} ,
\end{eqnarray*}
where 
\begin{eqnarray*}
\tilde{f}^{\eps}_{II} := [ \X_{\underline{\mathtt{v}}^{\eps}} -  \underline{\mathtt{v}}^{0,\flat}_d  , Z^l ] \frac{x_d}{\eps}  W^\eps_R 
+ Z^l ( \underline{\mathtt{V}}_d^\eps . W^\eps_R ) +  Z^l f^{\eps}_{II} .
\end{eqnarray*}
We estimate $[ \X_{\underline{\mathtt{v}}^{\eps}}  -  \underline{\mathtt{v}}^{0,\flat}_d ,  Z^l  ] (\frac{x_d}{\eps}  W^\eps_R )$ 
thanks to Lemma $\ref{Moser2}$ 
and $Z^l (  \underline{\mathtt{V}}_d^\eps . W^\eps_R )$ with Lemma $\ref{Moser2}$ and the previous subsection. 
We finally use a classic $L^2$ estimate.
\subsubsection{Estimate of $V^\eps_R$}
\label{coV}
We apply the derivative $Z^l$ to the boundary value problem $(\ref{Sx})$ and find that $Z^l V^\eps_R$ 
is solution of the following boundary value problem:
\begin{eqnarray*}
\left\{
\begin{array}{ccc}
&( \mathfrak{L}^\star (\underline{u}^\eps ,\partial)  + 
\underline{J}^\alpha ).(Z^l V^\eps_R )=
\tilde{f}^\eps_{III} \quad &\mathrm{when} \ x_d > 0 ,
\\ &Z^l V^\eps_{R,d} = 0  \quad &\mathrm{when} \ x_d = 0 ,
\end{array}
\right.
\end{eqnarray*}
 where 
 \begin{eqnarray*}
  \tilde{f}^\eps_{III} := 
 [S^\star (\underline{u}^\eps ) \X_{\underline{\mathtt{v}}^\eps}  + \L^\star (\partial_x ) + \underline{J}^\alpha , Z^l ] V^\eps_R +   Z^l
f^\eps_{III}.
\end{eqnarray*}
We estimate the commutator $[S^\star (\underline{u}^\eps ) \X_{\underline{\mathtt{v}}^\eps}  , Z^l ] V^\eps_R $ thanks to Proposition
$\ref{lapin}$ and the commutator $[ \underline{J}^\alpha , Z^l ] V^\eps_R $ thanks to lemma  $\ref{Moser2}$. 
We estimate $Z^l f^\eps_{III}$ thanks to the
subsection $\ref{coW}$ and $\ref{coXW}$. 
We now sketch the method to estimate $[\L^\star (\partial_x ), Z^l ] V^\eps_R $.
First notice that $ [\L^\star (\partial_x ), Z^l ] V^\eps_R = [\L^\star_d \partial_d , Z^l ] V^\eps_R$ is a sum of terms of the form $Z^{l'}
\L^\star_d \partial_d V^\eps_R $  with $l' \leqslant l$.
Then we extirpate $\L^\star_d \partial_d V^\eps_R $  by the first equation of $(\ref{Sx})$.
We can apply the previous methods for the resulting terms.
We use Lemma $\ref{carotte}$ to handle the most delicate point i.e. the estimate of $| \underline{\mathtt{v}}_d^\eps . \partial_d U^\eps_R |_{m-1,\lambda,T_0}$.
This supplies a good estimate of $[\L^\star (\partial_x ), Z^l ] V^\eps_R $.
We conclude as in subsection $\ref{estiV}$ with a $L^2$ estimate.
\subsection{Normal estimates}
\label{norm}

We look for estimates of $\lambda^{m-2k-l} | (\eps \partial)^k Z^l U^\eps_R |_{0,\lambda,T}$, with $k>0$ and $2k+l \leqslant m$. 
We begin to estimate for the noncharacteristic part: $\L_d^\star  V^\eps_{R}$.
This corresponds to the pressure $p^\eps_R$ and $\mathtt{v}^\eps_{R,d}$. 
Then we will estimate   $| W^\eps_{R} |^\mathbf{E}_{\eps,m,\lambda,T_0}$,   $\frac{x_d}{\eps}  | W^\eps_{R} |^\mathbf{E}_{\eps,m,\lambda,T_0}$ and finally
 $| \mathtt{v}^\eps_{R,y} |^\mathbf{E}_{\eps,m,\lambda,T_0}$ where we denote by $\mathtt{v}_y$ the tangential velocity.
\subsubsection{Normal estimates of $L_d^\star V^\eps_R$.}
\label{norm1}
We notice that 
 \begin{eqnarray*}
| \L_d^\star \eps \partial_d V^\eps_R |^\mathbf{E}_{\eps,m-2,\lambda,T} \leqslant 
\lambda^{-1} | \L_d^\star \eps \partial_d V^\eps_R |^\mathbf{E}_{\eps,m-1,\lambda,T}
\end{eqnarray*}
 and extirpate $\L_d^\star \eps \partial_d V^\eps_R$  from the equation $(\ref{AL})$. 
Among the resulting terms, the most delicate one to estimate is $| \underline{\mathtt{v}}^\eps_d . \partial_d 
U^\eps_R|^\mathbf{E}_{\eps,m-1,\lambda,T_0} $ which is given by Lemma  $\ref{carotte}$.
Notice that there is no problem with the term $\underline{J}^{\beta,2} \frac{x_d}{\eps} W^\eps_R $.

\subsubsection{Normal estimates of $W^\eps_R$ .}
\label{norm2}

 We apply the derivative $( \eps \partial_d )^k Z^l $ to the equation $\X_{\underline{\mathtt{v}}^{\eps}} W^\eps_R = f^\eps_I $.
 We find that $( \eps \partial_d )^k Z^l  W^\eps_R $ verify the following equation: 
 \begin{eqnarray*}
 \X_{\underline{\mathtt{v}}^{\eps}}  (\eps \partial_d
 )^k Z^l W^\eps_R = \hat{f}^\eps_I 
  \end{eqnarray*}
  where 
  \begin{eqnarray*}
 \hat{f}^\eps_I := [ \X_{\underline{\mathtt{v}}^{\eps}} , ( \eps \partial_d )^k Z^l ] W^\eps_R + ( \eps \partial_d )^k Z^l  f^\eps_I .
  \end{eqnarray*}
 We estimate the commutator $[ \X_{\underline{\mathtt{v}}^{\eps}} , ( \eps \partial_d )^k Z^l ] W^\eps_R $ thanks to Proposition $a)$ and
 $c)$ and use a $L^2$ estimate.
 
 \subsubsection{Normal estimates of $ \frac{x_d}{\eps} W^\eps_R$ .}
\label{norm3}

 We apply the derivative $( \eps \partial_d )^k Z^l $ to the equation $(\ref{coco})$ to find for $\hat{W}^\eps_R := ( \eps
 \partial_d )^k Z^l ( \frac{x_d}{\eps} W^\eps_R )$ the equation $(\X_{\underline{\mathtt{v}}^\eps} - \mathtt{v}^{0,\flat}_d )
 \hat{W}^\eps_R = \hat{f}^\eps_{III}$ where 
 \begin{eqnarray*}
 \hat{f}^\eps_{III} := [ \X_{\underline{\mathtt{v}}^\eps} - \mathtt{v}^{0,\flat}_d , ( \eps
 \partial_d )^k Z^l ] \hat{W}^\eps_R + ( \eps \partial_d )^k Z^l ( \underline{\mathtt{V}}^\eps_d . W^\eps_R )
 +  ( \eps \partial_d )^k Z^l f^\eps_{III} .
 \end{eqnarray*}
 We estimate the commutator $[ X_{\underline{\mathtt{v}}^\eps} , ( \eps \partial_d )^k Z^l ] \hat{W}^\eps_R$ thanks to 
 Proposition $\ref{lapin}$, the second term of the right member thanks to Lemma $\ref{Moser2}$ and the previous subsection. 
 We estimate $\hat{f}^\eps_{III}$ using that because $\beta_2$ is rapidly decreasing with respect to its last argument, 
 all the derivatives of
 $\frac{x_d}{\eps} \beta_2 (\eps,t,x,\frac{x_d}{\eps})$ are bounded in $L^\infty$ uniformly with respect to $\eps$ 
 and $(f^\eps_{III} ) \in \mathbf{E}^m$.
 We end with a $L^2$ estimate.
 
  \subsubsection{Normal estimates of $\mathtt{v}^\eps_{R,y} $}
 \label{norm4}

 Looking at the first $d-1$ equations of $(\ref{AL})$, we find for $\mathtt{v}^\eps_{R,y}$ an equation of the form 
 \begin{eqnarray*}
 (\underline{\rho}^\eps .  \X_{\underline{\mathtt{v}}^\eps}  + \underline{J}_1^\alpha ) . \mathtt{v}^\eps_{R,y} =
 f^\eps_{IV}
  \end{eqnarray*}
 where 
 \begin{eqnarray*}
 f^\eps_{IV} := - \underline{J}_2^\alpha . \mathtt{v}^\eps_{R,d} - \nabla_y p^\eps 
 - \underline{J}^{\beta,1} \, W^\eps_R - \underline{J}^{\beta,2} \, \frac{x_d}{\eps}) W^\eps_R  
 - R^\eps 
  \end{eqnarray*}
 and $\underline{J}_1^\alpha $, $\underline{J}_2^\alpha$ are some matrices extracted from $\underline{J}^\alpha$.
 We apply the derivative $(\eps \partial_d )^k Z^l$ to this equation and find for 
 $\hat{\mathtt{v}}^\eps_{R,y} := (\eps \partial_d )^k Z^l \mathtt{v}^\eps_{R,y}$ the equation 
  \begin{eqnarray*}
 ( \underline{\rho}^\eps  X_{\underline{\mathtt{v}}^\eps} + \underline{J}_1^\alpha ). \hat{\mathtt{v}}^\eps_{R,y} =
 \hat{f}^\eps_{IV} 
 \end{eqnarray*}
  where 
  \begin{eqnarray*}
 \hat{f}^\eps_{IV} := [ \underline{\rho}^\eps  \X_{\underline{\mathtt{v}}^\eps} +
 \underline{J}_1^\alpha  , (\eps \partial_d )^k  Z^l ] v^\eps_{R,y} + (\eps \partial_d )^k  Z^l f^\eps_{IV} .
  \end{eqnarray*}
 We estimate the commutator $[\underline{\rho}^\eps  \X_{\underline{\mathtt{v}}^\eps} , (\eps \partial_d )^k  Z^l ] v^\eps_{R,y}$
 thanks to Proposition $\ref{lapin}$, the commutator $[\underline{J}_1^\alpha , (\eps \partial_d )^k  Z^l  ]  v^\eps_{R,y} $
 thanks to Lemma $\ref{Moser2}$. 
 When looking at the term $(\eps \partial_d )^k Z^l f^\eps_{IV}$, the contributions of $- \underline{J}_2^\alpha
 .\underline{\mathtt{v}}^\eps_{R,d}$ and of $\underline{R}^\eps$ are easy to estimate.
 We estimate the contribution of $\underline{J}^{\beta,1}  \, W^\eps_R$ and  $\underline{J}^{\beta,2} \,  \frac{x_d}{\eps} W^\eps_R  $ 
  thanks to the previous subsections
 ($\ref{norm2}$ and  $\ref{norm3}$) and to Lemma $\ref{Moser2}$. 
 It remains to explain how to estimate 
  \begin{eqnarray*}
 \lambda^{m-2k-l}  \, |(\eps \partial_d )^k Z^l \nabla_y p^\eps |_{0,\lambda,T_0}. 
  \end{eqnarray*}
 Commuting $\eps \partial_d $ with $Z^l \nabla_y$, we are lead to estimate the terms 
 \begin{eqnarray*}
 \lambda^{m-2k-l}  \, | (\eps \partial_d )^{k'} Z^{l'} \nabla_y p^\eps |_{0,\lambda,T_0},
 \end{eqnarray*}
  with $k' \leqslant k-1$, $l' \leqslant
 l-1$. All these terms are bounded by $|\eps \partial_d p^\eps |^\mathbf{E}_{\eps,m-1,\lambda,T_0}$ which is estimated in subsection
 $\ref{norm2}$. 
 
 \subsection{An iterative scheme}
\label{Itsch}
 
 We define an iterative scheme $(U_R^{\eps,\nu} )_{\eps,\nu }$ by 
 \begin{eqnarray*}
 ( \mathfrak{L}^\star (u^{\eps,\nu} ,\partial) + J^{\alpha,\nu} ) V_R^{\eps,\nu+1} 
 + W^\eps_R  ( \underline{J}^{\beta,1,\nu} + \frac{x_d}{\eps}  .\underline{J}^{\beta,2,\nu}) = R^\eps_v 
  \quad \mathrm{when} \ x_d >0,
 \\ \X_{\mathtt{v}^{\eps,\nu}} W_R^{\eps,\nu+1} + \eps V_R^{\eps,\nu+1} .\nabla W_a^{\eps} = - \eps R^\eps_W \quad \mathrm{when}  \ x_d >0,
 \\ \mathtt{V}_{R,d}^{\eps,\nu+1} =0 \quad \mathrm{when}  \ x_d =0 , 
 \end{eqnarray*}
 when $ J^{\alpha,\nu}$ denotes the $(d+1) \times (d+1)$ matrix 
 \begin{eqnarray*}
 J^{\alpha,\nu} := J^{\alpha} ( \eps, V_a^{\eps} , Z_\eps V_a^{\eps} , W_a^{\eps,\flat} , \eps^M \, U_R^{\eps,\nu} ) ,
  \end{eqnarray*}
 for $i \in \{1,2\}$, $J^{\beta,i,\nu}$ denotes the $\R^{d+1}$-valued functions 
  \begin{eqnarray*}
  J^{\beta,i,\nu} := J^{\beta,i} (\eps, V_a^{\eps} , Z_\eps V_a^{\eps} , W_a^{\eps,\flat}, \eps^M \, U_R^{\eps,\nu} ) 
   \end{eqnarray*}
   and
   \begin{eqnarray*}
 u^{\eps,\nu} := u^0 + \eps U^{\eps,\nu} , \quad U^{\eps,\nu} := U^{\eps}_a + \eps^M U_R^{\eps,\nu}.
  \end{eqnarray*}
   Because this step is now very classic in BKW theory (cf. \cite{Gu93}, \cite{Gu92}, \cite{Gu95}, \cite{moi}, \cite{moi2},
 \cite{moi3}...), we only sketch as a preview how to deduce from the Sobolev embedding lemma $\ref{sobo}$ and from the linear estimates of
Theorem $\ref{pok}$  the uniform boundedness of $ (U_R^{\eps,\nu})_{\eps,\nu} $. 
 In order to do so, we fix two strictly positive real $h$ and $R$ , $\lambda$ and $\eps_0$ such that 
  \begin{eqnarray}
  \label{R}
 \lambda^{-1} \, \lambda_0 \leqslant \mathrm{min} (\frac{h}{2},\frac{1}{2}) ,
  \\ \eps_0^{M-\frac{1}{2}} \, c T_0 e^{\lambda T_0 } \leqslant  \frac{1}{2h} .
  \end{eqnarray}
 \begin{prop}
  
 If 
 \begin{eqnarray*}
 \eps^M || U^{\eps,\nu}_R ||^*_{\eps,T_0} \leqslant \mathfrak{R} \ \mathrm{and}  \ || U^{\eps,\nu}_R ||^\mathbf{E}_{\eps,m,\lambda,T_0}
 \leqslant h ,
  \end{eqnarray*}
 then  
 \begin{eqnarray*}
 \eps^M || U^{\eps,\nu+1}_R ||^*_{\eps,T_0} \leqslant  \mathfrak{R} \ \mathrm{and}  \ || U^{\eps,\nu+1}_R ||^\mathbf{E}_{\eps,m,\lambda,T_0} \leqslant h .
   \end{eqnarray*}
 \end{prop}
 
  \begin{proof}
   We will proceed in three steps.
   
   \begin{enumerate}
 \item We begin applying the Sobolev embedding lemma $\ref{sobo}$:
  \begin{eqnarray*}
\sqrt{\eps}  || U^{\eps,\nu+1}_R  ||^*_{\eps,T_0} \leqslant c T_0 e^{\lambda T_0 } \,  | U^{\eps,\nu+1}_R  |^E_{\eps,m,\lambda,T_0} .
\end{eqnarray*}
  Thanks to $(\ref{R})$, we obtain
  \begin{eqnarray*}
  \eps^M \  || U^{\eps,\nu+1}_R  ||^*_{\eps,T_0} \leqslant \frac{1}{2h} | U^{\eps,\nu+1}_R  |^E_{\eps,m,\lambda,T_0} .
  \end{eqnarray*}
 \item 
  We apply Theorem $\ref{pok}$ and find
\begin{eqnarray*}
|U_R^{\eps,\nu+1} |^\mathbf{E}_{\eps,m,\lambda,T_0} \leqslant
 \lambda^{-1} \lambda_0 .( 1+ \eps^M || U_R^{\eps,\nu+1} ||^*_{\eps,T_0} . | U_R^{\eps,\nu} |^\mathbf{E}_{\eps,m,\lambda,T_0}  ).
\end{eqnarray*}
 Thanks to $(\ref{R})$, we obtain
  \begin{eqnarray*}
  |U_R^{\eps,\nu+1} |^\mathbf{E}_{\eps,m,\lambda,T_0} \leqslant h.
  \end{eqnarray*}
 \item 
  Thanks to the Sobolev embedding lemma $\ref{sobo}$,  we obtain $ \eps^M || U^{\eps,\nu+1}_R ||^*_{\eps,T_0} \leqslant  \mathfrak{R}$.
\end{enumerate}
  \end{proof}
  \textit{I thank Olivier Gu\`es for suggesting me this subject; in September $2001$.}
 
\nocite{*}
\bibliographystyle{plain}
\bibliography{a}

$2000$ MSC: 76N20 Boundary-layer theory, 35Q35 Other equations arising in fluid mechanics, 35L50 Boundary value problems for
hyperbolic systems of first-order PDE.
\\
\\ Franck SUEUR
\\ \address{Laboratoire d'Analyse, de Topologie et de Probabilit\'e
\\ Centre de Math\'ematiques et d'Informatique 
\\ 39, rue F. Joliot Curie 
\\ 13453 Marseille Cedex 13 }
\\  \email{fsueur@cmi.univ-mrs.fr}
\\ \urladdr{http://www.cmi.univ-mrs.fr/~fsueur/}

 \end{document}